\documentclass[preprint,12pt]{elsarticle}
\usepackage{amssymb}
\usepackage{amsmath}
\usepackage{amsthm}
\usepackage{mathscinet}
\usepackage[colorlinks]{hyperref}
\AtBeginDocument{
  \hypersetup{
    linkcolor=blue,
    citecolor=green,
  }
}
\usepackage{cleveref}
\usepackage{geometry}
\usepackage[T1]{fontenc}

\usepackage[displaymath,mathlines,pagewise]{lineno}
\newcommand*\patchAmsMathEnvironmentForLineno[1]{%
  \expandafter\let\csname old#1\expandafter\endcsname\csname #1\endcsname
  \expandafter\let\csname oldend#1\expandafter\endcsname\csname end#1\endcsname
  \renewenvironment{#1}%
     {\linenomath\csname old#1\endcsname}%
     {\csname oldend#1\endcsname\endlinenomath}}%
\newcommand*\patchBothAmsMathEnvironmentsForLineno[1]{%
  \patchAmsMathEnvironmentForLineno{#1}%
  \patchAmsMathEnvironmentForLineno{#1*}}%
\AtBeginDocument{%
\patchBothAmsMathEnvironmentsForLineno{equation}%
\patchBothAmsMathEnvironmentsForLineno{align}%
\patchBothAmsMathEnvironmentsForLineno{flalign}%
\patchBothAmsMathEnvironmentsForLineno{alignat}%
\patchBothAmsMathEnvironmentsForLineno{gather}%
\patchBothAmsMathEnvironmentsForLineno{multline}%
}

\crefname{equation}{}{}

\newtheorem{theorem}{Theorem}[section]
\newtheorem{lemma}[theorem]{Lemma}

\theoremstyle{definition}
\newtheorem{definition}[theorem]{Definition}

\newtheorem{notation}[theorem]{Notation}

\theoremstyle{remark}
\newtheorem{remark}[theorem]{Remark}

\numberwithin{equation}{section}

\journal{~~}

\begin{document}

\begin{frontmatter}

\title{Boundary pointwise regularity for fully nonlinear parabolic equations and an application to regularity of free boundaries\tnoteref{t1}}
\author[rvt]{Yuanyuan Lian}
\ead{lianyuanyuan@sjtu.edu.cn; lianyuanyuan.hthk@gmail.com}

\author[rvt]{Kai Zhang\corref{cor1}}
\ead{zhangkaizfz@gmail.com}
\address[rvt]{School of Mathematical Sciences, Shanghai Jiao Tong University, Shanghai, China}
\cortext[cor1]{Corresponding author. ORCID: \href{https://orcid.org/0000-0002-1896-3206}{0000-0002-1896-3206}}
\tnotetext[t1]{This research is supported by the China Postdoctoral Science Foundation (Grant No.
2021M692086 and 2022M712081), the National Natural Science Foundation of China (Grant No. 12031012 and 11831003) and  the Institute of Modern Analysis-A Frontier Research Center of Shanghai.}

\begin{abstract}
In this paper, we prove boundary pointwise $C^{k,\alpha}$ regularity for any $k\geq 1$ for fully nonlinear parabolic equations. As an application, we give a direct and short proof of the higher regularity of the free boundaries in obstacle-type problems.
\end{abstract}

\end{frontmatter}

\section{Introduction}\label{S1}
We study the boundary pointwise regularity for viscosity solutions of the fully nonlinear parabolic equations:
\begin{equation}\label{e.Dirichlet}
\left\{\begin{aligned}
&u_t-F(D^2u,x,t)=f&& ~~\mbox{ in}~~\Omega\cap Q_1;\\
&u=g&& ~~\mbox{on}~~\partial \Omega\cap Q_1,
\end{aligned}\right.
\end{equation}
where $\Omega \subset \mathbb{R}^{n+1}$ is a bounded domain and $Q_1$ is the unit cylinder (see \Cref{no1.1}). Throughout this paper, we always assume that $F$ is uniformly elliptic with positive constants $\lambda$ and $\Lambda$: for any $M,N\in \mathcal{S}^n$ and $(x,t)\in \Omega\cap Q_1$,
\begin{equation}\label{e1.0}
\mathcal{M}^-(M,\lambda,\Lambda)\leq F(M+N,x,t)-F(N,x,t)
\leq \mathcal{M}^+(M,\lambda,\Lambda),
\end{equation}
where $\mathcal{S}^n$ denotes the set of $n\times n$ symmetric matrices; the $\mathcal{M}^-$ and $\mathcal{M}^+$ are the extremal Pucci's operators, i.e.,
\begin{equation*}
\mathcal{M}^-(M,\lambda,\Lambda)=\lambda\sum_{\lambda_i>0} \lambda_i+
\Lambda\sum_{\lambda_i<0} \lambda_i,~
\mathcal{M}^+(M,\lambda,\Lambda)=\Lambda\sum_{\lambda_i>0} \lambda_i+
\lambda\sum_{\lambda_i<0} \lambda_i~
\end{equation*}
and $\lambda_i$ are eigenvalues of $M$. In this paper, a solution always means the $L^{n+1}$-viscosity solution and we refer to \cite{MR1789919} for details.

In this paper, we prove a series of boundary pointwise regularity (see main results in this section) under suitable pointwise conditions. In particular, we show that $u$ has higher regularity if the derivatives of $u$ vanish at boundary, which was first found in the investigation of boundary regularity for elliptic equations (see \cite{MR4088470}). Moreover, we are interested in the application of pointwise regularity to free boundary problems.

For pointwise regularity, it can be tracked to the pioneering work of Caffarelli \cite{MR1005611}, in which the interior pointwise $C^{1,\alpha}$ and $C^{2,\alpha}$ regularity for fully nonlinear elliptic equations were proved. The boundary counterpart was proved by the authors \cite{MR4088470}. The pointwise regularity has been investigated extensively (see \cite{lian2020pointwise,MR3980853,MR2334822,MR2983006,MR3246039,MR2592289} etc.)
It reveals how the prescribed data influence the behavior of the solution near a concerned point.

The pointwise regularity has wide applications in other topics. For instance, the interior pointwise $C^{2,\alpha}$ regularity has been used to prove the optimal $C^{1,1}$ regularity for solutions in obstacle-type problems (see \cite[Proof of Theorem 1.2]{MR3198649}). Based on the boundary pointwise $C^{k,\alpha}$ ($k\geq 1$) regularity, one can give a direct and short proof of the higher regularity of free boundaries (see \cite{freeboundary} and \Cref{FreeBd}). The Liouville theorem on cones can be deduced with the aid of the boundary pointwise regularity (see \cite{MR_cones}). The boundary pointwise regularity has been employed in the establishment of the global $W^{2,p}$ regularity on a $C^{1,\alpha}$ domain for an appropriate $0<\alpha<1$ (see \cite{MR_LiLiZhang}).

With regard to parabolic equations, Krylov \cite{MR688919} first showed the boundary $C^{2,\alpha}$ a priori estimate. Wang \cite{MR1135923,MR1139064,MR1151267} obtained a series of fundamental estimates such as the Harnack inequality, interior pointwise $C^{1,\alpha}$ and $C^{2,\alpha}$ regularity and the boundary pointwise $C^{1,\alpha}$ regularity. This series of regularity for parabolic equations extended the results of Caffarelli \cite{MR1005611}. Wang's comprehensive research was a great help for our's study of the boundary pointwise regularity for parabolic equations. However, the proof of boundary $C^{1,\alpha}$ regularity in \cite{MR1139064} is a little complicated and the boundary pointwise $C^{2,\alpha}$ regularity was left. Adimurthi, Banerjee and Verma \cite{MR4073515} obtain the boundary $C^{2,\alpha}$ regularity by the classical method of flattening the boundary, which is not a pointwise regularity.


In this paper, we prove boundary pointwise $C^{k,\alpha}$ regularity for any $k\geq 1$ without flattening the boundary of the domain. Instead, we regard the curved boundary as a perturbation of a hyperplane. This viewpoint originates from \cite{MR3780142} and was further developed in \cite{MR4088470} and \cite{lian2020pointwise}.

Before attaining boundary regularity on a curved boundary, we need to derive enough regularity for problems on a flat boundary (called model problems, see \cref{model1} and \cref{e.C2a.F}). The $C^{1,\alpha}$ regularity has been obtained in \cite{MR1139064} and the $C^{2,\alpha}$ regularity was proved in \cite{MR4073515}. In this paper, we provide a clearer and simpler proof of the boundary $C^{1,\alpha}$ (see \Cref{Th1C1aF}) and $C^{2,\alpha}$ (see \Cref{Th1C2aF}) regularity for model problems.

The general problems are regarded as perturbations of model problems. Based on regularity for model problems, we prove the boundary pointwise regularity on general curved boundaries by the perturbation technique and the method of compactness. In addition, we present that if the derivatives of $u$ vanish at boundary, $u$ possesses higher boundary regularity. It was first observed in \cite{MR4088470} and used in \cite{freeboundary} to prove the higher regularity of free boundaries in obstacle-type problems for the Poisson's equation. As an extension, we will adopt this idea to give a direct and short proof of the higher regularity of free boundaries for fully nonlinear parabolic equations.

The most common method to prove the higher regularity of free boundaries is using the partial hodograph-Legendre transformation (see \cite{MR440187,MR0478079}). Here, we treat the obstacle problems based on the observation mentioned above. It is well-known that $u\in C^{k,\alpha}$ provided $\partial \Omega\in C^{k,\alpha}$ ($k\geq 1$). Besides, the regularity of $u$ may also lead to the regularity of $\partial \Omega$ because of $\varphi_i=u_i/u_n$ ($1\leq i< n$) and $\varphi_t=u_t/u_n$, where $\varphi$ is the representation function of $\partial \Omega$. Thus, we may have higher regularity for $u$ and then higher regularity for $\partial \Omega$ if we have more boundary conditions on $u$. The conditions $u=|\nabla u|=0$ on the boundary are enough exactly.

The results of this paper show the underlying relation between the regularity of solutions and the regularity of boundaries. The proofs demonstrate how overdetermined conditions lead to higher regularity of solutions and boundaries.

Before stating our main results, we introduce some notions and definitions for the pointwise regularity. In this paper, we use the parabolic terminology (see \Cref{no1.1} for details). For any $1\leq p <\infty$, let
\begin{equation*}
\|f\|^*_{L^p(\Omega)}= \left(\frac{1}{|\Omega|}\int_{\Omega} |f|^p\right)^{1/p},
\end{equation*}
where $\Omega\subset R^{n+1}$ is a bounded domain and $|\Omega|$ denotes its Lebesgue measure. The next is the standard pointwise characterization of a function in H\"{o}lder spaces.
\begin{definition}\label{d-f}
Let $\Omega\subset \mathbb{R}^{n+1}$ be a bounded set (may be not a domain) and $f:\Omega\rightarrow \mathbb{R}$ be a function. We say that $f$ is $C^{k,\alpha}$ ($k\geq 0, 0<\alpha\leq 1$) at $(x_0,t_0)\in \Omega$ or $f\in C^{k, \alpha}(x_0,t_0)$ if there exist constants $K,r_0>0$ and a polynomial $P\in \mathcal{P}_k$ (i.e., degree less than or equal to $k$) such that
\begin{equation}\label{m-holder}
  |f(x,t)-P(x,t)|\leq K|(x,t)-(x_0,t_0)|^{k+\alpha},~~\forall~(x,t)\in \Omega\cap Q_{r_0}(x_0,t_0).
\end{equation}
Then define for any $0\leq i\leq k$,
\begin{equation*}
  D^if(x_0,t_0)=D^iP(x_0,t_0)~~\mbox{and}~~
  \|f\|_{C^{k}(x_0,t_0)}=\sum_{i=0}^{k} |D^i P(x_0,t_0)|.
\end{equation*}
Moreover,
\begin{equation*}
\begin{aligned}
&[f]_{C^{k,\alpha}(x_0,t_0)}:=\min \left\{K\big | \cref{m-holder} ~\mbox{holds with}~P~\mbox{and}~K\right\}, \\
&\|f\|_{C^{k,\alpha}(x_0,t_0)}:=\|f\|_{C^{k}(x_0,t_0)}+[f]_{C^{k,\alpha}(x_0,t_0)}.
\end{aligned}
\end{equation*}

If $f\in C^{k, \alpha}(x,t)$ for any $(x,t)\in \Omega$ with the same $r_0$ and
\begin{equation*}
  \|f\|_{C^{k,\alpha}(\bar{\Omega})}= \sup_{(x,t)\in \Omega} \|f\|_{C^{k}(x,t)}+\sup_{(x,t)\in \Omega} [f]_{C^{k, \alpha}(x,t)}<+\infty,
\end{equation*}
we say that $f\in C^{k,\alpha}(\bar{\Omega})$.

In addition, we say that $f$ is $C^{-1,\alpha}$ at $(x_0,t_0)$ or $f\in C^{-1,\alpha}(x_0,t_0)$ if there exist constants $K,r_0>0$ such that
\begin{equation}\label{e.c-1}
\|f\|^*_{L^{n+1}(\Omega\cap Q_r(x_0,t_0) )}\leq Kr^{-1+\alpha}, ~\forall ~0<r<r_0.
\end{equation}
Then define
\begin{equation*}
\|f\|_{C^{-1,\alpha}(x_0,t_0)}=\min \left\{K \big | \cref{e.c-1} ~\mbox{holds with}~K\right\}.
\end{equation*}

If $f\in C^{-1, \alpha}(x,t)$ for any $(x,t)\in \Omega$ with the same $r_0$ and
\begin{equation*}
  \|f\|_{C^{-1,\alpha}(\bar{\Omega})}:= \sup_{(x,t)\in \Omega} \|f\|_{C^{-1,\alpha}(x,t)}<+\infty,
\end{equation*}
we say that $f\in C^{-1,\alpha}(\bar{\Omega})$.
\end{definition}

We also use the following classical definition of H\"{o}lder spaces (see \cite[Chap. 3]{Chen}).
\begin{definition}\label{de1.1}
Let $\Omega\subset \mathbb{R}^{n+1}$ be a bounded domain and $f:\Omega\rightarrow \mathbb{R}$. If $f\in C(\bar{\Omega})$, define
\begin{equation*}
\|f\|_{C(\bar{\Omega})}=\sup_{(x,t)\in \Omega}|f(x,t)|.
\end{equation*}
For $0<\alpha< 1$, if
\begin{equation*}
[f]_{C^{\alpha}(\bar{\Omega})}
:=\sup_{\mathop{(x,t),(y,s)\in \Omega,}\limits_{(x,t)\neq (y,s)}}\frac{|f(x,t)-f(y,s)|}{|(x,t)-(y,s)|^{\alpha}}<\infty,
\end{equation*}
we say that $f\in C^{\alpha}(\bar{\Omega})$. For $\alpha=1$, we have a similar definition (replace $\alpha$ by $1$ in the above equation). Then we write $f\in C^{0,1}(\bar{\Omega})$ and its semi-norm is denoted by $[f]_{C^{0,1}(\bar{\Omega})}$ (similarly hereinafter). We also introduce following semi-norm:
\begin{equation*}
[f]^{t}_{C^{\alpha}(\bar{\Omega})}:=\sup_{\mathop{(x,t),(x,s)\in \Omega,}\limits_{t\neq s}}
\frac{|f(x,t)-f(x,s)|}{|t-s|^{\alpha}}.
\end{equation*}

For $k\geq 1$, we say that $f\in C^k(\Omega)$ if $f$ has all derivatives of order $\leq k$ continuous in $\Omega$, i.e., all derivatives $D_{x}^{\sigma}D_{t}^{\gamma}f$ with $|\sigma|+2\gamma\leq k$ are continuous. Define
\begin{equation*}
  [f]_{C^{k}(\bar{\Omega})}=\left\{
  \begin{aligned}
&\sum_{|\sigma|+2\gamma= k}\|D_{x}^{\sigma}D_{t}^{\gamma}f\|_{C(\bar{\Omega})},~~
\mbox{if $k$ is even},~~\\
&\sum_{|\sigma|+2\gamma= k}\|D_{x}^{\sigma}D_{t}^{\gamma}f\|_{C(\bar{\Omega})}
+\sum_{|\sigma|+2\gamma= k-1}
[D_{x}^{\sigma}D_{t}^{\gamma}f]^{t}_{C^{1/2}(\bar{\Omega})},~~
\mbox{if $k$ is odd}.
  \end{aligned}
  \right.
\end{equation*}
We say that $f\in C^{k}(\bar\Omega)$ if
\begin{equation*}
\|f\|_{C^{k}(\bar{\Omega})}:=\|f\|_{C(\bar{\Omega})}+\sum_{i=1}^{k}[f]_{C^i(\bar{\Omega})}<\infty.
\end{equation*}

Similarly, for $k\geq 1$ and $0< \alpha< 1$, define
\begin{equation*}
  [f]_{C^{k,\alpha}(\bar{\Omega})}=\left\{
  \begin{aligned}
&\sum_{|\sigma|+2\gamma= k}[D_{x}^{\sigma}D_{t}^{\gamma}f]_{C^{\alpha}(\bar{\Omega})},~~
\mbox{if $k$ is even},~~\\
&\sum_{|\sigma|+2\gamma= k}[D_{x}^{\sigma}D_{t}^{\gamma}f]_{C^{\alpha}(\bar{\Omega})}
+\sum_{|\sigma|+2\gamma= k-1}
[D_{x}^{\sigma}D_{t}^{\gamma}f]^{t}_{C^{(1+\alpha)/2}(\bar{\Omega})},~~
\mbox{if $k$ is odd}.
  \end{aligned}
  \right.
\end{equation*}
We say that $f\in C^{k,\alpha}(\bar\Omega)$ if
\begin{equation*}
\|f\|_{C^{k,\alpha}(\bar{\Omega})}:=\|f\|_{C(\bar{\Omega})}+\sum_{i=1}^{k}[f]_{C^i(\bar{\Omega})}
+[f]_{C^{k,\alpha}(\bar{\Omega})}<\infty.
\end{equation*}
\end{definition}

\begin{remark}\label{re1.9}
It is known that if $\Omega$ is a smooth domain, the pointwise characterization \Cref{d-f} is equivalent to the classical \Cref{de1.1}. But we couldn't find a good reference for this fact.
\end{remark}
~\\

The next is a pointwise characterization of the smoothness of a domain's boundary. This definition is similar to \Cref{d-f}. That is, both definitions use polynomials to describe the smoothness. It was first introduced in \cite{MR4088470} for studying elliptic equations.

\begin{definition}\label{d-domain}
Let $\Omega$ be a bounded domain, $\Gamma\subset \partial \Omega$ be relatively open and $(x_0,t_0)\in \Gamma$. We say that $\Gamma$ is $C^{k,\alpha}$ ($k\geq 1, 0<\alpha\leq 1$) at $(x_0,t_0)$ or $\Gamma\in C^{k,\alpha}(x_0,t_0)$ if there exist constants $K,r_0>0$, a new coordinate system $\{x_1,...,x_n,t \}$ (by rotating and translating with respect to $x$ and only translating with $t$) and a polynomial $P\in \mathcal{P}_k$ with $P(0)=0$ and $DP(0)=0$ such that $(x_0,t_0)=0$ in this coordinate system,
\begin{equation}\label{e-re}
Q_{r_0} \cap \{(x',x_n,t)\big |x_n>P(x',t)+K|(x',t)|^{k+\alpha}\} \subset Q_{r_0}\cap \Omega
\end{equation}
and
\begin{equation}\label{e-re2}
Q_{r_0} \cap \{(x',x_n,t)\big |x_n<P(x',t)-K|(x',t)|^{k+\alpha}\} \subset Q_{r_0}\cap \Omega^c.
\end{equation}
Then define
\begin{equation*}
  \|\Gamma\|_{C^{k}(x_0,t_0)}=\sum_{i=2}^{k}|D^iP(0)|,
\end{equation*}
\begin{equation*}
[\Gamma]_{C^{k,\alpha}(x_0,t_0)}=\min \left\{K\big | \cref{e-re} ~\mbox{and}~ \cref{e-re2}~\mbox{hold with}~P~\mbox{and}~K\right\}
\end{equation*}
and
\begin{equation*}
\|\Gamma\|_{C^{k,\alpha}(x_0,t_0)}=\|\Gamma\|_{C^{k}(x_0,t_0)}+[\Gamma]_{C^{k,\alpha}(x_0,t_0)}.
\end{equation*}

If $\Gamma\in C^{k, \alpha}(x,t)$ for any $(x,t)\in \Gamma$ with the same $r_0$ and
\begin{equation*}
  \|\Gamma\|_{C^{k,\alpha}}:= \sup_{(x,t)\in \Gamma} \|\Gamma\|_{C^{k}(x,t)}+\sup_{(x,t)\in \Gamma}~[\Gamma]_{C^{k, \alpha}(x,t)}<+\infty,
\end{equation*}
we say that $\bar\Gamma\in C^{k,\alpha}$. If $\bar\Gamma'\in C^{k,\alpha}$ for any $\Gamma'\subset\subset\Gamma$, we denote $\Gamma\in C^{k,\alpha}$. If $\Gamma\in C^{k,\alpha}$ for any $k\geq1$ and $0<\alpha\leq 1$, we write $\Gamma\in C^{\infty}$.
\end{definition}

\begin{remark}\label{re1.13}
For a parabolic domain, there are two kinds of boundaries: the lateral boundary and the bottom boundary (see \cite{MR1135923}). If $\partial \Omega$ can be represented as a graph of some function $x_n=\varphi(x',t)$ near $(x_0,t_0)\in \partial \Omega$, $(x_0,t_0)$ must be on the lateral boundary. Similarly, the \cref{e-re} and \cref{e-re2} indicates that $(x_0,t_0)$ must belong to the lateral boundary of $\Omega$.

Note that in \Cref{d-domain}, $\partial \Omega$ doesn't need to be represented as a graph of some function and its geometrical property could be rather complicated. In fact, for any $k\geq 1$, one construct a domain $\Omega$ and $(x_0,t_0)\in \partial \Omega$ such that $\partial \Omega\in C^{k,1}(x_0,t_0)$ but any neighborhood of $(x_0,t_0)$ contains infinitely many points belonging to the bottom boundary.

If $\Omega=D\times (0,T]$ where $D\subset \mathbb{R}^n$ is a bounded domain, the definition of $\partial \Omega \in C^{k,\alpha}(x_0,t_0)$ is exactly equivalent to that for elliptic equations (see \cite[Definition 1.2]{MR4088470}).
\end{remark}

\begin{remark}\label{re1.2}
Throughout this paper, if we say that $f\in C^{k,\alpha}(x_0,t_0)$ ($\Gamma\in C^{k,\alpha}(x_0,t_0)$), we use $P_f$ ($P_{\Omega}$) to denote the corresponding polynomial in \Cref{d-f} (\Cref{d-domain}).
\end{remark}

\begin{remark}
We always assume that $0\in\partial \Omega$ and study the pointwise regularity at $0$ for \cref{e.Dirichlet}. In addition, if we use \Cref{d-f} and \Cref{d-domain} at $0$, we always assume that $r_0=1$, and \cref{e-re} and \cref{e-re2} hold if we say $\partial \Omega \cap Q_1\in C^{k,\alpha}(0)$.
\end{remark}
~\\

Now, we state our main results. Throughout this paper, a constant depending only on $n,\lambda$ and $\Lambda$ is called universal. The boundary pointwise $C^{1,\alpha}$ regularity is shown as follows.

\begin{theorem}[Boundary $C^{1,\alpha}$ regularity]\label{C1a}
  Let $0<\alpha<\bar\alpha$ and $u$ be a viscosity solution of
  \begin{equation*}
\left\{\begin{aligned}
&u\in S(\lambda,\Lambda,f)&& ~~\mbox{in}~~\Omega\cap Q_1;\\
&u=g&& ~~\mbox{on}~~\partial \Omega\cap Q_1.
\end{aligned}\right.
\end{equation*}
Suppose that $f\in C^{-1,\alpha}(0)$, $g\in C^{1,\alpha}(0)$ and $\partial \Omega \cap Q_1 \in C^{1,\alpha}(0)$.

Then $u\in C^{1,\alpha}(0)$, i.e., there exists $P\in \mathcal{P}_1$ such that
\begin{equation}\label{e.C1a}
|u(x,t)-P(x)|\leq C|(x,t)|^{1+\alpha}\left(\|u\|_{L^{\infty}(\Omega\cap Q_1)}+
\|f\|_{C^{-1,\alpha}(0)}+\|g\|_{C^{1,\alpha}(0)}\right),~~\forall~(x,t)\in \Omega \cap Q_{1},
\end{equation}
\begin{equation*}
  |Du(0)|\leq C\left(\|u\|_{L^{\infty}(\Omega\cap Q_1)}+
\|f\|_{C^{-1,\alpha}(0)}+\|g\|_{C^{1,\alpha}(0)}\right)
\end{equation*}
and
\begin{equation}\label{e1.7}
P(x',0)\equiv P_g(x',0),
\end{equation}
where $C$ depends only on $n, \lambda,\Lambda, \alpha$ and $\|\partial \Omega\cap Q_1\|_{C^{1,\alpha}(0)}$.
\end{theorem}

\begin{remark}\label{re1.1-ref3}
The constant $\bar\alpha$ (fixed throughout this paper) is universal, which is originated from the model problem (see \Cref{Th1C1aF}).
\end{remark}

\begin{remark}\label{re1.3}
Since $P$ is a linear polynomial, by the definition, $P$ is independent of $t$.
\end{remark}

\begin{remark}\label{re1.1-ref}
In this paper, the $L^{n+1}$-viscosity solution is employed. Thus, $f$ may not be continuous with respect to $x$ or $t$. For the notions of viscosity solutions and detailed theories, we refer to \cite{MR1351007,MR1789919,MR1118699,MR1135923}.
\end{remark}

\begin{remark}\label{re1.1-ref2}
Wang \cite[Theorem 2.1]{MR1139064} obtained the pointwise $C^{1,\alpha}$ regularity for \emph{some} $\alpha<\bar{\alpha}$ rather than for \emph{any} $0<\alpha< \bar{\alpha}$ here. In particular, for the heat equation $u_t-\Delta u=f$, we have $\bar{\alpha}=1$ and then obtain the boundary pointwise $C^{1,\alpha}$ regularity for any $0<\alpha<1$.
\end{remark}
~\\

For boundary $C^{2,\alpha}$ regularity, we need the following definition to estimate the oscillation of $F$ in $(x,t)$.
\begin{definition}\label{d-FP}
Let  $\Omega$ be a bounded domain and $p_0=(x_0,t_0)\in \bar{\Omega}$. Suppose that there exist $r_0>0$ and a fully nonlinear operator $F_{p_0}:\mathcal{S}^{n}\rightarrow \mathbb{R}$ such that for any $M\in \mathcal{S}^{n}$ and $(x,t)\in \bar\Omega\cap \bar Q_{r_0}(x_0,t_0)$,
\begin{equation}\label{e.C2a.beta}
|F(M,x,t)-F_{p_0}(M)|\leq \beta((x,t),(x_0,t_0))\|M\|+\gamma((x,t),(x_0,t_0)),
\end{equation}
where $F_{p_0}$ is uniformly elliptic with $\lambda$ and $\Lambda$.

If for some $K>0$,
\begin{equation}\label{e1.1-n}
\beta((x,t),(x_0,t_0))\leq K|(x,t)-(x_0,t_0)|^{\alpha}
\end{equation}
and
\begin{equation}\label{e1.1-2-n}
\gamma((x,t),(x_0,t_0))\leq K|(x,t)-(x_0,t_0)|^{\alpha},
\end{equation}
where $0<\alpha<\bar{\alpha}$, we say that $F\in C^{\alpha}(x_0,t_0)$ and define
\begin{equation*}
  \|\beta\|_{C^{\alpha}(x_0,t_0)}=\min \left\{K\big | \cref{e1.1-n} ~\mbox{holds with}~K\right\}
\end{equation*}
and
\begin{equation*}
  \|\gamma\|_{C^{\alpha}(x_0,t_0)}=\min \left\{K\big | \cref{e1.1-2-n} ~\mbox{holds with}~K\right\}.
\end{equation*}

If $F\in C^{\alpha}(x,t)$ for any $(x,t)\in \bar{\Omega}$ with the same $r_0$ and
\begin{equation*}
\|\beta\|_{C^{\alpha}(\bar{\Omega})}:=\sup_{(x,t)\in \bar\Omega} \|\beta\|_{C^{\alpha}(x,t)}<+\infty,
\end{equation*}
\begin{equation*}
\|\gamma\|_{C^{\alpha}(\bar{\Omega})}:=\sup_{(x,t)\in \bar\Omega} \|\gamma\|_{C^{\alpha}(x,t)}<+\infty,
\end{equation*}
\begin{equation*}
\|F\|_{\bar{\Omega}}:=\sup_{p\in \bar\Omega} |F_{p}(0)|<+\infty,
\end{equation*}
we say that $F\in C^{\alpha}(\bar{\Omega})$.
\end{definition}

\begin{remark}\label{re.ocs.F1}
We remark here that if $F$ is continuous in $(x,t)$, we may use $F(M,x_0,t_0)$ to measure the oscillation of $F$ in $(x,t)$ in \cref{e.C2a.beta} instead of $F_{p_0}(M)$.
\end{remark}

\begin{remark}\label{re.ocs.F2}
As before, we always assume that $r_0=1$ throughout this paper when we use \Cref{d-FP}.
\end{remark}
~\\

Now, we give the boundary pointwise $C^{2,\alpha}$ regularity.
\begin{theorem}[Boundary $C^{2,\alpha}$ regularity]\label{C2a-1}
Let $0<\alpha<\bar\alpha$ and $u$ be a viscosity solution of
  \begin{equation}\label{e.C2a.u-1}
\left\{\begin{aligned}
&u_t-F(D^2u,x,t)=f&& ~~\mbox{in}~~\Omega\cap Q_1;\\
&u=g&& ~~\mbox{on}~~\partial \Omega \cap Q_1.
\end{aligned}\right.
\end{equation}
Suppose that $F\in C^{\alpha}(0)$ (i.e., $F$ satisfies \cref{e.C2a.beta} with $\beta,\gamma\in C^{\alpha}(0)$), $f\in C^{\alpha}(0)$, $g\in C^{2,\alpha}(0)$ and $\partial \Omega\cap Q_1 \in C^{1,\alpha}(0)$. In addition, assume that $u\in C^{1,\alpha}(0)$ and
\begin{equation*}
u(0)=|Du(0)|=|Dg(0)|=0.
\end{equation*}

Then $u\in C^{2,\alpha}(0)$, i.e., there exists $P\in \mathcal{HP}_2$ such that
\begin{equation}\label{e.C2a-1}
|u(x,t)-P(x,t)|\leq CU|(x,t)|^{2+\alpha},~~\forall~(x,t)\in \Omega \cap Q_{1},
\end{equation}
\begin{equation}\label{e.C2a.2-1}
|D^2u(0)| \leq CU,
\end{equation}
\begin{equation*}
U=\|u\|_{L^{\infty}(\Omega\cap Q_1)}+|F_0(0)|+
\|f\|_{C^{\alpha}(0)}+\|g\|_{C^{2,\alpha}(0)}
+\|\gamma\|_{C^{\alpha}(0)}
\end{equation*}
and
\begin{equation}\label{e.C2a.3-1}
\begin{aligned}
&P_t-F_0(D^2 P)-P_f=0,\\
&P(x',0,t)\equiv P_g(x',0,t),
\end{aligned}
\end{equation}
where $F_0$ is from \cref{e.C2a.beta} with $p_0=0$ and $C$ depends only on $n, \lambda,\Lambda,\alpha$, $\|\beta\|_{C^{\alpha}(0)}$ and $\|\partial \Omega\cap Q_1\|_{C^{1,\alpha}(0)}$.
\end{theorem}

\begin{remark}\label{re1.5}
One feature of the theorem is that $u\in C^{2,\alpha}(0)$ even for $\partial \Omega\in C^{1,\alpha}(0)$, which is key to the higher regularity of free boundaries (see \Cref{FreeBd}). This was first shown in \cite{MR4088470} and the higher order counterpart was obtained in \cite{freeboundary}.
\end{remark}
~\\

As an easy corollary, we have
\begin{theorem}[Boundary $C^{2,\alpha}$ regularity]\label{C2a}
Let $0<\alpha<\bar\alpha$ and $u$ be a viscosity solution of
  \begin{equation}\label{e.C2a.u}
\left\{\begin{aligned}
&u_t-F(D^2u,x,t)=f&& ~~\mbox{in}~~\Omega\cap Q_1;\\
&u=g&& ~~\mbox{on}~~\partial \Omega \cap Q_1.
\end{aligned}\right.
\end{equation}
Suppose that $F\in C^{\alpha}(0)$, $f\in C^{\alpha}(0)$, $g\in C^{2,\alpha}(0)$ and $\partial \Omega\cap Q_1 \in C^{2,\alpha}(0)$.

Then $u\in C^{2,\alpha}(0)$, i.e., there exists $P\in \mathcal{P}_2$ such that
\begin{equation}\label{e.C2a}
|u(x,t)-P(x,t)|\leq CU|(x,t)|^{2+\alpha},~~\forall~(x,t)\in \Omega \cap Q_{1},
\end{equation}
\begin{equation}\label{e.C2a.2}
  |Du(0)|+|D^2u(0)| \leq CU,
\end{equation}
\begin{equation*}
U=\|u\|_{L^{\infty}(\Omega\cap Q_1)}+|F_0(0)|+
\|f\|_{C^{\alpha}(0)}+\|g\|_{C^{2,\alpha}(0)}
+\|\gamma\|_{C^{\alpha}(0)}
\end{equation*}
and
\begin{equation}\label{e.C2a.3}
\begin{aligned}
&P_t-F_0(D^2 P)-P_f=0,\\
&\mathbf{\Pi}_{2}\left(P(x',P_{\Omega}(x',t),t)\right)\equiv
\mathbf{\Pi}_{2}\left(P_g(x',P_{\Omega}(x',t),t)\right),
\end{aligned}
\end{equation}
where $C$ depends only on $n, \lambda,\Lambda,\alpha$, $\|\beta\|_{C^{\alpha}(0)}$ and $\|\partial \Omega\cap Q_1\|_{C^{2,\alpha}(0)}$.
\end{theorem}


\begin{remark}\label{re.C2a1}
For the boundary $C^{2,\alpha}$ regularity of the parabolic equation, Adimurthi, Banerjee and Verma \cite{MR4073515} attained the regularity by means of a transformation in order to flatten the boundary and then use the result for the homogeneous equation. This method is classical but makes the equation complicated and the result is not pointwise.

In this paper, we use the idea of perturbation to obtain the pointwise $C^{2,\alpha}$ regularity in a curved boundary without flattening the boundary. This idea was first used for oblique derivative problems in \cite{MR3780142} and for Dirichlet problems in \cite{MR4088470}. Additionally, our proof of the $C^{2,\alpha}$ regularity for the homogenous equation is much simpler than that of \cite{MR4073515} (see \Cref{Th1C2aF}).
\end{remark}

\begin{remark}\label{re1.6}
The \cref{e.C2a.3} demonstrates the peculiarity of pointwise regularity. Note that $P,P_f,P_g$ and $P_{\Omega}$ are the ``Taylor polynomials'' of $u,f,g$ and $\partial\Omega$ at $0$ respectively. Hence, \cref{e.C2a.3} is a pointwise version of \cref{e.C2a.u} at $0$ and it shows clearly the relation between the solution and the prescribed data at $0$.
\end{remark}

\begin{remark}\label{re1.7}
Note that the degree of the polynomial $P(x',P_{\Omega}(x',t),t)$ may be greater than $2$. Since we only have $C^{2,\alpha}$ estimate at $0$, we need to take a projection to $\mathcal{P}_2$ in the second equality of \cref{e.C2a.3}.
\end{remark}
~\\


For higher $C^{k,\alpha}$ $(k\geq 2)$ regularity, we introduce the following definition to estimate the oscillation of $F$ in $(x,t)$.
\begin{definition}\label{d-FP-Ck}
Let  $\Omega$ be a bounded domain and $p_0=(x_0,t_0)\in \bar{\Omega}$. Suppose that there exist $r_0>0$ and a fully nonlinear operator $F_{p_0}:\mathcal{S}^{n}\times \bar\Omega\cap \bar Q_{r_0}(x_0,t_0)\rightarrow \mathbb{R}$ such that for any $M\in \mathcal{S}^{n}$ and $(x,t)\in \bar\Omega\cap \bar Q_{r_0}(x_0,t_0)$,
\begin{equation}\label{e.Cka.beta}
|F(M,x,t)-F_{p_0}(M,x,t)|\leq \beta((x,t),(x_0,t_0)) \|M\|+\gamma((x,t),(x_0,t_0)),
\end{equation}
where \\
(i) $F_{p_0}$ is convex in $M$ and has the same uniform ellipticity as $F$;\\
(ii) $F_{p_0} \in C^{k,1}(\mathcal{S}^n\times \bar{\Omega}\cap \bar{Q}_{r_0}(x_0,t_0))$;\\
(iii) there exists a constant $K_0>0$ such that for any $M_1,M_2\in \mathcal{S}^{n}$ and $(x_1,t_1),(x_2,t_2)\in \bar\Omega\cap \bar Q_{r_0}(x_0,t_0)$,
\begin{equation}\label{e1.8}
|F_{p_0,ij}(M_1,x_1,t_1)-F_{p_0,ij}(M_2,x_2,t_2)|\leq K_0 (\|M_1-M_2\|+|(x_1,t_1)-(x_2,t_2)|).
\end{equation}

If for some $K>0$,
\begin{equation}\label{e1.1}
\beta((x,t),(x_0,t_0))\leq K|(x,t)-(x_0,t_0)|^{k+\alpha}
\end{equation}
and
\begin{equation}\label{e1.1-2}
\gamma((x,t),(x_0,t_0))\leq K|(x,t)-(x_0,t_0)|^{k+\alpha},
\end{equation}
where $0<\alpha<1$, we say that $F\in C^{k,\alpha}(x_0,t_0)$ and define
\begin{equation*}
  \|\beta\|_{C^{k,\alpha}(x_0,t_0)}=\min \left\{K\big | \cref{e1.1} ~\mbox{holds with}~K\right\}
\end{equation*}
and
\begin{equation*}
  \|\gamma\|_{C^{k,\alpha}(x_0,t_0)}=\min \left\{K\big | \cref{e1.1-2} ~\mbox{holds with}~K\right\}.
\end{equation*}

In addition, for $r>0$, denote
\begin{equation*}
\omega(r)=\|F_{p_0}\|_{C^{k,1}(\bar{\textbf{B}}_r\times \bar\Omega\cap \bar Q_{r_0}(x_0,t_0))},
\end{equation*}
where $\textbf{B}_r=\left\{M\big| \|M\|<r\right\}$.

If $F\in C^{k, \alpha}(x,t)$ for any $(x,t)\in \bar{\Omega}$ with the same $r_0$ and $\omega$, and
\begin{equation*}
\|\beta\|_{C^{k,\alpha}(\bar{\Omega})}:=\sup_{(x,t)\in \bar\Omega} \|\beta\|_{C^{k,\alpha}(x,t)}<+\infty,
\end{equation*}
\begin{equation*}
\|\gamma\|_{C^{k,\alpha}(\bar{\Omega})}:=\sup_{(x,t)\in \bar\Omega} \|\gamma\|_{C^{k,\alpha}(x,t)}<+\infty,
\end{equation*}
we say that $F\in C^{k, \alpha}(\bar{\Omega})$. If $F \in C^{k,\alpha}$ for any $k\geq1$ and $0<\alpha< 1$, we write $F \in C^{\infty}(\bar{\Omega})$.
\end{definition}

\begin{remark}
Note that the conditions (i) and (ii) are natural to derive higher regularity. In addition, the linearization method is used to prove higher regularity and (iii) states that the linearized operator is Lipschitz continuity. It seems more natural to replace (iii) by the following condition:
\begin{equation*}
|F_{p_0}(M,x_1,t_1)-F_{p_0}(M,x_2,t_2)|\leq K_0 (\|M\|+1)|(x_1,t_1)-(x_2,t_2)|.
\end{equation*}
However, we can't prove the higher regularity under the above condition.
\end{remark}

\begin{remark}\label{re1.4}
Note that $F_{p_0} \in C^{k,1}(\mathcal{S}^n\times \bar{\Omega})$ is used in the parabolic terminology. That is, for any $A\subset\subset \mathcal{S}^n$, the semi-norm $[F]_{C^{k,1}(\bar A\times\bar{\Omega})}$ is defined as
\begin{equation*}
\left\{
  \begin{aligned}
&\sum_{|\zeta|+|\sigma|+2\gamma= k}
[D_M^{\zeta}D_{x}^{\sigma}D_{t}^{\gamma}F]_{C^{0,1}(\bar A\times\bar{\Omega})},~~
\mbox{if $k$ is even},~~\\
&\sum_{|\zeta|+|\sigma|+2\gamma= k}
[D_M^{\zeta}D_{x}^{\sigma}D_{t}^{\gamma}F]_{C^{0,1}(\bar A\times\bar{\Omega})}
+\sum_{|\zeta|+|\sigma|+2\gamma= k-1}
[D_M^{\zeta}D_{x}^{\sigma}D_{t}^{\gamma}F]^{t}_{C^{0,1}(\bar A\times\bar{\Omega})},~~
\mbox{if $k$ is odd}.
  \end{aligned}
  \right.
\end{equation*}
\end{remark}
~\\

Next, we give the boundary pointwise $C^{k,\alpha}$ regularity.

\begin{theorem}[Boundary $C^{k,\alpha}$ regularity]\label{Ckla}
Let $0<\alpha<1$, $k,l\geq 1$, $k+l\geq 3$ and $u$ be a viscosity solution of
\begin{equation}\label{e.Ckla.1}
 \left\{ \begin{aligned}
&u_t-F(D^2 u,x,t)=f ~~&&\mbox{in}~~\Omega\cap Q_1;\\
&u=g ~~&&\mbox{on}~~\partial \Omega\cap Q_1.
 \end{aligned}\right.
\end{equation}
Suppose that $F\in C^{k+l-2,\alpha}(0)$, $f\in C^{k+l-2,\alpha}(0)$, $g\in C^{k+l,\alpha}(0)$ and $\partial \Omega\cap Q_1\in C^{l,\alpha}(0)$. Moreover, assume that $u\in C^{k,\alpha}(0)$ and
\begin{equation*}
u(0)=\cdots=|D^ku(0)|=|Dg(0)|=\cdots=|D^kg(0)|=0.
\end{equation*}

Then $u\in C^{k+l,\alpha}(0)$, i.e., there exists $P\in \mathcal{P}_{k+l}$ such that
\begin{equation*}
  \begin{aligned}
&|u(x,t)-P(x,t)|\leq C |(x,t)|^{k+l+\alpha},~\forall ~(x,t)\in \Omega\cap Q_{1},\\
&|D^{k+1}u(0)|+\cdots+|D^{k+l}u(0)|\leq C
  \end{aligned}
\end{equation*}
and
\begin{equation}\label{e1.2}
\begin{aligned}
&|P_t-F_0(D^2 P,x,t)-P_f|\leq C|(x,t)|^{k+l-1},\\
&\mathbf{\Pi}_{k+l}\left(P(x',P_{\Omega}(x',t),t)\right)\equiv
\mathbf{\Pi}_{k+l}\left(P_g(x',P_{\Omega}(x',t),t)\right),
\end{aligned}
\end{equation}
where $C$ depends only on $k,l,n,\lambda,\Lambda,\alpha,K_0$, $\omega$, $\|\beta\|_{C^{k+l-2,\alpha}(0)}$, $\|\partial \Omega\cap Q_1\|_{C^{l,\alpha}(0)}$, $\|u\|_{L^{\infty}(\Omega_1)}$, $\|f\|_{C^{k+l-2,\alpha}(0)}, \|g\|_{C^{k+l,\alpha}(0)}$ and  $\|\gamma\|_{C^{k+l-2,\alpha}(0)}$.
\end{theorem}



\begin{remark}\label{re1.12}
In fact, the expression of the polynomial $P$ can be written explicitly:
\begin{equation}\label{e1.5}
  P(x,t)=P_g(x,t)+\mathbf{\Pi}_{k+l}\left(\sum_{\mathop{k+1\leq |\sigma|+2\gamma\leq k+l,}\limits_{\sigma_n\geq 1}} \frac{a_{\sigma\gamma}}{\sigma!\gamma !}x^{\sigma-e_n}t^{\gamma}\left(x_n-P_{\Omega}(x',t)\right)\right),
\end{equation}
where $a_{\sigma\gamma}$ are constants.
\end{remark}


\begin{remark}\label{re1.14}
As before, the \cref{e1.2} can be regarded as a polynomial version of \cref{e.Ckla.1}. Note that $F_0$ is a fully nonlinear operator and hence $F_0(D^2P,x,t)$ may be not a polynomial. Thus, we can only obtain the estimate as in \cref{e1.2} rather than an equality.
\end{remark}

\begin{remark}\label{re1.8}
For the usual local estimates, once the $C^{2,\alpha}$ regularity is established, the higher regularity can be deduced easily by taking difference quotient (or derivative) on both sides of the equation. However, for pointwise regularity, the above method fails since the prescribed data are smooth only in the sense of pointwise definition. Hence, the higher pointwise regularity is not trivial. Instead, its proof is more complicated.
\end{remark}
~\\

As a consequence, we have the following boundary pointwise regularity.
\begin{theorem}[Boundary $C^{k,\alpha}$ regularity]\label{Cka}
Let $0<\alpha<1$, $k\geq 3$ and $u$ be a viscosity solution of
\begin{equation}\label{e.Cka.uF}
 \left\{ \begin{aligned}
&u_t-F(D^2u,x,t)=f ~~&&\mbox{in}~~\Omega\cap Q_1;\\
&u=g ~~&&\mbox{on}~~\partial \Omega\cap Q_1.
 \end{aligned}\right.
\end{equation}
Suppose that $F\in C^{k-2,\alpha}(0)$, $f\in C^{k-2,\alpha}(0)$, $g\in C^{k,\alpha}(0)$ and $\partial \Omega\cap Q_1\in C^{k,\alpha}(0)$.

Then $u\in C^{k,\alpha}(0)$, i.e., there exists $P\in \mathcal{P}_k$ such that
\begin{equation*}
  \begin{aligned}
&|u(x,t)-P(x,t)|\leq C |(x,t)|^{k+\alpha},~\forall ~(x,t)\in \Omega\cap Q_{1},\\
&|Du(0)|+\cdots+|D^{k}u(0)|\leq C
  \end{aligned}
\end{equation*}
and
\begin{equation*}
\begin{aligned}
&|P_t-F_0(D^2 P,x,t)-P_f|\leq C|(x,t)|^{k-1},\\
&\mathbf{\Pi}_{k}\left(P(x',P_{\Omega}(x',t),t)\right)\equiv
\mathbf{\Pi}_{k}\left(P_g(x',P_{\Omega}(x',t),t)\right).
\end{aligned}
\end{equation*}
where $C$ depends only on $k,n,\lambda,\Lambda,\alpha,K_0,\omega$, $\|\beta\|_{C^{k-2,\alpha}(0)}$, $\|\partial \Omega\cap Q_1\|_{C^{k,\alpha}(0)}$, $\|u\|_{L^{\infty}(\Omega\cap Q_1)}$, $\|f\|_{C^{k-2,\alpha}(0)}$, $\|g\|_{C^{k,\alpha}(0)}$ and $\|\gamma\|_{C^{k-2,\alpha}(0)}$.
\end{theorem}

As an application of \Cref{C2a-1} and \Cref{Ckla} to the higher regularity of free boundaries in obstacle-type problems, we have
\begin{theorem}\label{FreeBd}
Let $u$ be a viscosity solution of
\begin{equation}\label{e.FBd.1}
\left\{\begin{aligned}
   &u_t-F(D^2 u,x,t)=1~~~~&&\mbox{in}~~\Omega \cap Q_1;\\
   &u=|Du|=0~~~~&&\mbox{on}~~\partial \Omega\cap Q_1.\\
\end{aligned}\right.
\end{equation}
Suppose that $F\in C^{\infty}(\mathcal{S}^n\times \bar\Omega\cap \bar Q_{1})$ is convex in $M$, $F(0,0,0)=0$ and for some $0<\alpha<\bar{\alpha}$, for any $M\in \mathcal{S}^n$, $ (x_1,t_1),(x_2,t_2)\in \bar{\Omega}\cap \bar Q_1$,
\begin{equation}\label{e1.10}
|F(M,x_1,t_1)-F(M,x_2,t_2)|\leq K|(x_1,t_1)-(x_2,t_2)|^{\alpha}(\|M\|+1).
\end{equation}
Assume further that and $\partial \Omega\cap Q_1\in C^{1,\alpha}$. Then $u\in C^{\infty} (\bar\Omega\cap \bar Q_{1})$ and $\partial \Omega\cap Q_1\in C^{\infty}$.
\end{theorem}

\begin{remark}\label{re.freeb.1}
Generally, for the obstacle-type problem \cref{e.FBd.1}, the higher regularity of free boundaries is obtained by employing the partial hodograph-Legendre transformation. The \cref{e.FBd.1} is transformed to another fully nonlinear parabolic equations with a flat boundary. Then the regularity of the free boundary is improved by the higher regularity of the solution of the transformed equation \cite{MR440187}. This method requires that the free boundary $\partial \Omega\cap Q_1 \in C^1$ and the solution $u \in C^2$. In addition, the transformed equation is in a more complicated form.

As a comparison, \Cref{FreeBd} requires that $\partial \Omega\cap Q_1 \in C^{1,\alpha}$ and $u$ is a viscosity solution. Another feature of \Cref{FreeBd} is that its proof is direct and rather simple.

De Silva and Savin \cite{MR3393271} gave a direct and simple proof of higher regularity of solutions and free boundaries based on a higher order boundary Harnack inequality. However, this method is not applicable to nonlinear equations.
\end{remark}
~\\

At the end of this section, we remark that all results in this paper can be extended by a similar method to equations in the general forms
\begin{equation*}
  u_t-F(D^2u, Du,u,x,t)=f
\end{equation*}
under appropriate structure conditions (see \cite{lian2020pointwise}).

This paper is organized as follows. \Cref{sec:C1a} is devoted to the boundary pointwise $C^{1,\alpha}$ regularity for the parabolic equations. \Cref{sec:C2a} is concerned with the boundary pointwise $C^{2,\alpha}$ regularity. In \Cref{Sec:Cka_freeb}, we prove the boundary pointwise $C^{k,\alpha}$ regularity for Dirichlet problems and the higher regularity for free boundaries of obstacle-type problems. Notations used in this paper are listed below.
\begin{notation}\label{no1.1}
\begin{enumerate}~~\\
\item $\{e_i\}^{n}_{i=1}$: the standard basis of $\mathbb{R}^n$, i.e., $e_i=(0,...0,\underset{i^{th}}{1},0,...0)$.
\item $x'=(x_1,x_2,...,x_{n-1})$, $x=(x',x_n) \in \mathbb{R}^{n}$ and $(x,t)=(x',x_n,t)\in \mathbb{R}^{n+1}$.
\item $|x|=\left(\sum_{i=1}^{n} x_i^2\right)^{1/2}$ for $x\in \mathbb{R}^n$. The parabolic norm (or the parabolic distance from $(x,t)$ to $0$) is defined as
\begin{equation*}
|(x,t)|=\left\{\begin{aligned}
&(|x|^2+|t|)^{1/2} &&~~t\leq 0;\\
&\infty && ~~t>0.
\end{aligned}\right.
\end{equation*}
for $(x,t)\in \mathbb{R}^{n+1}$.
\item $\mathbb{R}^n_+=\{x\in \mathbb{R}^n\big|x_n>0\}$ and $\mathbb{R}^{n+1}_+=\{(x,t)\in \mathbb{R}^{n+1}\big|x_n>0\}$.
\item $B_r(x_0)=B(x_0,r)=\{x\in \mathbb{R}^{n}\big| |x-x_0|<r\}$ and $B_r=B_r(0)$.
\item $B_r^+(x_0)=B_r(x_0)\cap \mathbb{R}^n_+$ and $B_r^+=B^+_r(0)$.
\item $T_r(x_0)\ =\{x' \in \mathbb{R}^{n-1}\big| |x'-x_0'|<r\}$ and $T_r=T_r(0)$.
\item $Q_r(x_0,t_0)=Q((x_0,t_0),r)=B_r(x_0)\times (t_0-r^2,t_0]$ and $Q_r=Q_r(0)$.
\item $Q_r^+(x_0,t_0)=Q_r(x_0,t_0)\cap \mathbb{R}^{n+1}_+$ and $Q_r^+=Q^+_r(0)$.
\item $S_r(x_0,t_0)\ =T_r(x_0)\times (t_0-r^2,t_0]$ and $S_r=S_r(0)$.
\item $\Omega_r=\Omega\cap Q_r$ and $(\partial\Omega)_r=\partial\Omega\cap Q_r$.
\item $A^c$: the complement of $A$; $\bar A $: the closure of $A$, where $ A\subset \mathbb{R}^{n+1}$.
\item $\mathcal{S}^{n}$: the set of $n\times n$ symmetric matrices and $\|A\|=$ the spectral radius of $A$ for any $A\in \mathcal{S}^{n}$.
\item $D^k \varphi$: the $k$-th order derivatives of a function $\varphi$, i.e.,
     \begin{equation*}
    D^k \varphi=\left\{D_{x}^{\sigma}D_{t}^{\gamma} \varphi:|\sigma|+2\gamma=k\right\},
     \end{equation*}
     where the standard multi-index notation is used. Define
     \begin{equation*}
       |D^k\varphi |
    =\left(\sum_{|\sigma|+2\gamma=k}|D_{x}^{\sigma}D_{t}^{\gamma}\varphi|^2\right)^{1/2}.
     \end{equation*}
\item $\varphi_t=\partial  \varphi/\partial t$, $\varphi _i=\partial \varphi/\partial x _{i}$ $(1\leq i\leq n)$, $\varphi _{ij} =\partial ^{2}\varphi/\partial x_{i}\partial x_{j}$ $(1\leq i,j\leq n)$ and we also use similar notations for higher order derivatives.
\item For $F(M,x,t):\mathcal{S}^n\times \Omega\rightarrow \mathbb{R}$, we use $F_i,F_t,F_{ij}$ to denote the derivative with respect to $x_i,t$ and $M_{ij}$ respectively.
\item $\mathcal{P}_k (k\geq 0):$ the set of polynomials of degree less than or equal to $k$. That is, any $P\in \mathcal{P}_k$ can be written as
\begin{equation*}
P(x,t)=\sum_{|\sigma|+2\gamma\leq k}\frac{a_{\sigma\gamma}}{\sigma!\gamma!}x^{\sigma}t^{\gamma}
\end{equation*}
where $a_{\sigma\gamma}$ are constants. Define
\begin{equation*}
\|P\|= \sum_{|\sigma|+2\gamma \leq k}|a_{\sigma\gamma}|.
\end{equation*}
\item $\mathcal{HP}_k (k\geq 0):$ the set of homogeneous polynomials of degree $k$. That is, any $P\in \mathcal{HP}_k$ can be written as
\begin{equation*}
P(x,t)=\sum_{|\sigma|+2\gamma= k}\frac{a_{\sigma\gamma}}{\sigma!\gamma!}x^{\sigma}t^{\gamma}.
\end{equation*}
\item $\mathbf{\Pi}_k:$ The projection from $\mathcal{P}_l$ to $\mathcal{P}_k$ for $l\geq k$. That is, if $P\in \mathcal{P}_l$ is written as
\begin{equation*}
P(x,t)=\sum_{|\sigma|+2\gamma\leq l}\frac{a_{\sigma\gamma}}{\sigma!\gamma!}x^{\sigma}t^{\gamma},
\end{equation*}
then
\begin{equation*}
\mathbf{\Pi}_kP(x,t)=\sum_{|\sigma|+2\gamma\leq k}\frac{a_{\sigma\gamma}}{\sigma!\gamma!}x^{\sigma}t^{\gamma}.
\end{equation*}
\end{enumerate}
\end{notation}
~\\

\section{Boundary $C^{1,\alpha}$ regularity}\label{sec:C1a}

In this section, we give the proof of boundary pointwise $C^{1,\alpha}$ regularity. First, we develop the boundary regularity for the model problem (i.e. homogeneous equation with a flat boundary). We prove it by comparing the solution with $x_n$ and an iteration procedure (\Cref{Th1C1aF}).

Then we use the idea of perturbation to get the regularity for general boundaries. Roughly speaking, we first obtain a ``$C^{1,\alpha}$ like'' estimate on some scale with small prescribed data by the compactness method (\Cref{L.K1}). Next, we get the boundary pointwise $C^{1,\alpha}$ regularity by an iteration procedure (\Cref{L.Sca1}) and some appropriate normalizations.


\subsection{Boundary $C^{1,\alpha}$ regularity on flat boundaries}

In this subsection, we show the boundary $C^{1,\alpha}$ regularity on a flat boundary. The strategy is comparing the solution with $x_n$, which is similar to that for elliptic equations (see \cite{MR_Israel}). First, we prove the Hopf lemma and the boundary Lipschitz regularity.

\begin{lemma}[Hopf Lemma]\label{L3}
  Let $u(e_n/2,-3/4)=1$ and $u\geq 0$ be a viscosity solution of
  \begin{equation}\label{e.S0}
\left\{\begin{aligned}
&u\in S(\lambda,\Lambda,0)&& ~~\mbox{in}~~Q_1^+;\\
&u=0&& ~~\mbox{on}~~S_1.
\end{aligned}\right.
\end{equation}
Then
\begin{equation}\label{e.Hopf}
u(x,t)\geq Cx_n~~~~\mbox{ in}~~Q_{1/2}^+,
\end{equation}
where $C$ is universal.
\end{lemma}
\proof Set $r=1/2$ and $Q=B_{r/2}(re_n,0)\times (-r^2,0]$. By the Harnack inequality \cite[Theorem 4.18]{MR1135923} and noting $u(e_n{\AE}/2,-3/4)=1$, we have
\begin{equation*}
   u\geq C_0~~~~\mbox{ on}~~\partial Q,
\end{equation*}
where $C_0>0$ is universal.

Let
\begin{equation}\label{e1.3}
  v(x,t)=C\left( \left(\frac{1}{t+r^2}e^{-\frac{|x-re_n|^2}{t+r^2}}\right)^{\beta}
-\frac{e^{-\beta}}{r^{2\beta}}\right).
\end{equation}
By taking  $\beta$ large enough and $C$ small enough, $v$ satisfies
\begin{equation}\label{e.Hopf.v}
\left\{\begin{aligned}
&v_t-\mathcal{M}^-(D^2v)\leq 0&& ~~\mbox{in}~~Q_{r}(re_n,0)\backslash Q;\\
&v\leq C_0&& ~~\mbox{on}~~\partial Q\cap Q_{r}(re_n,0);\\
&v\leq 0&& ~~\mbox{on}~~\partial Q_{r}(re_n,0)\backslash \bar{Q}.
\end{aligned}\right.
\end{equation}
From the comparison principle and noting that $v(0)=0$, $v\geq 0$ for $t=0$,
\begin{equation}\label{e1.4}
  u\geq v\geq Cx_n~~~~\mbox{ on}~~~~\{x'=0,0< x_n< r/2, t=0\},
\end{equation}
where $C$ is universal.

By considering $v(x'-x'_0,x_n,t-t_0)$  for $(x'_0,t_0)\in S_{1/2}$ and similar arguments, we obtain
\begin{equation}\label{e2.1}
u\geq Cx_n~~~~\mbox{ in}~~\left\{(x'_0,0,t_0)\in S_{1/2},0<x_n<r/2\right\}.
\end{equation}
Finally, by the Harnack inequality again,
\begin{equation*}
u(x,t)\geq Cu(e_n/2,-3/4)\geq Cx_n,~~\forall ~(x,t)\in Q^+_{1/2}~~\mbox{and}~~x_n\geq 1/4.
\end{equation*}
By combining with \cref{e2.1}, \cref{e.Hopf} follows.\qed~\\

\begin{remark}\label{Hopf1}
Obviously, the construction of $v$ (see \cref{e1.3}) is key to the proof. Its construction is based on the fundamental solution for the heat equation and is analogous to the elliptic equations (see \cite{MR_Israel}). Here, we adopt the more concise form $v_0(x,t)=t^{-1}e^{-|x|^2/t}$, which is inspired by \cite[$\Gamma^N$ in P. 154]{MR1139064}.

The basic properties of $v_0$ are the following. The $v_0 > 0$ for $t>0$ and $v_0\equiv 0$ for $t=0$ except the singular point $0$. For any fixed $t$, $v_0$ is radial symmetry and decreasing in $|x|$. For any fixed $|x|$, $v_0$ is increasing for $0\leq t\leq |x|^2$ and decreasing afterwards. This indicates that the maximum of $v_0$ on $\partial Q_1(0,1)$ is attained at $\{t=1\}$. Finally, for $\beta$ large enough, $v_0^{\beta}$ is a subsolution.

Hence, by translating and stretching $v_0^{\beta}$ properly, $v(0)=0$ and \cref{e.Hopf.v} can be guaranteed.
\end{remark}

\begin{remark}\label{Hopf2}
Instead of proving \cref{e.Hopf} directly, we only need to prove the estimate on the line $\{x'=0,0< x_n< r/2,t=0\}$ (see \cref{e1.4}). The benefit is that we can construct the auxiliary function $v$ easily. The general estimate \cref{e.Hopf} can be obtained by translating $v$ properly.
\end{remark}

\begin{remark}\label{Hopf3}
In this proof, we can use $B_r(re_n)\times (-t_0,0]$ to replace $Q_r(re_n,0)$ for $0\leq t_0\leq r^2$, which guarantees the monotonicity of $v$ in $(-t_0,0)$ for any fixed $x\in \partial B_r(re_n)$.
\end{remark}
~\\

\begin{lemma}[Boundary Lipschitz regularity]\label{L4}
  Let $u$ be a viscosity solution of
  \begin{equation}\label{e.S0}
\left\{\begin{aligned}
&u\in S(\lambda,\Lambda,0)&& ~~\mbox{in}~~Q_1^+;\\
&u=0&& ~~\mbox{on}~~S_1.
\end{aligned}\right.
\end{equation}
Then
\begin{equation}\label{e.Lipschitz}
|u(x,t)|\leq Cx_n\|u\|_{L^{\infty}(Q_1^+)},~~~~\forall ~~(x,t) \in Q_{1/2}^+,
\end{equation}
where $C$ is universal.
\end{lemma}

\proof Without loss of generality, we assume that $\|u\|_{L^{\infty}(Q_1^+)}\leq 1$. Otherwise, we can consider $\tilde{u}=u/\|u\|_{L^{\infty}(Q_1^+)}$ and then $\|\tilde{u}\|_{L^{\infty}(Q_1^+)}=1$.

Let
\begin{equation*}
  v(x,t)=C\left(-\left(\frac{1}{t+1}e^{-\frac{|x+e_n|^2}{t+1}}\right)^{\beta}
+e^{-\beta}\right).
\end{equation*}
Then by taking $\beta$ and $C$ large enough, $v$ satisfies
\begin{equation}\label{e.Lip.v}
\left\{\begin{aligned}
&v_t-\mathcal{M}^-(D^2v)\geq 0&& ~~\mbox{in}~~Q^+_{1};\\
&v\geq 0&& ~~\mbox{on}~~S_{1};\\
&v\geq 1&&~~\mbox{in}~~\partial Q^+_{1}\backslash S_{1}.
\end{aligned}\right.
\end{equation}
According to the comparison principle and a direct calculation, we have (noting $v(0)=0$)
\begin{equation*}
-Cx_n\leq -v\leq u\leq v\leq Cx_n ~~~~\mbox{ on}~~\{x'=0,0\leq x_n\leq 1/2,t=0\}.
\end{equation*}
As in the proof of \Cref{L3}, by considering $v(x'-x'_0,x_n,t-t_0)$  for $(x'_0,t_0)\in S_{1/2}$ and similar arguments,
\begin{equation*}
-Cx_n\leq u\leq Cx_n ~~~~\mbox{ in}~~ Q^+_{1/2},
\end{equation*}
where $C$ is universal.
\qed~\\

\begin{remark}\label{Lip_re1}
Similar to the elliptic equations, the Hopf lemma and the boundary Lipschitz regularity can be obtained under the more general interior $C^{1,\mathrm{Dini}}$ condition and exterior $C^{1,\mathrm{Dini}}$ condition respectively \cite{lian2018boundary}.
\end{remark}

\begin{remark}\label{re2.2}
Since we place the singular point of the auxiliary function $v$ at $(-e_n,-1)$ (outside $Q_1^+$) rather than $(e_n/2,-1)$ (inside $Q_1^+$) as in \Cref{e1.3}, we consider the equation in $Q_1^+$ directly. In contrast, we have to subtract a set from $Q_1^+$ in \Cref{L3}. Hence, the proof of boundary Lipschitz regularity is shorter although the idea is similar to that of the Hopf lemma.
\end{remark}
~\\

With the aid of an iteration procedure, the Hopf lemma and the boundary Lipschitz regularity imply the boundary $C^{1,\alpha}$ regularity for the following model problem.
\begin{theorem}[Boundary $C^{1,\alpha}$ regularity on a flat boundary]\label{Th1C1aF}
  Let $u$ be a viscosity solution of
  \begin{equation}\label{model1}
\left\{\begin{aligned}
&u\in S(\lambda,\Lambda,0)&& ~~\mbox{in}~~Q_1^+;\\
&u=0&& ~~\mbox{on}~~S_1.
\end{aligned}\right.
\end{equation}
Then $u\in C^{1,\bar{\alpha}}(0)$, i.e., there exists a constant $a$ such that
\begin{equation}\label{e.fbC1a}
|u(x,t)-ax_n|\leq C_1 x_n|(x,t)|^{\bar{\alpha}}\|u\|_{L^{\infty}(Q_1^+)},~~\forall~(x,t)\in Q^+_{1/2}
\end{equation}
and
\begin{equation*}
  |a|\leq C_1,
\end{equation*}
where $0<\bar{\alpha}<1$ and $C_1$ are universal.
\end{theorem}

\begin{remark}\label{re.xa1}
Throughout this paper, $\bar{\alpha}$ and $C_1$ are fixed constants.
\end{remark}

\proof Assume that $\|u\|_{L^{\infty}(B_1^+)}\leq 1$. To show \cref{e.fbC1a}, we only need to prove the following: there exist a nonincreasing sequence $\{a_m\}$ ($m\geq 0$) and a nondecreasing sequence $\{b_m\}$ ($m\geq 0$) such that for all $m\geq 1$,
\begin{equation}\label{e.fbC1a.k}
\begin{aligned}
  &b_mx_n\leq u\leq a_mx_n~~~~\mbox{   in}~~~~Q^+_{2^{-m}},\\
  &0\leq a_m-b_m\leq (1-\mu)(a_{m-1}-b_{m-1}),\\
  &|a_0|\leq 2C~~\mbox{ and }~~|b_0|\leq 2C,
\end{aligned}
\end{equation}
where $0<\mu<1/2$ and $C$ are universal.

We prove the above by induction. From \Cref{L4},
\begin{equation*}
  -Cx_n\leq u \leq Cx_n~~~~\mbox{ in}~~Q^+_{1/2}.
\end{equation*}
Thus, by taking $a_1=C, b_1=-C$ and $a_0=2C, b_0=-2C$, \cref{e.fbC1a.k} holds for $m=1$.

Assume that \cref{e.fbC1a.k} holds for $m$ and we need to prove it for $m+1$. Let $r=2^{-m}$. Since \cref{e.fbC1a.k} holds for $m$, $u(re_n/2,-3r^2/4)$ is greater than  or less than $(a_m+b_m)/2\cdot r/2$. That is, there are two possible cases:
\begin{equation*}
\begin{aligned}
&\mbox{\textbf{Case 1}:}~~&&~~u(re_n/2,-3r^2/4)\geq \frac{a_m+b_m}{2}\cdot \frac{r}{2},\\
&\mbox{\textbf{Case 2}:}~~&&~~u(re_n/2,-3r^2/4)< \frac{a_m+b_m}{2}\cdot \frac{r}{2}.
\end{aligned}
\end{equation*}
Without loss of generality, we suppose that \textbf{Case 1} holds. Let $w=u-b_mx_n$ and then
\begin{equation*}
  \left\{\begin{aligned}
    &w\in S(\lambda,\Lambda,0)&& ~~\mbox{in}~~Q_{r}^+;\\
    &w\geq 0&& ~~\mbox{in}~~Q_{r}^+;\\
    &w(re_n/2,-3r^2/4)\geq \frac{a_m-b_m}{2}\cdot \frac{r}{2}.
  \end{aligned}\right.
\end{equation*}

From \Cref{L3}, for some universal constant $0<\mu<1$,
\begin{equation*}\label{e.fbC1a-w}
w(x)\geq \mu(a_m-b_m)x_n~~\mbox{ in}~~Q_{r/2}^+.
\end{equation*}
Thus,
\begin{equation*}\label{e.fbC1a-w2}
u(x)\geq (b_m+\mu(a_m-b_m))x_n~~\mbox{ in}~~Q_{r/2}^+.
\end{equation*}
Let $a_{m+1}=a_m$ and $b_{m+1}=b_m+\mu(a_m-b_m)$. Then
\begin{equation*}\label{e.reg-ak}
a_{m+1}-b_{m+1}=(1-\mu)(a_m-b_m).
\end{equation*}
Hence, \cref{e.fbC1a.k} holds for $m+1$. By induction, the proof is completed. \qed~\\

\begin{remark}\label{re.flatC1a.1}
The idea of the proof is to compare $u$ with the special solution ``$x_n$'', which has been used widely (see \cite[Theorem 2.1]{MR1139064}, \cite[Lemma 2.12]{MR_Israel} and \cite[Theorem 1.4]{MR_cones}).
\end{remark}

\begin{remark}\label{re2.3}
In fact, we can obtain $u\in C^{1,\bar\alpha}(x_0,t_0)$ for any $(x_0,t_0)\in S_{1/2}$ by almost the same proof. Hence, for any $p_0=(x_0,t_0)\in S_{1/2}$, there exists $a_{p_0}\in \mathbb{R}$ such that
\begin{equation*}
|u(x,t)-a_{p_0}x_n|\leq C_1 x_n|(x,t)-(x_0,t_0)|^{\bar{\alpha}}\|u\|_{L^{\infty}(Q_1^+)},~~\forall~(x,t)\in Q^+_{1/2}.
\end{equation*}
Then for any $p_0=(x_0,t_0),p_1=(x_1,t_1)\in S_{1/2}$, by the triangle inequality and eliminating $x_n$, we have
\begin{equation*}
|a_{p_0}-a_{p_1}|\leq C_1|(x_0,t_0)-(x_1,t_1)|^{\bar{\alpha}}\|u\|_{L^{\infty}(Q_1^+)}.
\end{equation*}
That is, $u_n$ exists and is $C^{\bar{\alpha}}$ on $S_{1/2}$.
\end{remark}
~\\

\subsection{Boundary pointwise $C^{1,\alpha}$ regularity on general boundaries}

We will attain the pointwise $C^{1,\alpha}$ regularity on a general boundary in this subsection. The idea of perturbation is employed in the proof. Precisely, a $C^{1,\alpha}$ estimate on some scale will be obtained first by the compactness method with small prescribed data. Then we use scaling technique and some proper normalizations to get the pointwise $C^{1,\alpha}$ regularity on the boundary.

The equicontinuity is necessary for the compactness method. The equicontinuity in the interior of the domain can be deduced from the Harnack inequality directly. In the following, we make some uniform estimate, which can be regarded as an estimate of modulus of continuity up to the boundary.

\begin{lemma}\label{L21}
Let $0<\delta<1/4$. Suppose that $u$ satisfies
\begin{equation*}
\left\{\begin{aligned}
&u\in S(\lambda,\Lambda,f)&& ~~\mbox{in}~~\Omega_1;\\
&u=g&& ~~\mbox{on}~~(\partial \Omega)_1
\end{aligned}\right.
\end{equation*}
with $\|u\|_{L^{\infty}(\Omega_1)}\leq 1$, $\|f\|_{L^{n+1}(\Omega_1)}\leq\delta$, $\|g\|_{L^{\infty}((\partial \Omega)_1)}\leq \delta$
and $\underset{Q_1}{\mathrm{osc}}~\partial\Omega \leq \delta$, where
\begin{equation*}
\underset{Q_{r}}{\mathrm{osc}}~\partial\Omega:= \underset{(x,t)\in (\partial\Omega)_r}{\sup} x_n -\underset{(x,t)\in (\partial\Omega)_r}{\inf} x_n~~\mbox{ for }~~0<r<1.
\end{equation*}

Then
\begin{equation}\label{e1.6}
 |u(x,t)|\leq C(x_n+\delta),~~~~\forall ~~(x,t) \in \Omega_{1/4},
\end{equation}
where $C$ is universal.
\end{lemma}

\proof Let $\tilde{B}^{+}_{1}=B^{+}_{1}-\delta e_n $, $\tilde{T}_1=T_1-\delta e_n$, $\tilde{Q}^{+}_{1}=Q^{+}_{1}-\delta (e_n,0) $ and $\tilde{S}_1=S_1-\delta (e_n,0)$. Then $(\partial \Omega)_{1/4}\subset \tilde{Q}^{+}_{1}$.

Set
\begin{equation*}
  v(x,t)=C\left(-\left(\frac{1}{t+1}e^{-\frac{|x+(1+\delta)e_n|^2}{t+1}}\right)^{\beta}
+e^{-\beta}\right).
\end{equation*}
Then by taking $\beta$ and $C$ large enough, $v$ satisfies $v(-\delta e_n,0)=0$ and
\begin{equation}\label{e.uniform.v}
\left\{\begin{aligned}
&v_t-\mathcal{M}^-(D^2v)\geq 0&& ~~\mbox{in}~~\tilde Q^+_{1};\\
&v\geq 0&& ~~\mbox{on}~~\tilde S_{1};\\
&v\geq 1&&~~\mbox{in}~~\partial \tilde Q^+_{1}\backslash \tilde S_{1}.
\end{aligned}\right.
\end{equation}
By the boundary Lipschitz regularity (\Cref{L4}),
\begin{equation*}
  |v(x,t)|\leq C(x_n+\delta)~~~~\mbox{ in}~~\tilde{Q}^+_{1/2}.
\end{equation*}

Next, let $w=u-v$ and then $w$ satisfies (noting $v\geq 0$ in $\tilde Q^+_{1}$)
\begin{equation*}
    \left\{
    \begin{aligned}
      &w\in \underline{S}(\lambda,\Lambda , f) &&\mbox{in}~~ \Omega \cap \tilde{Q}^{+}_{1}; \\
      &w\leq g &&\mbox{on}~~\partial \Omega \cap \tilde{Q}^{+}_{1};\\
      &w\leq 0 &&\mbox{on}~~\partial \tilde{Q}^{+}_{1}\cap \bar{\Omega}.
    \end{aligned}
    \right.
\end{equation*}
According to the Aleksandrov-Bakel’man-Pucci-Krylov-Tso maximum principle \cite[Theorem 3.14]{MR1135923},
\begin{equation*}
  \begin{aligned}
\sup_{\Omega \cap \tilde{Q}^{+}_{1}} w& \leq\|g\|_{L^{\infty }(\partial \Omega \cap \tilde{Q}^{+}_{1})}+C\|f\|_{L^{n+1}(\Omega \cap \tilde{Q}^{+}_{1})}\leq C\delta.
  \end{aligned}
\end{equation*}
Hence,
\begin{equation*}
\sup_{\Omega_{1/4}} u\leq \sup_{\Omega\cap \tilde{Q}_{1/2}^+} u \leq \|v\|_{L^{\infty}(\tilde{Q}_{1/2}^+)}+\sup_{\Omega \cap \tilde{Q}^{+}_{1}} w\leq C(x_n+\delta).
\end{equation*}

The proof for
\begin{equation*}
  \inf_{\Omega_{1/4}} u \geq -C(x_n+\delta)
\end{equation*}
is similar and we omit it here. Hence, the proof is completed.\qed~\\
\begin{remark}\label{re2.1}
Strictly speaking, \cref{e1.6} is not a modulus of continuity up to the boundary. In fact, since the boundary is not smooth, we can't obtain the true continuity up to the boundary. Nevertheless, it is enough for the compactness method (see \Cref{L.K1}), which has been inspired by \cite[Lemma 3.4]{MR3246039} and used in \cite[Lemma 2.4]{MR4088470}. The compactness method based on \cref{e1.6} allows us to treat non-smooth boundaries, especially for the pointwise regularity.
\end{remark}

To apply the compactness method, we also need the following closedness result for viscosity solutions, which was proved by Crandall, Kocan and \'{S}wi\polhk ech \cite{MR1789919}.

\begin{lemma}\label{L23}
Let $F$ and $\{F_m\}_{m\geq 1}$ be uniformly elliptic operators with $\lambda$ and $\Lambda$, $f$ and $f_m \in L^p(\Omega)$ ($p\geq n+1$), and $u_m$ be $L^p$-viscosity subsolutions (supersolutions) of
\begin{equation*}
  (u_m)_t-F_m(D^2u_m,x,t)=f_m~~~~\mbox{in}~~\Omega.
\end{equation*}
Suppose that $u_m\rightarrow u$ locally uniformly in $\Omega$ as $m\rightarrow \infty$. Moreover, assume that for any $Q_r(x_0,t_0)\subset \subset \Omega$ and $\varphi \in W^{2,1,p}(Q_r(x_0,t_0))$,
\begin{equation}\label{L23g}
  \|(\psi-\psi_m)^+\|_{L^p(Q_r(x_0,t_0))}\rightarrow 0 ~~\left(\|(\psi-\psi_m)^-\|_{L^p(Q_r(x_0,t_0))}\rightarrow 0\right),
\end{equation}
where
\begin{equation*}
  \psi_m(x,t)=-F_m(D^2\varphi,x,t)-f_m(x,t)~~\mbox{and}~~\psi(x,t)=-F(D^2\varphi,x,t)-f(x,t).
\end{equation*}

Then $u$ is an $L^p$-viscosity subsolution (supersolution) of
\begin{equation*}
u_t-F(D^2u,x,t)=f~~\mbox{ in}~~\Omega.
\end{equation*}
Furthermore, if $F$ and $f$ are continuous, then $u$ is a $C$-viscosity subsolution (supersolution) and we only require that \cref{L23g} holds for $\varphi \in C^2(Q_r(x_0,t_0))$.
\end{lemma}

\begin{remark}\label{re.closed.1}
In this lemma, a function $u\in W^{2,1,p}$ means that $u$ is second order weak differentiable with $x$ and first order differentiable with $t$ and all weak derivatives belong to $L^p$.
\end{remark}

%
%
%

Next, we prove that the solution can be approximated by a linear function on some scale provided that the prescribed data are small enough.

\begin{lemma}\label{L.K1}
For any $0<\alpha<\bar{\alpha}$, there exists $\delta>0$ depending only on $n, \lambda,\Lambda$ and $\alpha$ such that if $u$ satisfies
\begin{equation*}
\left\{\begin{aligned}
&u\in S(\lambda,\Lambda,f)&& ~~\mbox{in}~~\Omega_1;\\
&u=g&& ~~\mbox{on}~~(\partial \Omega)_1
\end{aligned}\right.
\end{equation*}
with $\|u\|_{L^{\infty}(\Omega_1)}\leq 1$, $\|f\|^{*}_{L^{n+1}(\Omega_1)}\leq \delta$, $\|g\|_{L^{\infty}((\partial \Omega)_1)}\leq \delta$
and $\underset{Q_1}{\mathrm{osc}}~\partial\Omega \leq \delta$, then there exists a constant $a$ such that
\begin{equation*}\label{e.lK1.1}
  \|u-ax_n\|_{L^{\infty}(\Omega_{\eta})}\leq \eta^{1+\alpha}
\end{equation*}
and
\begin{equation*}\label{e.lK1.2}
|a|\leq C_1,
\end{equation*}
where $0<\eta<1/8$ depends only on $n,\lambda,\Lambda$ and $\alpha$.
\end{lemma}

\proof We prove the lemma by contradiction. Suppose that the lemma is false. Then there exist $0<\alpha<\bar{\alpha}$ and sequences of $u_m,f_m,g_m,\Omega_m$ such that
\begin{equation*}
\left\{\begin{aligned}
&u_m\in S(\lambda,\Lambda,f_m)&& ~~\mbox{in}~~\Omega_m\cap Q_1;\\
&u_m=g_m&& ~~\mbox{on}~~\partial \Omega_m\cap Q_1
\end{aligned}\right.
\end{equation*}
with $\|u_m\|_{L^{\infty}(\Omega_m\cap Q_1)}\leq 1$,
$\|f_m\|^{*}_{L^{n+1}(\Omega_m\cap Q_1)}\leq 1/m$, $\|g_m\|_{L^{\infty}(\partial \Omega_m\cap Q_1)}\leq 1/m$
and $\underset{Q_1}{\mathrm{osc}}~\partial\Omega_m\leq 1/m$. Furthermore,
\begin{equation}\label{e.key1}
  \|u_m-ax_n\|_{L^{\infty}(\Omega_{m}\cap Q_{\eta})}> \eta^{1+\alpha}, \forall~|a|\leq C_1,
\end{equation}
where $0<\eta<1/8$ is taken small such that
\begin{equation}\label{e.key2}
C_1\eta^{\bar{\alpha}-\alpha}<1/2.
\end{equation}

Clearly, $u_m$ are uniformly bounded. In addition, by the interior H\"{o}lder regularity (see \cite[Theorem 4.19]{MR1135923}), $u_m$ are equicontinuous in any compact subset of $Q_1^+$. Hence, there exist a subsequence (denoted by $u_m$ again) and $\bar u$ such that $u_m\rightarrow \bar u$ locally uniformly in $Q_1^+$. Note also that $\|f_m\|_{L^{n+1}(\Omega_m\cap Q_1)}\rightarrow 0$. Hence, by \Cref{L23}, $\bar u\in C(Q_1^+)$ is a viscosity solution of
\begin{equation*}
\begin{aligned}
&\bar u\in S(\lambda,\Lambda,0)&& ~~\mbox{in}~~Q_{1/4}^+.
\end{aligned}
\end{equation*}

Next, by \Cref{L21},
\begin{equation*}
|u_m(x,t)|\leq C(x_n+1/m),~~~~\forall ~~(x,t) \in \Omega_m\cap Q_{1/4}.
\end{equation*}
For any $(x,t)\in Q^+_{1/4}$, by taking $m\rightarrow \infty$, we have
\begin{equation*}
  |\bar u(x,t)|\leq Cx_n.
\end{equation*}
Therefore, $u$ is continuous up to $S_{1/4}$ and $u\equiv 0$ on $S_{1/4}$.


By \Cref{Th1C1aF}, there exists a constant $\bar{a}$ such that
\begin{equation*}
  |\bar u(x,t)-\bar{a}x_n|\leq C_1 |(x,t)|^{1+\bar{\alpha}}, ~~\forall ~(x,t)\in Q_{1/8}^+
\end{equation*}
and
\begin{equation*}
  |\bar{a}|\leq C_1.
\end{equation*}
Hence, by noting \cref{e.key2},
\begin{equation}\label{e.key3}
  \|\bar u-\bar{a}x_n\|_{L^{\infty}(Q_{\eta}^+)}\leq \eta^{1+\alpha}/2.
\end{equation}
Set $a=\bar{a}$ and let $m\rightarrow \infty$ in \Cref{e.key1}. Then
\begin{equation*}
    \|\bar u-\bar{a}x_n\|_{L^{\infty}(Q_{\eta})}\geq \eta^{1+\alpha},
\end{equation*}
which contradicts with \cref{e.key3}.  ~\qed~\\

\begin{remark}\label{re.closed.1}
We point out here that there is an error in the proof of \cite[Lemma 2.8]{MR4088470}. The $\delta$ in that proof depends on $k$ and $k$ depends on $\delta$ in turn.
\end{remark}

Now, by a scaling argument, we can prove the boundary pointwise $C^{1,\alpha}$ regularity under rather strict conditions.
\begin{theorem}\label{L.Sca1}
Let $0<\alpha<\bar{\alpha}$ and $u$ be a viscosity solution of
\begin{equation*}
\left\{\begin{aligned}
&u\in S(\lambda,\Lambda,f)&& ~~\mbox{in}~~\Omega_1;\\
&u=g&& ~~\mbox{on}~~(\partial \Omega)_1.
\end{aligned}\right.
\end{equation*}
Assume that
\begin{equation}\label{e.Sca.C1a1}
  \begin{aligned}
    &\|u\|_{L^{\infty}(\Omega_1)}\leq 1;\\
    &\|f\|_{C^{-1,\alpha}(0)}\leq \delta;\\
    &|g(x,t)|\leq \frac{\delta}{2} |(x,t)|^{1+\alpha},~~~~\forall~~(x,t)\in (\partial\Omega)_1;\\
    &\|(\partial \Omega)_1\|_{C^{1,\alpha}(0)}\leq \frac{\delta}{2C_0},
  \end{aligned}
\end{equation}
where $\delta$ is as in \Cref{L.K1} and $C_0$ depends only on $n,\lambda,\Lambda$ and $\alpha$.

Then $u\in C^{1,\alpha}(0)$, i.e., there exists a constant $a$ such that
\begin{equation*}\label{e.Sca1.1}
  |u(x,t)-ax_n|\leq C|(x,t)|^{1+\alpha},~~\forall~~(x,t)\in \Omega_1
\end{equation*}
and
\begin{equation*}\label{e.Sca1.2}
|a|\leq C,
\end{equation*}
where $C$ depends only on $n,\lambda,\Lambda$ and $\alpha$.
\end{theorem}

\proof It is well-known that, to prove that $u$ is $C^{1,\alpha}$ at $0$, we only need to prove the following. There exists a sequence of $a_m$ ($m\geq -1$) such that for all $m\geq 0$,
\begin{equation}\label{e.ScaC1a.1}
\|u-a_mx_n\|_{L^{\infty }(\Omega_{\eta^{m}})}\leq \eta ^{m(1+\alpha )}
\end{equation}
and
\begin{equation}\label{e.ScaC1a.2}
|a_m-a_{m-1}|\leq C_1\eta ^{m\alpha},
\end{equation}
where $C_1$ is the universal constant as in \Cref{Th1C1aF} and $\eta$ depending only on $n,\lambda,\Lambda$ and $\alpha$, is as in \Cref{L.K1}.

We prove the above by induction. For $m=0$, by setting $a_0=a_{-1}\equiv 0$, \crefrange{e.ScaC1a.1}{e.ScaC1a.2} hold clearly. Suppose that the conclusion holds for $m=m_0$. We need to prove that the conclusion holds for $m=m_0+1$.

Let $r=\eta ^{m_0}$, $y=x/r$, $s=t/r^2$ and
\begin{equation}\label{e.Sca.C1a.v}
  v(y,s)=\frac{u(x,t)-a_{m_0}x_n}{r^{1+\alpha}}.
\end{equation}
Then $v$ satisfies
\begin{equation}\label{e.Sca.C1a.v2}
\left\{\begin{aligned}
&v\in S(\lambda,\Lambda,\tilde{f})&& ~~\mbox{in}~~\tilde{\Omega}\cap Q_1;\\
&v=\tilde{g}&& ~~\mbox{on}~~\partial \tilde{\Omega}\cap Q_1,
\end{aligned}\right.
\end{equation}
where
\begin{equation*}\label{e.Sca.C1a.v3}
  \begin{aligned}
&\tilde{f}(y,s)=\frac{f(x,t)}{r^{\alpha-1}},~~\tilde{g}(y,s)=\frac{g(x,t)-a_{m_0}x_n}{r^{1+\alpha}}
~~\mbox{ and }~~\tilde{\Omega}=\frac{\Omega}{r}.\\
  \end{aligned}
\end{equation*}

By \cref{e.ScaC1a.2}, there exists a constant $C_0>1$ depending only on $n,\lambda,\Lambda$ and $\alpha$ such that $|a_m|\leq C_0$ ($\forall~0\leq m\leq m_0$). In addition, from \cref{e.Sca.C1a1}, it is easy to verify that
\begin{equation*}
\begin{aligned}
&\|v\|_{L^{\infty}(\tilde{\Omega}\cap Q_1)}\leq 1,\\
&\|\tilde{f}\|^*_{L^{n+1}(\tilde{\Omega}\cap Q_{1})}=\frac{\|f\|^*_{L^{n+1}(\Omega_r)}}{r^{\alpha-1}}\leq \delta,\\
&\|\tilde{g}\|_{L^{\infty}(\partial \tilde \Omega \cap Q_1)}
\leq \frac{1}{r^{1+\alpha}}\left(\frac{\delta}{2}r^{1+\alpha}+C_0\cdot \frac{\delta}{2C_0}r^{1+\alpha}\right)
\leq \delta,\\
 &\underset{Q_1}{\mathrm{osc}}~\partial\tilde{\Omega}=
\frac{1}{r}\underset{Q_r}{\mathrm{osc}}~\partial\Omega \leq 2\|\partial \Omega\|_{C^{1,\alpha}(0)}r^{\alpha} \leq \delta.
\end{aligned}
\end{equation*}

By \Cref{L.K1}, there exists a constant $\bar{a}$ such that
\begin{equation*}
\begin{aligned}
    \|v-\bar{a}y_n\|_{L^{\infty }(\tilde{\Omega} _{\eta})}&\leq \eta ^{1+\alpha}
\end{aligned}
\end{equation*}
and
\begin{equation*}
|\bar{a}|\leq C_1.
\end{equation*}
Let $a_{m_0+1}=a_{m_0}+r^{\alpha}\bar{a}$. Then \cref{e.ScaC1a.2} hold for $m_0+1$. Recalling \cref{e.Sca.C1a.v}, we have
\begin{equation*}
  \begin{aligned}
&\|u-a_{m_0+1}x_n\|_{L^{\infty}(\Omega_{\eta^{m+1}})}\\
&= \|u-a_{m_0}x_n-r^{\alpha}\bar{a}x_n\|_{L^{\infty}(\Omega_{\eta r})}\\
&= \|r^{1+\alpha}v-r^{1+\alpha}\bar{a}y_n\|_{L^{\infty}(\tilde{\Omega}_{\eta})}\\
&\leq r^{1+\alpha}\eta^{1+\alpha}=\eta^{(m_0+1)(1+\alpha)}.
  \end{aligned}
\end{equation*}
Hence, \cref{e.ScaC1a.1} holds for $m=m_0+1$. By induction, the proof is completed.\qed
~\\

\noindent\textbf{Proof of \Cref{C1a}.} In fact, \Cref{L.Sca1} has already contained the essential ingredients for the boundary pointwise $C^{1,\alpha}$ regularity. The following proof is just a normalization.

Let $u_1(x,t)=u(x,t)-P_g(x,t)$. Thus $u_1$ satisfies
\begin{equation*}
\left\{\begin{aligned}
&u_1\in S(\lambda,\Lambda,f)&& ~~\mbox{in}~~\Omega_1;\\
&u_1=g_1&& ~~\mbox{on}~~(\partial \Omega)_1,
\end{aligned}\right.
\end{equation*}
where $g_1(x,t)=g(x,t)-P_g(x,t)$. Hence,
\begin{equation*}
  |g_1(x,t)|\leq [g]_{C^{1,\alpha}(0)}|(x,t)|^{1+\alpha}, ~~\forall ~(x,t)\in (\partial \Omega)_1.
\end{equation*}

Next, let $u_2=u_1/U$ where $U=\|u\|_{L^{\infty}(\Omega_1)}+\|g\|_{C^{1,\alpha}(0)}
+\|f\|_{C^{-1,\alpha}(0)}$. Then we have
\begin{equation*}
\left\{\begin{aligned}
&u_2\in S(\lambda,\Lambda,f_1)&& ~~\mbox{in}~~\Omega_1;\\
&u_2=g_2&& ~~\mbox{on}~~(\partial \Omega)_1,
\end{aligned}\right.
\end{equation*}
where $f_1=f/U$ and $g_2=g_1/U$. Hence,
\begin{equation*}
  \|u_2\|_{L^{\infty}(\Omega_1)}\leq 1,~~~~\|f_1\|_{C^{-1,\alpha}(0)}\leq 1
\end{equation*}
and
\begin{equation*}
  |g_2(x,t)|\leq  |(x,t)|^{1+\alpha}, ~~\forall ~(x,t)\in (\partial \Omega)_1.
\end{equation*}

Finally, let $y=x/\rho$, $s=t/\rho^2$ and $u_3(y,s)=u_2(x,t)$. Then,
\begin{equation*}
\left\{\begin{aligned}
&u_3\in S(\lambda,\Lambda,f_2)&& ~~\mbox{in}~~\tilde \Omega_1;\\
&u_3=g_3&& ~~\mbox{on}~~(\partial \tilde \Omega)_1,
\end{aligned}\right.
\end{equation*}
where $f_2(y,s)=\rho^2f_1(x,t)$, $g_3(y,s)=g_2(x,t)$ and $\tilde{\Omega}=\Omega/\rho$.

By a direct calculation,
\begin{equation*}
  \begin{aligned}
    &\|f_2\|^*_{L^{n+1}(\tilde{\Omega}_r)}=\rho^2\|f_1\|^*_{L^{n+1}(\Omega_{\rho r})}
    \leq \rho^{1+\alpha} r^{-1+\alpha},~~\forall ~0<r<1,\\
    &|g_3(y,s)|\leq \rho^{1+\alpha}|(y,s)|^{1+\alpha},~~\forall~~(y,s)\in(\partial \tilde\Omega)_1.
  \end{aligned}
\end{equation*}
In addition, since $\partial \Omega \in C^{1,\alpha}(0)$,
\begin{equation*}
  |x_n|\leq [\partial \Omega]_{C^{1,\alpha}(0)}|(x',t)|^{1+\alpha},~~\forall~~(x,t)\in(\partial \Omega)_1.
\end{equation*}
Hence,
\begin{equation*}
  \begin{aligned}
    &|y_n|\leq \rho^{\alpha}[\partial \Omega]_{C^{1,\alpha}(0)}|(y',s)|^{1+\alpha},
    ~~~~\forall~~(y,s)\in (\partial\tilde{\Omega})_1.
  \end{aligned}
\end{equation*}

Then by taking $\rho$ small enough (depending only on $n,\lambda,\Lambda,\alpha$ and $\|\partial \Omega\|_{C^{1,\alpha}(0)}$), \cref{e.Sca.C1a1} can be guaranteed. By \Cref{L.Sca1}, $u_3 \in C^{1,\alpha}(0)$ and therefore, $u\in C^{1,\alpha}(0)$ and the estimate \cref{e.C1a} holds. Furthermore, from above arguments, the $P$ in \cref{e.C1a} has the following form:
\begin{equation*}
P(x)=P_g(x)+\frac{aU}{\rho} x_n,
\end{equation*}
where $a$ is obtained by applying \Cref{L.Sca1} to $u_3$. Hence, $P(x',0)\equiv P_g(x',0)$. That is, \cref{e1.7} holds. \qed~\\

\section{Boundary $C^{2,\alpha}$ regularity}\label{sec:C2a}

In this section, we prove the boundary pointwise $C^{2,\alpha}$ regularity. Firstly, as \Cref{sec:C1a}, we study the pointwise $C^{2,\alpha}$ regularity on a flat boundary. Then we attack the boundary $C^{2,\alpha}$ regularity by the perturbation skill. An interesting observation is that if $u(0)=0$ and $Du(0)=0$, the boundary  $C^{2,\alpha}$ regularity holds even if $\partial \Omega\in C^{1,\alpha}(0)$ (originated from \cite{MR4088470}). There are similar observations for higher order regularity (see \Cref{Sec:Cka_freeb}) and that is the point to obtain the regularity of free boundaries.

Based on the $C^{1,\alpha}$ regularity, we give the proof of the $C^{2,\alpha}$ regularity on a flat boundary (i.e. the model problem) by the method of difference quotient.
\begin{theorem}[Boundary $C^{2,\alpha}$ regularity on a flat boundary]\label{Th1C2aF}
  Let $u$ be a viscosity solution of
  \begin{equation}\label{e.C2a.F}
\left\{\begin{aligned}
&u_t-F(D^2u)=0&& ~~\mbox{in}~~Q_1^+;\\
&u=0&& ~~\mbox{on}~~S_1.
\end{aligned}\right.
\end{equation}
Then $u\in C^{2,\alpha}(0)$ for any $0<\alpha<\bar{\alpha}$, i.e., there exists $P\in\mathcal{P}_2$ such that
\begin{equation}\label{e.C2a.F2}
|u(x,t)-P(x)|\leq C_2 x_n|(x,t)|^{1+\alpha}\left(\|u\|_{L^{\infty}(Q_1^+)}+|F(0)|\right),~~\forall~(x,t)\in Q^+_{1},
\end{equation}
\begin{equation*}
  \|P\|\leq C_2
\end{equation*}
and
\begin{equation}\label{e.C2a.F2.2}
F(D^2P)=0,
\end{equation}
where $C_2$ depends only on $n, \lambda$, $\Lambda$ and $\alpha$. The $P$ can be written as
\begin{equation*}
P(x)=ax_n-b_{in}x_ix_n,
\end{equation*}
where $a,b_{in}$ are constants.

Moreover, $u_t \in C^{1,\bar{\alpha}}(0)$, i.e., there exists a constant $b$ such that
\begin{equation*}
  |u_t(x,t)-bx_n|\leq Cx_n|(x,t)|^{\bar{\alpha}}\left(\|u\|_{L^{\infty}(Q_1^+)}+|F(0)|\right),~\forall ~(x,t)\in Q_{1}^+,
\end{equation*}
where $C$ is universal.
\end{theorem}
\proof Without loss of generality, we assume that
\begin{equation*}
\|u\|_{L^{\infty}(Q_1^+)}+|F(0)|\leq 1.
\end{equation*}
Since $u\in S(\lambda,\Lambda,-F(0))$, according to \Cref{C1a} and the interior $C^{1,\alpha}$ regularity (see \cite[Theorem 4.8]{MR1139064}), $u\in C^{1,\tilde\alpha}(Q_1^+ \cup S_1)$ for some small $\tilde{\alpha}<\bar{\alpha}$. In addition, note that if $u$ and $v$ are two viscosity solutions of \cref{e.C2a.F}, then $u-v\in S(\lambda,\Lambda,0)$ \cite[Theorem 4.6]{MR1139064}. Hence, by taking difference quotient in the direction $e_i$ ($1\leq i\leq n-1$), i.e.,
\begin{equation*}
  \Delta^{h}_{i} u (x,t)=\frac{u(x+he_i,t)-u(x,t)}{h}, ~~0<h<1/4,
\end{equation*}
we have
\begin{equation*}
\left\{\begin{aligned}
&\Delta^{h}_{i} u\in S(\lambda,\Lambda,0)&& ~~\mbox{in}~~Q_{3/4}^+;\\
&\Delta^{h}_{i} u=0&& ~~\mbox{on}~~S_{3/4}.
\end{aligned}\right.
\end{equation*}
By the closedness of viscosity solutions and sending $h\rightarrow 0$, we derive
\begin{equation*}
\left\{\begin{aligned}
&u_{i}\in S(\lambda,\Lambda,0)&& ~~\mbox{in}~~Q_{3/4}^+;\\
&u_{i}=0&& ~~\mbox{on}~~S_{3/4}.
\end{aligned}\right.
\end{equation*}
From \Cref{Th1C1aF}, $u_i\in C^{1,\bar{\alpha}}(0)$. That is, there exist constants $b_{in}$ such that
\begin{equation}\label{e.C2a.F3}
|u_i(x,t)-b_{in}x_n|\leq Cx_n|(x,t)|^{\bar{\alpha}},~\forall ~(x,t)\in Q_{1/2}^+,
\end{equation}
where $C$ is universal. We write $u_{in}(0)=b_{in}$ in the following proof.

Next, by the interior regularity, $u_t\in C^{\tilde\alpha}(Q_1^+)$ (see \cite[Theorem 4.9]{MR1139064}). Then by considering the difference quotient along direction $t$, i.e.,
\begin{equation*}
  \Delta^{h}_{t} u(x,t)=\frac{u(x,t+h)-u(x,t)}{h}
\end{equation*}
and taking the limit,
\begin{equation*}
\left\{\begin{aligned}
&u_{t}\in S(\lambda,\Lambda,0)&& ~~\mbox{in}~~Q_{3/4}^+;\\
&u_{t}=0&& ~~\mbox{on}~~S_{3/4}.
\end{aligned}\right.
\end{equation*}
Similarly, by \Cref{Th1C1aF}, $u_t \in C^{1,\bar{\alpha}}(0)$, i.e., there exists a constant $b$
\begin{equation}\label{e.C2a.F3.1}
|u_t(x,t)-bx_n|\leq Cx_n|(x,t)|^{\bar{\alpha}},~\forall ~(x,t)\in Q_{1/2}^+,
\end{equation}
where $C$ is universal. In what follows, we write $u_{tn}(0)=b$.

In the following, we show that $u_n\in C^{1,\bar\alpha}(S_1)$. Let $h,\tilde{h}, \hat{h}<1/2$ be three positive constants. Then by \cref{e.C2a.F3} and \cref{e.C2a.F3.1},
\begin{equation}\label{e.C2a.F4}
  \begin{aligned}
u(he_n&+\tilde{h}e_i, \hat{h})-u(he_n,0)-\left(u(\tilde{h}e_i, \hat{h})-u(0)\right)\\
=&u(he_n+\tilde{h}e_i, \hat{h})-u(he_n,0)\\
=&\tilde{h}u_i(he_n+\xi \tilde{h}e_i, \xi \hat{h})
+\hat{h} u_t(he_n+\xi \tilde{h}e_i, \xi \hat{h})\\
= & \tilde{h}\left(u_{in}(0)h+O (h(h^{\bar{\alpha}}+\tilde{h}^{\bar{\alpha}}+\hat{h}^{\frac{\bar{\alpha}}{2}}))\right)
+\hat{h}\left(u_{tn}(0) h+O (h(h^{\bar{\alpha}}+\tilde{h}^{\bar{\alpha}}+\hat{h}^{\frac{\bar{\alpha}}{2}}))\right)\\
= &u_{in}(0)h\tilde{h}+u_{tn}(0)h\hat{h}+
O (h(\tilde{h}+\hat{h})(h^{\bar{\alpha}}+\tilde{h}^{\bar{\alpha}}+\hat{h}^{\frac{\bar{\alpha}}{2}})).\\
  \end{aligned}
\end{equation}
In addition, by noting that $u_n\in C^{\tilde\alpha}(Q_1^+\cup S_1)$,
\begin{equation}\label{e.C2a.F5}
  \begin{aligned}
u(he_n&+\tilde{h}e_i,\hat{h})-u(\tilde{h}e_i,\hat{h})-\left(u(he_n,0)-u(0)\right)\\
=&h\left(u_n(\xi_1 he_n+\tilde{h}e_i,\hat{h})-u_n(\xi_2 he_n,0)\right)\\
=&h\left(u_n(\tilde{h}e_i,\hat{h})-u_n(0)+O(h^{\tilde \alpha})\right).\\
  \end{aligned}
\end{equation}
By comparing \cref{e.C2a.F4} with \cref{e.C2a.F5}, we have
\begin{equation*}
  \begin{aligned}
    &|u_n(\tilde{h}e_i,\hat{h})-u_n(0)-u_{in}(0)\tilde{h}-u_{tn}(0)\hat{h}|\\
\leq&  C\left(\tilde{h}^{1+\bar{\alpha}}+\tilde{h}\hat{h}^{\frac{\bar{\alpha}}{2}}\right)
+C\left(\hat{h}\tilde{h}^{\bar{\alpha}}+\hat{h}^{1+\frac{\bar{\alpha}}{2}}\right)
+C\tilde{h}h^{\bar{\alpha}}+C\hat{h}h^{\bar{\alpha}}+Ch^{\tilde \alpha}.
  \end{aligned}
\end{equation*}
Let $h\rightarrow 0$ and we have
\begin{equation*}
  \begin{aligned}
    &|u_n(\tilde{h}e_i,\hat{h})-u_n(0)-u_{in}(0)\tilde{h}-u_{tn}(0)\hat{h}|\\
\leq  &C\left(\tilde{h}^{1+\bar{\alpha}}+\tilde{h}\hat{h}^{\frac{\bar{\alpha}}{2}}\right)
+C\left(\hat{h}\tilde{h}^{\bar{\alpha}}+\hat{h}^{1+\frac{\bar{\alpha}}{2}}\right).
  \end{aligned}
\end{equation*}
Thus,
\begin{equation*}
  |u_n(\tilde{h}e_i,\hat{h})-u_n(0)-u_{in}(0)\tilde{h}|
  \leq C\left(\tilde{h}^{1+\bar{\alpha}}+\hat{h}^{\frac{1+\bar{\alpha}}{2}}\right).
\end{equation*}
That is, $u_n\in C^{1,\bar{\alpha}}(0)$ (regard $u_n$ as a function on $S_1$). Hence, $u_n\in C^{1,\bar{\alpha}}(S_{1})$ and $u_{ni}=u_{in}$ on $S_{1}$.

Similar to $u_i$ ($1\leq i\leq n-1$),
\begin{equation*}
\begin{aligned}
&u_{n}\in S(\lambda,\Lambda,0)&& ~~\mbox{in}~~Q_{3/4}^+.
\end{aligned}
\end{equation*}
Then by the boundary $C^{1,\alpha}$ regularity (\Cref{C1a}), we have $u_n\in C^{1,\alpha}(x_0)$ for any $0<\alpha<\bar{\alpha}$ and any $x_0\in S_{1}$. Hence,
\begin{equation*}
  \begin{aligned}
    u(x,t)=&\int_{0}^{1}(u_i(\tau x,\tau t)x_i+u_n(\tau x,\tau t)x_n+u_t(\tau x,\tau t)t)d\tau \\
=&\int_{0}^{1}(\tau u_{in}(0)x_n+\tau ^{1+\alpha}x_nO(|(x,t)|^{\alpha}))x_id\tau\\
    +&\int_{0}^{1}(u_n(0)+\tau u_{ni}(0)x_i
    +\tau u_{nn}(0)x_n+\tau u_{nt}(0)t
    +\tau^{1+\alpha}O(|(x,t)|^{1+\alpha}))x_nd\tau\\
    +&\int_{0}^{1}(\tau u_{tn}(0)x_n+\tau^{1+\alpha}x_nO(|(x,t)|^{\alpha}))t d\tau\\
=&u_n(0)x_n+u_{in}(0)x_ix_n+\frac{1}{2}u_{nn}(0)x_n^2+u_{tn}(0)tx_n+x_nO(|x|^{1+\alpha})\\
=&u_n(0)x_n+u_{in}(0)x_ix_n+\frac{1}{2}u_{nn}(0)x_n^2+x_nO(|x|^{1+\alpha}).
  \end{aligned}
\end{equation*}
That is, $u\in C^{2,\alpha}(0)$ and \cref{e.C2a.F2} holds. Finally, by \cref{e.C2a.F2} and the definition of viscosity solution, \cref{e.C2a.F2.2} holds (see \cite[Lemma 2.5]{MR1351007}).
~\qed~\\

\begin{remark}\label{re.C2a.F3}
From the viewpoint of scaling, $t\simeq x^2$. For instance, $u\in C^{1,\alpha}$ means that $u$ is $C^{1,\alpha}$ with respect to $x$ and only $C^{(1+\alpha)/2}$ with respect to $t$. However, for \cref{e.C2a.F}, we have higher regularity with respect to $t$ than that of $x$ since $u_t$ and $u_i$ satisfy the same equation $S(\lambda,\Lambda,0)$.
\end{remark}

\begin{remark}\label{re.C2a.F1}
Note that $u$ is approximated by a quadratic polynomial in the form $ax_n+b_{in}x_ix_n$, which is the same as the boundary $C^{2,\alpha}$ regularity on a flat boundary for elliptic equations. The reason is that $u_i(0)=u_t(0)=u_{ij}(0)=0$ for any $1\leq i,j<n$.
\end{remark}

\begin{remark}\label{re.C2a.F2}
For fully nonlinear elliptic equations, Silvestre and Sirakov \cite[Lemma 4.1]{MR3246039} proved the boundary $C^{2,\alpha}$ regularity on a flat boundary. They first deduced $u\in C^{2,\alpha}(T_1)$ by comparing two equalities similar to \cref{e.C2a.F4} and \cref{e.C2a.F5}. Next, the only unknown second order derivative $u_{nn}(0)$ was determined by the uniform ellipticity of the operator $F$. Then they obtained that $u\in C^{2,\alpha}(0)$ along $e_n$ based on the definition of viscosity solution and the Hopf lemma. Finally, they completed the proof by combining above results together. Adimurthi, Banerjee and Verma \cite[Lemma 3.1]{MR4073515} extended the result to fully nonlinear parabolic equations by an analogous method. Our proof is much simpler by employing the boundary $C^{1,\alpha}$ regularity.
\end{remark}

\begin{remark}\label{re.C2a.F4}
The exponent for the $C^{2,\alpha}$ regularity is strictly less than $\bar{\alpha}$. Note that in \cite{MR3246039}, these two exponents can be the same.
\end{remark}

\begin{remark}\label{re3.1}
Deriving $u_n\in C^{1,\alpha}(S_1)$ by \cref{e.C2a.F4} and \cref{e.C2a.F5} is inspired by \cite[Lemma 4.1]{MR3246039}. The idea is similar to that proving $f_{ij}=f_{ji}$ in mathematical analysis.
\end{remark}
~\\

In the following proof, we will use a kind of homogeneous polynomial in a special form. We call $H\in \mathcal{HP}_k$ a $k$-form ($k\geq 1$) if $H$ can be written as
\begin{equation}\label{k.form}
H(x,t)=\sum_{|\sigma|+2\gamma=k,\sigma_n\geq 1}\frac{a_{\sigma\gamma}}{\sigma!\gamma!}x^{\sigma}t^{\gamma}.
\end{equation}
That is, $x_n$ appears in the expression of $H$ at least once (thus $H\equiv 0$ on $S_1$), which turns out to be vital for the boundary regularity.~\\

On the basis of the $C^{2,\alpha}$ regularity on a flat boundary, we next obtain the boundary $C^{2,\alpha}$ regularity on a general boundary. Similar to \Cref{L.K1}, the following lemma is an estimate of $u$ on some scale when the prescribed data are small.

\begin{lemma}\label{L-C2a.K1}
Suppose that $F$ satisfies \cref{e.C2a.beta} at $0$. For any $0<\alpha<\bar{\alpha}$, there exists $\delta>0$ depending only on $n,\lambda,\Lambda$ and $\alpha$ such that if $u$ satisfies
\begin{equation*}
\left\{\begin{aligned}
&u_t-F(D^2u,x,t)=f&& ~~\mbox{in}~~\Omega_1;\\
&u=g&& ~~\mbox{on}~~(\partial \Omega)_1
\end{aligned}\right.
\end{equation*}
with $\|u\|_{L^{\infty}(\Omega_1)}+|F_0(0)|\leq 1$, $u(0)=|Du(0)|=0$, $\|\beta\|^*_{L^{n+1}(\Omega_1)}\leq \delta$, $\|\gamma\|_{C^{-1,\alpha}(0)}\leq \delta$, $\|f\|_{C^{-1,\alpha}(0)}\leq \delta$, $\|g\|_{C^{1,\alpha}(0)}\leq \delta$ and $\|(\partial \Omega)_1\|_{C^{1,\alpha}(0)} \leq \delta$, then there exists a $2$-form $H$ such that
\begin{equation*}
  \|u-H\|_{L^{\infty}(\Omega_{\eta})}\leq \eta^{2+\alpha},
\end{equation*}
\begin{equation*}
\|H\|\leq C_2+1
\end{equation*}
and
\begin{equation*}
F_0(D^2  H)=0,
\end{equation*}
where $\eta$ depends only on $n,\lambda,\Lambda$ and $\alpha$.
\end{lemma}

\proof As the proof of \Cref{L.K1}, we prove the lemma by contradiction. Suppose that the lemma is false. Then there exist $0<\alpha<\bar{\alpha}$ and sequences of $F_m,u_m,f_m,g_m,\Omega_m$ satisfying
\begin{equation*}
\left\{\begin{aligned}
&(u_m)_t-F_m(D^2u_m,x,t)=f_m&& ~~\mbox{in}~~\Omega_m\cap Q_1;\\
&u_m=g_m&& ~~\mbox{on}~~\partial \Omega_m\cap Q_1.
\end{aligned}\right.
\end{equation*}
In addition, $F_m$ satisfies \cref{e.C2a.beta} with $F_{0m},\beta_m$ and $\gamma_m$. Furthermore, $\|u_m\|_{L^{\infty}(\Omega_m\cap Q_1)}+|F_{0m}(0)|\leq 1$, $u_m(0)=|Du_m(0)|=0$, $\|\beta_m\|^*_{L^{n+1}(\Omega_m\cap Q_1)}\leq 1/m$, $\|\gamma_m\|_{C^{-1,\alpha}(0)}\leq 1/m$, $\|f_m\|_{C^{-1,\alpha}(0)}\leq 1/m$, $\|g_m\|_{C^{1,\alpha}(0)}\leq 1/m$ and $\|\partial \Omega_m\cap Q_1\|_{C^{1,\alpha}(0)} \leq 1/m$ and for any $2$-form $H$ satisfying $\|H\|\leq C_2+1$ and
\begin{equation}\label{e.lC2a.K1}
F_{0m}(D^2 H)=0,
\end{equation}
we have
\begin{equation}\label{e.lC2a.K2}
  \|u_m-H\|_{L^{\infty}(\Omega_{m}\cap Q_{\eta})}> \eta^{2+\alpha},
\end{equation}
where $0<\eta<1$ is taken small such that
\begin{equation}\label{e.lC2a.K3}
C_2\eta^{(\bar{\alpha}-\alpha)/2}<1/2.
\end{equation}

As in the proof of \Cref{L.K1}, $u_m$ are uniformly bounded and equicontinuous in any compact subset of $Q_1^+$. Thus, there exists a subsequence (denoted by $u_m$ again) and $\bar u: Q_1^+\rightarrow \mathbb R$ such that $u_m\rightarrow \bar u$ uniformly in any compact subset of $Q^+_1$. In addition, since $|F_{0m}(0)|\leq 1$ and $F_{0m}$ are Lipschitz continuous in $M$ with a uniform Lipschitz constant, there exists a subsequence (denoted by $F_{0m}$ again) and $\bar F:\mathcal {S}^{n}\rightarrow \mathbb R$ such that $F_{0m}\rightarrow \bar F$ uniformly in any compact subset of $\mathcal {S}^{n}$.

Furthermore, for any $Q\subset\subset Q_1^+$ and $\varphi\in C^2(\bar{Q})$, set $\psi_m(x,t)=-F_m(D^2\varphi,x,t)-f_m$ and $\psi(x,t)=-\bar F(D^2\varphi)$. Then
\begin{equation*}
  \begin{aligned}
\|\psi&_m-\psi\|_{L^{n+1}(Q)}\\
\leq& \|F_m(D^2\varphi,x,t)-F_{0m}(D^2\varphi)\|_{L^{n+1}(Q)}\\
    &+\|F_{0m}(D^2\varphi)-\bar F(D^2\varphi)\|_{L^{n+1}(Q)}+\|f_m\|_{L^{n+1}(Q)}\\
\leq& \|\beta_{m}\|_{L^{n+1}(Q)}\|D^2\varphi\|_{L^{\infty}(Q)}
    +\|\gamma_{m}\|_{L^{n+1}(Q)}\\
    &+\|F_{0m}(D^2\varphi)-\bar F(D^2\varphi)\|_{L^{\infty}(Q)}+\|f_m\|_{L^{n+1}(Q)}.
\end{aligned}
\end{equation*}
Thus, $\|\psi_m-\psi\|_{L^{n+1}(Q)}\rightarrow 0$ as $m\rightarrow \infty$. Then by combining \Cref{L21} with \Cref{L23}, $\bar u$ is a viscosity solution of
\begin{equation*}
\left\{\begin{aligned}
&\bar u_t-\bar F(D^2\bar u)=0&& ~~\mbox{in}~~Q_{1}^+;\\
&\bar u=0&& ~~\mbox{on}~~S_{1}.
\end{aligned}\right.
\end{equation*}

Since $u_m\in S(\lambda,\Lambda,f_m-F_m(0,x,t))$ and
\begin{equation*}
|F_m(0,x,t)-F_{0m}(0)|\leq \gamma_m(x,t),~~\forall ~(x,t)\in \Omega_m\cap Q_1,
\end{equation*}
by the $C^{1,\alpha}$ regularity (\Cref{C1a}) for $u_m$ and noting $u_m(0)=|Du_m(0)|=0$, we have
\begin{equation*}
\|u_m\|_{L^{\infty }(\Omega_m\cap Q_r)}\leq Cr^{1+\alpha}, ~~~~\forall~0<r<1,
\end{equation*}
where $C$ depends only on $n,\lambda,\Lambda$ and $\alpha$. Since $u_m$ converges to $\bar u$ uniformly,
\begin{equation*}
\|\bar u\|_{L^{\infty }(Q_r^+)}\leq Cr^{1+\alpha}, ~~~~\forall~0<r<1.
\end{equation*}
Hence, $\bar u(0)=|D\bar u(0)|=0$.

By \Cref{Th1C2aF}, there exists a $2$-form $\bar{H}$ such that
\begin{equation}\label{e.lC2a.K4}
  |\bar u(x,t)-\bar H|\leq C_2 x_n|(x,t)|^{2+(\alpha+\bar{\alpha})/2}, ~~\forall ~(x,t)\in Q_{1/2}^+,
\end{equation}
\begin{equation*}
\|\bar H\|\leq C_2
\end{equation*}
and
\begin{equation*}
  \bar F(D^2 \bar H)=0.
\end{equation*}
Combining \cref{e.lC2a.K3} and \cref{e.lC2a.K4}, we have
\begin{equation}\label{e.lC2a.K5}
  \|\bar u-\bar H\|_{L^{\infty}(Q_{\eta}^+)}\leq \eta^{2+\alpha}/2.
\end{equation}

In addition, since $F_{0m}(D^2 \bar H)\rightarrow \bar F(D^2 \bar H)=0$, there exist $|\tau_m|\leq \eta^{2+\alpha}/4$ and $\tau_m\rightarrow 0$ such that
\begin{equation*}
  \begin{aligned}
F_{0m}(D^2 \bar H+\tau_m \delta_{nn})=0,
  \end{aligned}
\end{equation*}
where $\delta_{nn}$ denotes the matrix $A_{ij}$ whose elements are all $0$ except $A_{nn}=1$ (similarly hereinafter).

Hence, \cref{e.lC2a.K2} holds for $\bar H+\tau_m x^2_n/2$, i.e.,
\begin{equation*}
  \|u_m-\bar H-\tau_m x_n^2/2\|_{L^{\infty}(\Omega_m\cap Q_{\eta})}> \eta^{2+\alpha}.
\end{equation*}
Let $m\rightarrow \infty$, we have
\begin{equation*}
    \|\bar u-\bar H\|_{L^{\infty}(Q_{\eta}^+)}\geq \eta^{2+\alpha},
\end{equation*}
which contradicts with \cref{e.lC2a.K5}.  ~\qed~\\

Next, we get a special $C^{2,\alpha}$ estimate (also the essential) under rather strict conditions.
\begin{theorem}\label{t-C2a.S}
Let $0<\alpha <\bar{\alpha}$ and $F$ satisfy \cref{e.C2a.beta} at $0$. Suppose that $u$ satisfies
\begin{equation*}
\left\{\begin{aligned}
&u_t-F(D^2u,x,t)=f&& ~~\mbox{in}~~\Omega_1;\\
&u=g&& ~~\mbox{on}~~(\partial \Omega)_1.
\end{aligned}\right.
\end{equation*}
Assume that
\begin{equation}\label{e.C2a.S.be}
\begin{aligned}
&\|u\|_{L^{\infty}(\Omega_1)}+|F_0(0)|\leq 1,~ u(0)=|Du(0)|=0, \\
&|\beta(x,t)|\leq \frac{\delta}{2C_0} |(x,t)|^{\alpha}, ~~\forall ~(x,t)\in \Omega_1,\\
&|\gamma(x,t)|\leq \frac{\delta}{2} |(x,t)|^{\alpha}, ~~\forall ~(x,t)\in \Omega_1,\\
&|f(x,t)|\leq \frac{\delta}{2} |(x,t)|^{\alpha}, ~~\forall ~(x,t)\in \Omega_1,\\
&|g(x,t)|\leq \frac{\delta}{2}|(x,t)|^{2+\alpha}, ~~\forall ~(x,t)\in (\partial \Omega)_1,\\
&\|\partial \Omega\cap Q_1\|_{C^{1,\alpha}(0)} \leq \frac{\delta}{2C_0},\\
\end{aligned}
\end{equation}
where $\delta$ is as in \Cref{L-C2a.K1} and $C_0$ depends only on $n,\lambda,\Lambda$ and $\alpha$.

Then $u\in C^{2,\alpha}(0)$, i.e., there exists a $2$-form $H$ such that
\begin{equation}\label{e.C2a.S1}
  |u(x,t)-H(x)|\leq C |(x,t)|^{2+\alpha}, ~~\forall ~(x,t)\in \Omega_{1},
\end{equation}
\begin{equation}\label{e.C2a.S3}
\|H\| \leq C
\end{equation}
and
\begin{equation}\label{e.C2a.S2}
F_0(D^2 H)=0,
\end{equation}
where $C$ depends only on $n, \lambda, \Lambda$ and $\alpha$.
\end{theorem}

\proof To prove that $u$ is $C^{2,\alpha}$ at $0$, we only need to prove the following. There exists a sequence of $2$-forms $H_m$ ($m\geq 0$) such that for all $m\geq 1$,
\begin{equation}\label{e.C2a.S4}
\|u-H_m\|_{L^{\infty }(\Omega_{\eta^{m}})}\leq \eta ^{m(2+\alpha )},
\end{equation}
\begin{equation}\label{e.C2a.S6}
\|H_m-H_{m-1}\|\leq (C_2+1)\eta ^{m\alpha},
\end{equation}
and
\begin{equation}\label{e.C2a.S5}
F_0(D^2 H_m)=0,
\end{equation}
where $\eta$ depending only on $n,\lambda,\Lambda$ and $\alpha$, is as in \Cref{L-C2a.K1}.

We prove the above by induction. For $m=1$, by setting $H_0\equiv 0$ and \Cref{L-C2a.K1}, \crefrange{e.C2a.S4}{e.C2a.S5} hold clearly. Suppose that the conclusion holds for $m=m_0$. We need to prove that the conclusion holds for $m=m_0+1$.

Let $r=\eta ^{m_{0}}$, $y=x/r$, $s=t/r^2$ and
\begin{equation}\label{e.C2a.S.v}
  v(y,s)=\frac{u(x,t)-H_{m_0}(x)}{r^{2+\alpha}}.
\end{equation}
Then $v$ satisfies
\begin{equation}\label{e.C2a.S.F}
\left\{\begin{aligned}
&v_s-\tilde{F}(D^2v,y,s)=\tilde{f}&& ~~\mbox{in}~~\tilde{\Omega}\cap Q_1;\\
&v=\tilde{g}&& ~~\mbox{on}~~\partial \tilde{\Omega}\cap Q_1,
\end{aligned}\right.
\end{equation}
where for $M \in \mathcal{S}^{n}$ and $(y,s)\in\tilde{\Omega}\cap Q_1$,
\begin{equation*}\label{e.c2a.S.b}
  \begin{aligned}
&\tilde{F}(M,y,s)=\frac{1}{r^{\alpha}}
F(r^{\alpha}M+D^2 H_{m_0},x,t),~~\tilde{F_0}(M)=\frac{1}{r^{\alpha}}
F_0(r^{\alpha}M+D^2 H_{m_0}),\\
&\tilde{f}(y,s)=\frac{f(x,t)}{r^{\alpha}},~~\tilde{g}(y,s)=\frac{g(x,t)-H_{m_0}(x)}{r^{2+\alpha}}
~~\mbox{ and }~~\tilde{\Omega}=\frac{\Omega}{r}.\\
  \end{aligned}
\end{equation*}

In the following, we show that \cref{e.C2a.S.F} satisfies the assumptions of \Cref{L-C2a.K1}. First, it is easy to verify that
\begin{equation*}
\begin{aligned}
&\|v\|_{L^{\infty}(\tilde{\Omega}\cap Q_1)}\leq 1, \tilde{F}_0(0)=0, v(0)=|Dv(0)|=0,\\
 &\|\partial \tilde{\Omega}\cap Q_1\|_{C^{1,\alpha}(0)} \leq r^{\alpha}\|\partial\Omega\cap Q_1\|_{C^{1,\alpha}(0)}\leq \delta
\end{aligned}
\end{equation*}
and for any $0<\rho<1$,
\begin{equation*}
\begin{aligned}
&\|\tilde{f}\|^*_{L^{n+1}(\tilde{\Omega}\cap Q_{\rho})}
=\frac{\|f\|^*_{L^{n+1}(\Omega\cap Q_{\rho r})}}{r^{\alpha}}
\leq \frac{\|f\|_{L^{\infty}(\Omega\cap Q_{\rho r})}}{r^{\alpha}}
\leq \delta \leq \delta \rho^{-1+\alpha}.
\end{aligned}
\end{equation*}
That is, $\|\tilde f\|_{C^{-1,\alpha}(0)}\leq \delta$.

From \cref{e.C2a.S6}, there exists a constant $C_0$ depending only on $n,\lambda,\Lambda$ and $\alpha$ such that $\|H_m\|\leq C_0$ ($\forall~0\leq m\leq m_0$). Thus, we have
\begin{equation*}
  \begin{aligned}
    |\tilde{F}(M,y,s)-\tilde{F_0}(M)|\leq& \frac{1}{r^{\alpha}}
    \left|F(r^{\alpha}M+D^2 H_{m_0},x,t)-F_0(r^{\alpha}M+D^2 H_{m_0})\right|\\
    \leq& \frac{1}{r^{\alpha}}\left(\beta(x,t)\left(r^{\alpha}\|M\|+C_0\right)
    +\gamma(x,t)\right)\\
    \leq& \delta \left(\|M\|+1\right)
    :=\tilde{\beta}(y,s)\|M\|+\tilde{\gamma}(y,s).
  \end{aligned}
\end{equation*}
Then
\begin{equation*}
\|\tilde{\beta}\|_{L^{\infty}(\tilde{\Omega}\cap Q_1)} \leq \delta,~~
\|\tilde{\gamma}\|_{L^{\infty}(\tilde{\Omega}\cap Q_1)} \leq \delta.
\end{equation*}
In addition,
\begin{equation*}
\begin{aligned}
|\tilde{g}(y,s)|
&\leq \frac{1}{r^{2+\alpha}}\left(\frac{\delta}{2}|(x,t)|^{2+\alpha}+
C_0\cdot \frac{\delta}{2C_0}|(x,t)|^{2+\alpha}\right)\\
&\leq  \delta |(y,s)|^{2+\alpha}\leq \delta |(y,s)|^{1+\alpha},~~\forall ~(y,s)\in \partial \tilde{\Omega}\cap Q_{1}.
\end{aligned}
\end{equation*}
Hence,
\begin{equation*}
  \|\tilde{g}\|_{C^{1,\alpha}(0)}\leq \delta.
\end{equation*}

Therefore, from above calculations, \cref{e.C2a.S.F} satisfies the assumptions of \Cref{L-C2a.K1}. Then there exists a $2$-form $\bar H (y)$ such that
\begin{equation*}
\begin{aligned}
    \|v-\bar H\|_{L^{\infty }(\tilde{\Omega} _{\eta})}&\leq \eta ^{2+\alpha},
\end{aligned}
\end{equation*}
\begin{equation*}
\tilde{F}_0(D^2 \bar H)=0
\end{equation*}
and
\begin{equation*}
\|\bar H\|\leq C_2+1.
\end{equation*}

Let $H_{m_0+1}(x)=H_{m_0}(x)+r^{\alpha}\bar H(x)$. Then \cref{e.C2a.S6} and \cref{e.C2a.S5} hold for $m_0+1$. Recalling \cref{e.C2a.S.v}, we have
\begin{equation*}
  \begin{aligned}
&\|u-H_{m_0+1}(x)\|_{L^{\infty}(\Omega_{\eta^{m_0+1}})}\\
&= \|u-H_{m_0}(x)-r^{\alpha}\bar H(x)\|_{L^{\infty}(\Omega_{\eta r})}\\
&= \|r^{2+\alpha}v-r^{2+\alpha}\bar H(y)\|_{L^{\infty}(\tilde{\Omega}_{\eta})}\\
&\leq r^{2+\alpha}\eta^{2+\alpha}=\eta^{(m_0+1)(2+\alpha)}.
  \end{aligned}
\end{equation*}
Hence, \cref{e.C2a.S4} holds for $m=m_0+1$. By induction, the proof is completed.\qed~\\

Now, we give the~\\
\noindent\textbf{Proof of \Cref{C2a-1}.} As before, the following proof is mere a normalization procedure. In the following proof, $C$ always denotes a constant depending only on $n,\lambda,\Lambda,\alpha,\|\beta\|_{C^{\alpha}(0)}$ and $\|(\partial \Omega)_1\|_{C^{1,\alpha}(0)}$. Denote
\begin{equation*}
U=2(\|u\|_{L^{\infty}(\Omega_1)}+|F_0(0)|+\|f\|_{C^{\alpha}(0)}+\|g\|_{C^{2,\alpha}(0)}
+\|\gamma\|_{C^{\alpha}(0)}).
\end{equation*}

First, let
\begin{equation*}
F_1(M,x,t)=F(M,x,t)+f(0),~~ F_{10}(M)=F_0(M)+f(0)
\end{equation*}
for $M \in \mathcal{S}^n$ and $(x,t)\in \Omega_1$. In the following, $F_0$ is transformed in the same way as $F$ and we omit the details. Then $u$ satisfies
\begin{equation*}
\left\{\begin{aligned}
&u_t-F_1(D^2u,x,t)=f_1&& ~~\mbox{in}~~\Omega_1;\\
&u=g&& ~~\mbox{on}~~(\partial \Omega)_1,
\end{aligned}\right.
\end{equation*}
where $f_1(x,t)=f(x,t)-f(0)$. Thus,
\begin{equation*}
  |f_1(x,t)|\leq [f]_{C^{\alpha}(0)}|(x,t)|^{\alpha}, ~~\forall ~(x,t)\in \Omega_1.
\end{equation*}

Next, set $u_1(x,t)=u(x,t)-P_g(x,t)$ and
\begin{equation*}
F_2(M,x,t)=F_1(M+D^{2}P_g,x,t)-(P_g)_t.
\end{equation*}
Then $u_1$ satisfies
\begin{equation*}
\left\{\begin{aligned}
&(u_1)_t-F_2(D^2u_1,x,t)=f_1&& ~~\mbox{in}~~\Omega_1;\\
&u_1=g_1&& ~~\mbox{on}~~(\partial \Omega)_1,
\end{aligned}\right.
\end{equation*}
where $g_1(x,t)=g(x,t)-P_g(x,t)$. Hence,

\begin{equation*}
  |g_1(x,t)|\leq [g]_{C^{2,\alpha}(0)}|(x,t)|^{2+\alpha}, ~~\forall ~(x,t)\in (\partial \Omega)_1.
\end{equation*}


Finally, for $0<\rho<1$, let
\begin{equation*}
y=\frac{x}{\rho}, s=\frac{t}{\rho^2}, u_2(y,s)=\frac{u_1(x,t)}{U}
\end{equation*}
and
\begin{equation*}
F_3(M,y,s)= \frac{\rho^2}{U} F_2\left(\frac{UM}{\rho^2},x,t\right).
\end{equation*}
Then $u_2$ satisfies
\begin{equation}\label{F5-1}
\left\{\begin{aligned}
&(u_2)_t-F_3(D^2u_2,y,s)=f_2&& ~~\mbox{in}~~\tilde{\Omega}_1;\\
&u_2=g_2&& ~~\mbox{on}~~(\partial \tilde\Omega)_1,
\end{aligned}\right.
\end{equation}
where $f_2(y,s)=\rho^2 f_1(x,t)/U$, $g_2(y,s)=g_1(x,t)/U$ and $\tilde{\Omega}=\Omega/\rho$.

Now, we can check that \cref{F5-1} satisfies the conditions of \Cref{t-C2a.S} by choosing a proper $\rho$. First, it can be checked easily that
\begin{equation*}
  \begin{aligned}
  & \|u_2\|_{L^{\infty}(\tilde{\Omega}_1)}
  \leq \frac{1}{U}\left(\|u\|_{L^{\infty}(\Omega_1)}+\|g\|_{C^{2,\alpha}(0)}
  \right)\leq 1/2,\\
&|F_{30}(0)|=\frac{\rho^2}{U}\left|F_0(D^2P_g)-(P_g)_t+f(0)\right|\\
&~~~~~~~~\leq \frac{\rho^2}{U}\left(|F_0(0)|+(\Lambda+1)\|g\|_{C^{2,\alpha}(0)}+\|f\|_{C^{\alpha}(0)}\right)
\leq \rho^2\left(\Lambda+1\right)\\
  &u_2(0)=|Du_2(0)|=0,\\
  &|f_2(y,s)|= \frac{\rho^2|f_1(x,t)|}{U}\leq \rho^{2+\alpha}|(y,s)|^{\alpha},~~\forall ~(y,s)\in\tilde{\Omega}_1,\\
  &|g_3(y,s)|=\frac{|g_2(x,t)|}{U}
  \leq \rho^{2+\alpha}|(y,s)|^{2+\alpha},~~\forall ~(y,s)\in\partial \tilde{\Omega}_1,\\
  & \|(\partial \tilde{\Omega})_1\|_{C^{1,\alpha}(0)} \leq \rho\|(\partial \Omega)_1\|_{C^{1,\alpha}(0)}.
  \end{aligned}
\end{equation*}
Moreover,
\begin{equation*}
  \begin{aligned}
    &|F_3(M,y,s)-F_{30}(M)|\\
    &=\frac{\rho^2}{U} \left|F\left(\frac{UM}{\rho^2}+D^2P_g,x,t\right)
    -F_0\left(\frac{UM}{\rho^2}+D^2P_g\right)\right|\\
&\leq \frac{\rho^2}{U} \left(\beta(x,t)\left(\frac{U}{\rho^2}\|M\|+\|D^2P_g\|\right)+\gamma(x,t)\right)\\
&\leq \left(\|\beta\|_{C^{\alpha}(0)}+1\right)\rho^{\alpha}|(y,s)|^{\alpha}(\|M\|+1).
  \end{aligned}
\end{equation*}

From above arguments, we can choose $\rho$ small enough (depending only on $n, \lambda,\Lambda,\alpha$, $\|\beta\|_{C^{\alpha}(0)}$ and $\|(\partial \Omega)_1\|_{C^{1,\alpha}(0)}$) such that \cref{e.C2a.S.be} holds. By \Cref{t-C2a.S}, $u_2 \in C^{2,\alpha}(0)$, i.e., there exists a $2$-form $H$ such that
\begin{equation*}
  |u_2(y,s)-H(y)|\leq C |(y,s)|^{2+\alpha}, ~~\forall ~(y,s)\in \tilde \Omega_{1},
\end{equation*}
\begin{equation*}
  \|H\|\leq C
\end{equation*}
and
\begin{equation*}
  F_{30}(D^2 H)=0.
\end{equation*}
That is, $u\in C^{2,\alpha}(0)$, i.e.,
\begin{equation*}
  \left|u(x,t)-P(x,t)\right|\leq CU \frac{|(x,t)|^{2+\alpha}}{\rho^{2+\alpha}} , ~~\forall ~(x,t)\in \Omega_{\rho},
\end{equation*}
where
\begin{equation*}
P(x,t)=P_g(x,t)+\frac{U}{\rho^2}H(x).
\end{equation*}
Furthermore,
\begin{equation*}
  \begin{aligned}
    0=&F_{30}(D^2H)\\
    =& \frac{\rho^2}{U} \left(F_0\left(D^2P_g
    +\frac{U}{\rho^2} D^2H \right)-(P_g)_t+P_f\right)\\
    =&\frac{\rho^2}{U} \left(F_0\left(D^2 P\right)-P_t+P_f\right).
  \end{aligned}
\end{equation*}
Therefore, \crefrange{e.C2a-1}{e.C2a.3-1} hold.

~\qed~\\

The \Cref{C2a} is an easy consequence of \Cref{C2a-1}:~\\
\noindent\textbf{Proof of \Cref{C2a}.} From \Cref{C1a}, $u\in C^{1,\alpha}(0)$. Let
\begin{equation*}
v(x,t)=u(x,t)-P_g(x,t)+((P_{g})_n(0)-u_n(0))\left(x_n-P_{\Omega}(x',t)\right)
\end{equation*}
and
\begin{equation*}
\begin{aligned}
\tilde{F}(M,x,t)=&F(M+D^2P_g+((P_{g})_n(0)-u_n(0))D^2P_{\Omega},x,t)\\
&-(P_{g})_t-((P_{g})_n(0)-u_n(0))(P_{\Omega})_t
\end{aligned}
\end{equation*}
for $(M,x,t) \in \mathcal{S}^n\times \Omega_1$.
Then $v$ satisfies
\begin{equation*}
\left\{\begin{aligned}
&v_t-\tilde{F}(D^2v,x,t)=f&& ~~\mbox{in}~~\Omega_1;\\
&v=\tilde{g}&& ~~\mbox{on}~~(\partial \Omega)_1,
\end{aligned}\right.
\end{equation*}
where
\begin{equation*}
  \tilde{g}(x,t)=g(x,t)-P_g(x,t)+((P_{g})_n(0)-u_n(0))\left(x_n-P_{\Omega}(x',t)\right).
\end{equation*}
By \cref{e1.7} and noting $P_{\Omega}\in \mathcal{HP}_2$,
\begin{equation*}
v(0)=|Dv(0)|=0.
\end{equation*}
In addition, since $g\in C^{2,\alpha}(0)$ and $(\partial \Omega)_1\in C^{2,\alpha}(0)$,
\begin{equation*}
|\tilde{g}(x,t)|\leq C|(x,t)|^{2+\alpha},~~\forall ~(x,t)\in (\partial \Omega)_1.
\end{equation*}
By applying \Cref{C2a-1} directly, \crefrange{e.C2a}{e.C2a.3} hold.~\qed~\\

\section{Boundary $C^{k,\alpha}$ regularity and the application to the regularity of free boundaries }\label{Sec:Cka_freeb}
In this section, we first prove the boundary pointwise $C^{k,\alpha}$ regularity for any $k\geq 3$ (\Cref{Ckla} and \Cref{Cka}) and then give the proof of the higher regularity of free boundaries in obstacle-type problems (\Cref{FreeBd}). For the boundary pointwise $C^{k,\alpha}$ regularity, firstly, we establish the regularity for the model problems (i.e., problems with homogeneous conditions and flat boundaries) in \Cref{subsec:Cka.flat}. Then we obtain the boundary pointwise $C^{k,\alpha}$ regularity by means of the perturbation skill as before.

\subsection{Boundary $C^{k, \alpha}$ regularity on a flat boundary}\label{subsec:Cka.flat}

The following lemma states the local $C^{2,\alpha}$ regularity up to the boundary, which can be derived by combining the interior pointwise regularity and the boundary pointwise regularity and is similar to the that of elliptic equations as shown in \cite[Lemma 9.1]{lian2020pointwise}.

\begin{lemma}\label{L.C2a.bar}
Let $0<\alpha<\bar\alpha$ and $u$ be a viscosity solution of
\begin{equation*}
\left\{\begin{aligned}
&u_t-F(D^2u,x,t)=0&& ~~\mbox{in}~~Q_1^+;\\
&u=0&& ~~\mbox{on}~~S_1.
\end{aligned}\right.
\end{equation*}
Suppose that $F\in C^{\alpha}(\bar{Q}_1^+)$ (i.e., $F$ satisfies \cref{e.C2a.beta} at every point of $\bar{Q}_1^+$ with the same $\beta$ and $\gamma$) and is convex in $M$.

Then $u\in C^{2,\tilde\alpha}(\bar{Q}^+_{1/2})$ for some $0<\tilde{\alpha}\leq \alpha$ and
\begin{equation*}
 \|u\|_{C^{2,\tilde\alpha}(\bar{Q}^+_{1/2})}\leq C \left(\|u\|_{L^{\infty }(Q_1^+)}+\|F\|_{\bar{\Omega}}+\|\gamma\|_{C^{\alpha}(\bar Q^+_1)}\right),
\end{equation*}
where $C$ depends only on $n,\lambda, \Lambda, \alpha$ and $\|\beta\|_{C^{\alpha}(\bar Q^+_1)}$.
\end{lemma}

\proof Similar to the $C^{1,\alpha}$ regularity, the conclusion follows by combining the interior regularity and the boundary regularity. Since the proof is a little complicated, we give it in details. As before, we assume
\begin{equation*}
\|u\|_{L^{\infty }(Q_1^+)}+\|F\|_{\bar{\Omega}}+\|\gamma\|_{C^{\alpha}(\bar Q^+_1)}\leq 1.
\end{equation*}
To prove the lemma, we only need to show that for any $p=(x,t)\in Q^+_{1/2}$, there exists a quadratic polynomial $P_p$ such that
\begin{equation}\label{e.l61-3}
  |u(y,s)-P_p(y,s)|\leq C|(y,s)-(x,t)|^{2+\tilde \alpha},~~\forall ~(y,s)\in \bar{Q}^+_{1/2}(p),
\end{equation}
and
\begin{equation}\label{e.glo.p}
  |D P_{p}(x,t)|+|D^2 P_{p}(x,t)|\leq C,
\end{equation}
where $C$ depends only on $n,\lambda, \Lambda,\alpha$ and $\|\beta\|_{C^{\alpha}(\bar{Q}^+_1)}$. Unless otherwise stated, $C$ always has the same dependence in this proof.

Let $p'=(x',0,t)\in S_1$ be the projection of $p$ to $S_1$. By the boundary $C^{2,\alpha}$ regularity (\Cref{C2a}), there exists a quadratic polynomial $P_{p'}$ such that
\begin{equation}\label{e4.1}
  \begin{aligned}
    &|u(y,s)-P_{p'}(y,s)|\leq C|(y,s)-(x',0,t)|^{2+\alpha}, ~\forall ~(y,s)\in \bar{Q}_{1/2}^{+}(p'),\\
    &|D P_{p'}(x',0,t)|+|D^2 P_{p'}(x',0,t)|\leq C
  \end{aligned}
\end{equation}
and
\begin{equation}\label{e4.2}
  (P_{p'})_s-F_{p'}(D^2 P_{p'})=0,
\end{equation}
where $F_{p'}$ is the fully nonlinear operator in \Cref{e.C2a.beta} corresponding to $p'$.

Set $v(y,s)=u(y,s)-P_{p'}(y,s)$ and $v$ satisfies
\begin{equation}\label{e.Cka-4}
  v_s-\tilde{F}(D^2v,y,s)=0~~~~\mbox{    in}~Q_{x_n}(p),
\end{equation}
where
\begin{equation*}
\tilde{F}(M,y,s)=F(M+D^2P_{p'},y,s)-(P_{p'})_s.
\end{equation*}
Define
\begin{equation*}
\tilde{F}_{p}(M)=F_{p}(M+D^2P_{p'})-(P_{p'})_s.
\end{equation*}
Then
\begin{equation*}
  \begin{aligned}
\left|\tilde{F}(M,y,s)-\tilde{F}_{p}(M)\right|&=\left|F(M+D^2P_{p'},y,s)
-F_{p}(M+D^2P_{p'})\right|\\
&\leq \beta\left((y,s),(x,t)\right)\left(\|M\|+\|D^2P_{p'}\|\right)+
    \gamma\left((y,s),(x,t)\right)\\
    &:=\beta\left((y,s),(x,t)\right)\|M\|+\tilde{\gamma}\left((y,s),(x,t)\right),
  \end{aligned}
\end{equation*}
where $\tilde{\gamma}\in C^{\alpha}(\bar Q^+_{1})$. Then $\tilde{F}$ satisfies \Cref{e.C2a.beta} with $\tilde{F}_p$ corresponding to $p$.

By virtue of \cref{e4.1},
\begin{equation*}
\|v\|_{L^{\infty}(Q_{x_n}(p))}\leq Cx_n^{2+\alpha}.
\end{equation*}
Moreover, with the aid of the triangle inequality, \cref{e.C2a.beta} implies
\begin{equation*}
|F_{p}(M)-F_{p'}(M)|\leq Cx_n^{\alpha}(\|M\|+1),~~\forall ~M\in \mathcal{S}^n.
\end{equation*}
Hence, by noting \cref{e4.2},
\begin{equation*}
|\tilde{F}_{p}(0)|=|F_{p}(D^2P_{p'})-(P_{p'})_s|= |F_{p}(D^2P_{p'})-F_{p'}(D^2P_{p'})|
\leq Cx_n^{\alpha}.
\end{equation*}
Therefore, from the interior $C^{2,\alpha}$ regularity (\cite[Theorem 1.1, Theorem 4.13]{MR1139064}), there exist $0<\tilde \alpha\leq \alpha$ and a quadratic polynomial $P$ such that
\begin{equation*}
  \begin{aligned}
|v(y,s)-P(y,s)|&\leq C\frac{|(y,s)-(x,t)|^{2+\tilde \alpha}}{|x_n|^{2+\tilde \alpha}}
\left(\|v\|_{L^{\infty}(Q_{x_n}(p))}+x_n^2|\tilde{F}_{p}(0)|\right)\\
  &\leq C |(y,s)-(x,t)|^{2+\tilde \alpha},~~\forall~~(y,s)\in Q_{x_n/2}(p)
  \end{aligned}
\end{equation*}
and
\begin{equation}\label{e.gol.p2}
  \begin{aligned}
&|P(x,t)|=|v(x,t)|\leq  Cx_n^{2+\alpha},\\
&|D P(x,t)|\leq \frac{C}{x_n}
\left(\|v\|_{L^{\infty}(Q_{x_n}(p))}+x_n^2|\tilde{F}_{p}(0)|\right)\leq Cx_n^{1+\alpha},\\
&|D^2P(x,t)|\leq \frac{C}{x_n^2}
\left(\|v\|_{L^{\infty}(Q_{x_n}(p))}+x_n^2|\tilde{F}_{p}(0)|\right)\leq Cx_n^{\alpha}.
  \end{aligned}
\end{equation}

Let $P_p=P_{p'}+P$. By combining \cref{e4.1} and \cref{e.gol.p2}, \cref{e.glo.p} holds. In addition, if $|(y,s)-(x,t)|<x_n/2$,
\begin{equation*}
  |u(y,s)-P_p(y,s)|=|v(y,s)-P(y,s)|\leq C|(y,s)-(x,t)|^{2+\tilde \alpha}.
\end{equation*}
If $|(y,s)-(x,t)|\geq x_n/2$,
\begin{equation*}
  \begin{aligned}
    |u&(y,s)-P_p(y,s)|\\
    &\leq |u(y,s)-P_{p'}(y,s)|+|P(y,s)|\\
    &\leq C|(y,s)-(x',0,t)|^{2+\alpha}\\
    &~~~~~~~~+ \left(|P(x,t)|+|DP(x,t)||(y,s)-(x,t)|
    +|D^2P(x,t)||(y,s)-(x,t)|^2\right)\\
    &\leq C|(y,s)-(x',0,t)|^{2+\alpha}+ C(x_n^{2+\alpha}+|(y,s)-(x,t)|x_n^{1+\alpha}+|(y,s)-(x,t)|^2x_n^{\alpha})\\
    &\leq  C|(y,s)-(x,t)|^{2+\alpha}.
  \end{aligned}
\end{equation*}
That is, \cref{e.l61-3} holds and the proof is completed. ~\qed \\

\begin{remark}\label{re.Cka.UtB.2}
There exists a universal constant $0<\hat{\alpha}<1$ such that the interior $C^{2,\alpha}$ regularity holds for any $0<\alpha<\hat{\alpha}$ ($\hat \alpha$ is exactly the $\alpha$ in \cite[Theorem 4.13]{MR1139064}). By taking this fact into account, we can obtain $u\in C^{2,\alpha}$ for any $0<\alpha<\min (\bar{\alpha},\hat{\alpha})$ if $F\in C^{\alpha}(\bar{Q}_1^+)$.
\end{remark}

Now, we establish the $C^{k,\alpha}$ regularity up the boundary for the model problem.

\begin{lemma}\label{L.Cka.F0}
Let $u$ be a viscosity solution of
\begin{equation*}
\left\{\begin{aligned}
&u_t-F(D^2u,x,t)=0&& ~~\mbox{in}~~Q_1^+;\\
&u=0&& ~~\mbox{on}~~S_1.
\end{aligned}\right.
\end{equation*}
Suppose that $F\in C^{k-2,1}(\mathcal{S}^{n}\times \bar Q^+_1)$ ($k\geq 3$) is convex in $M$ and satisfies \cref{e1.8}.

Then $u\in C^{k,\alpha}(\bar{Q}^+_{1/2})$ for any $0<\alpha<1$ and
\begin{equation*}
\|u\|_{C^{k,\alpha}(\bar{Q}^+_{1/2})}\leq C_k,
\end{equation*}
where $C_k$ depends only on $k,n,\lambda, \Lambda,\alpha,K_0,\omega$ and $\|u\|_{L^{\infty }(Q_1^+)}$. Here,
\begin{equation*}
\omega(r):=\|F\|_{C^{k-2,1}(\bar{\textbf{B}}_r\times \bar Q_{1}^{+})},~~\forall ~r>0,
\end{equation*}
where $\textbf{B}_r=\left\{M\big| \|M\|<r\right\}$.

In particular, $u\in C^{k,\alpha}(0)$ and there exists $P\in \mathcal{P}_k$ with the following expression
\begin{equation}\label{e4.3}
P(x,t)=\sum_{i=1}^{k}H_i(x,t),~H_i \mbox{ is an } i\mbox{-form }
\end{equation}
such that
\begin{equation}\label{e.l62-1}
  |u(x,t)-P(x,t)|\leq C_k |(x,t)|^{k+\alpha}, ~~\forall ~(x,t)\in Q_{1}^+,
\end{equation}
\begin{equation}\label{e.l62-3}
  \|P\|\leq C_k
\end{equation}
and
\begin{equation}\label{e.l62-2}
  |P_t-F(D^2P,x,t)|\leq C_k |(x,t)|^{k-1}, ~~\forall ~(x,t)\in Q_{1}^+.
\end{equation}
\end{lemma}

\proof We prove this lemma by the standard technique of difference quotient. First, we show that $F$  satisfies \cref{e.C2a.beta} at every point of $\bar{Q}_1^+$. Indeed, for any $p_0=(x_0,t_0)\in \bar{Q}_1^+$, define $F_{p_0}(M)=F(M,x_0,t_0)$. Since $F$ satisfies \cref{e1.8}, for any $M\in \mathcal{S}^{n}$ and $(x,t)\in \bar{Q}_1^+\cap \bar Q_{1}(x_0,t_0)$,
\begin{equation*}
  \begin{aligned}
|F&(M,x,t)-F_{p_0}(M)|\\
=&|F(M,x,t)-F(M,x_0,t_0)|\\
= &|F(M,x,t)-F(0,x,t)+F(0,x,t)-F(0,x_0,t_0)+F(0,x_0,t_0)-F(M,x_0,t_0)|\\
=&\left|\int_{0}^{1}\left(F_{ij}(\tau M,x,t)-F_{ij}(\tau M,x_0,t_0)\right)M_{ij}d\tau+F(0,x,t)-F(0,x_0,t_0)\right|\\
\leq&n^2K_0|(x,t)-(x_0,t_0)|\|M\|+C|(x,t)-(x_0,t_0)|,
  \end{aligned}
\end{equation*}
where $C$ depends only on $\omega$. Hence, $F$ satisfies \cref{e.C2a.beta} at every point of $\bar{Q}_1^+$.

By \Cref{L.C2a.bar}, for some universal constant $\tilde{\alpha}>0$,
\begin{equation*}
\|u\|_{C^{2,\tilde\alpha}(\bar{Q}_{3/4}^+)}\leq C,
\end{equation*}
where $C$ depends only on $n,\lambda, \Lambda,K_0,\omega$ and $\|u\|_{L^{\infty }(Q_1^+)}$. Unless stated otherwise, $C$ always has the same dependence (may also depending on $k$) in the following proof.

Now, we show that $u\in C^{3,\alpha}$. Let $h>0$ be small and $1\leq l\leq n-1$. By taking the difference quotient along $e_l$ on both sides of the equation,
\begin{equation*}
  (\Delta_l^h u)_{t}-a_{ij}(\Delta_l^h u)_{ij}=G ~~~~\mbox{ in}~~ Q^+_{5/8},
\end{equation*}
where
\begin{equation*}
  \begin{aligned}
&\Delta_l^h u(x,t)=\left(u(x+he_l,t)-u(x,t)\right)/h,\\
&a_{ij}(x,t)=\int_{0}^{1}F_{ij}(\xi,x+\tau h e_l,t)d \tau,\\
&G(x,t)=\int_{0}^{1}F_{x_l}(\xi,x+\tau h e_l,t)d \tau \\
  \end{aligned}
\end{equation*}
and
\begin{equation*}
  \begin{aligned}
\xi=\tau D^2u(x+he_l,t)+(1-\tau)D^2u(x,t).
  \end{aligned}
\end{equation*}
Thus, $a_{ij}$ is uniformly elliptic with $\lambda$ and $\Lambda$. Since $F\in C^{1,1}(\mathcal{S}^{n}\times \bar Q^+_1)$, we have $F_{ij},F_{x_l}\in C^{0,1}(\mathcal{S}^{n}\times \bar Q^+_{5/8})$. Moreover, note $D^2 u \in C^{\tilde \alpha}(\bar Q^+_{5/8})$. Hence, it is easy to check that $\|a_{ij}\|_{C^{\tilde \alpha}(\bar{Q}^+_{5/8})}\leq C$ and $\|G\|_{C^{\tilde \alpha}(\bar{Q}^+_{5/8})}\leq C$.

By the Schauder's estimate for linear equations, we have $\Delta_l^h u \in C^{2,\tilde \alpha}(\bar{Q}^+_{1/2})$. Hence, $u_l \in C^{2,\tilde \alpha}(\bar{Q}^+_{1/2})$ and
\begin{equation*}
\|u_l\|_{C^{2,\tilde \alpha}(\bar{Q}_{1/2}^+)}\leq C.
\end{equation*}
It follows that $u_{i l}\in C^{1,\tilde \alpha}$ ($i=1,2,\cdots,n; l=1,2,\cdots,n-1$). In a similar way, by taking the difference quotient with respect to $t$, we have $u_t \in C^{1,\tilde \alpha}(\bar Q^+_{1/2})$. Then by combining with $u_t-F(D^2u,x,t)=0$ and the implicit function theorem, $u_{nn}\in C^{1,\tilde \alpha}$. Thus, $u\in C^{3,\tilde \alpha}$ and
\begin{equation*}
\|u\|_{C^{3,\tilde \alpha}(\bar Q_{1/2}^+)}\leq C.
\end{equation*}

Since $u\in C^{3,\tilde \alpha}$, for $1\leq l \leq n-1$, $u_l$ satisfies
\begin{equation}\label{uksat}
  (u_l)_t-a_{ij}(u_l)_{ij}=G(x,t),
\end{equation}
where
\begin{equation}\label{aijuk}
  a_{ij}(x,t)=F_{ij}(D^2u,x,t), ~~G(x,t)=F_{l}(D^2u,x,t).
\end{equation}

Then it can be seen that $a_{ij}\in C^{0,1}(\bar{Q}^+_{1/2})$ and $G\in C^{0,1}(\bar{Q}^+_{1/2})$. By the Schauder's estimates for linear equations, $u_l\in C^{2,\alpha}$ for any $0<\alpha<1$. Similarly, by taking difference quotient with $t$, we have $u_t\in C^{1,\alpha}$. Then $u\in C^{3,\alpha}$ and
\begin{equation*}
\|u\|_{C^{3,\alpha}(\bar{Q}_{1/4}^+)}\leq C.
\end{equation*}

Since $F\in C^{k-2,1}$, by considering \cref{uksat} iteratively and a covering argument, we obtain
\begin{equation*}
\|u\|_{C^{k,\alpha}(\bar{Q}^+_{1/2})}\leq C.
\end{equation*}

In particular, $u\in C^{k,\alpha}(0)$ and there exists $P\in \mathcal{P}_k$ such that \cref{e.l62-1} holds. Note that $u(0)=0$, $u_i(0)=0~(1\leq i<n \text{ or } i=t)$, $u_{ij}(0)=0 ~(i,j\neq n)$ etc. Hence, $P$ can be written as \cref{e4.3}. Finally, from $u\in C^{k,\alpha}(\bar{Q}^+_{1/2})$ and \cref{e.l62-1},
\begin{equation*}
  \begin{aligned}
    &|P_t-F(D^2P,x,t)|\\
    &=|P_t-F(D^2P,x,t)-\left(u_t-F(D^2u,x,t)\right)|\\
    &\leq \left|P_t-u_t\right|+\left|\int_{0}^{1} F_{ij}(\tau D^2P+(1-\tau) D^2u,x,t)
    \left(P_{ij}-u_{ij}\right)d\tau \right|\\
    &\leq C |(x,t)|^{k-2+\alpha}.
  \end{aligned}
\end{equation*}
That is,
\begin{equation*}
  D^{i} \Big(P_t-F(D^2P,x,t)\Big)\Big|_{(x,t)=0}=0,~\forall ~0\leq i\leq k-2.
\end{equation*}
Combining with $F\in C^{k-2,1}$, \cref{e.l62-2} holds. \qed~\\

\begin{remark}\label{re4.2}
Since the $C^{k,\alpha}$ regularity is derived from the Schauder's estimate for linear equations and the coefficients $a^{ij}$ depend on $u$, one can't obtain an estimate depending on $u$ explicitly as the
$C^{1,\alpha}$ and $C^{2,\alpha}$ regularity.
\end{remark}

\begin{remark}\label{re4.1}
For boundary $C^{1,\alpha}$ regularity, the model equation is $u\in S(\lambda,\Lambda,0)$. The boundary regularity is obtained by comparing $u$ with $x_n$ and we have $C^{1,\bar{\alpha}}$ regularity for some universal $0<\bar{\alpha}<1$.

For boundary $C^{2,\alpha}$ regularity, the model equation is $u_t-F(D^2u)=0$. The boundary regularity is obtained by differentiating the equation and using the boundary $C^{1,\alpha}$ regularity. As a result, we have $C^{2,\alpha}$ regularity for any $0<\alpha<\bar{\alpha}$.

For higher order regularity, the model equation is $u_t-F(D^2u,x,t)=0$ where $F$ is smooth. However, the boundary regularity is derived from Schauder estimate for linear equations (see \cref{uksat}). Hence, we have $C^{k,\alpha}$ regularity for any $0<\alpha<1$.
\end{remark}
~\\


\subsection{Boundary pointwise $C^{k, \alpha}$ regularity on general boundaries}
\label{subsec:Cka.free}

In this subsection, we first prove an important immediate result (see \Cref{t.Cka.1}) and then \Cref{Ckla} and \Cref{Cka} can be obtained by normalization procedures.

%

\begin{theorem}\label{t.Cka.1}
Let $k\geq 1$, $0<\alpha<1$ ($0<\alpha<\bar{\alpha}$ if $k=1$) and $u$ be a viscosity solution of
\begin{equation*}
\left\{\begin{aligned}
&u_t-F(D^2 u,x,t)=f&& ~~\mbox{in}~~\Omega_1;\\
&u=g&& ~~\mbox{on}~~(\partial \Omega)_1.
\end{aligned}\right.
\end{equation*}
Suppose that $F\in C^{k-1,\alpha}(0)$, $f\in C^{k-1,\alpha}(0)$, $g\in C^{k+1,\alpha}(0)$ and $(\partial \Omega)_1\in C^{1,\alpha}(0)$. Furthermore, we assume $u\in C^{k,\alpha}(0)$ and
\begin{equation}\label{e.Dkg}
u(0)=\cdots=|D^ku(0)|=|Dg(0)|\cdots=|D^kg(0)|=0.
\end{equation}

Then $u\in C^{k+1,\alpha}(0)$, i.e., there exists a $(k+1)$-form $H$ such that
\begin{equation}\label{e.Cka.u1}
  \begin{aligned}
&|u(x,t)-H(x,t)|\leq C |(x,t)|^{k+1+\alpha},~~\forall ~(x,t)\in \Omega_{1},\\
&|D^{k+1}u(0)|\leq C
  \end{aligned}
\end{equation}
and
\begin{equation}\label{e.Cka.u2}
  \begin{aligned}
    &|H_t-F_0(D^2 H,x,t)-P_f|\leq C|(x,t)|^{k},\\
    &H(x',0,t)\equiv P_g(x',0,t),
  \end{aligned}
\end{equation}
where $C$ depends only on $k,n,\lambda, \Lambda, \alpha,K_0,\omega$, $\|\beta\|_{C^{k-1,\alpha}(0)}$, $\|(\partial \Omega)_1\|_{C^{1,\alpha}(0)}$, $\|u\|_{L^{\infty}(\Omega_{1})}$, $\|f\|_{C^{k-1,\alpha}(0)}$, $\|g\|_{C^{k+1,\alpha}(0)}$ and $\|\gamma\|_{C^{k-1,\alpha}(0)}$.
\end{theorem}

We prove the theorem by induction. For $k=1$, \Cref{t.Cka.1} reduces to \Cref{C2a-1}, which has been proved in last section. Suppose that \Cref{t.Cka.1} holds for $k-1$ and we only need to prove it for $k$ in the following.

First, we prove a key step.
\begin{lemma}\label{L.Cka.2}
Let $0<\alpha<1$ and $F\in C^{k-1,\alpha}(0)$. There exists $\delta>0$ depending only on $k,n,\lambda,\Lambda,\alpha,K_0$ and $\omega$ such that if $u$ satisfies
\begin{equation*}
\left\{\begin{aligned}
&u_t-F(D^2u,x,t)=f&& ~~\mbox{in}~~\Omega_1;\\
&u=g&& ~~\mbox{on}~~(\partial \Omega)_1
\end{aligned}\right.
\end{equation*}
with $u\in C^{k,\alpha}(0)$ and
\begin{equation*}
  \begin{aligned}
    &\|u\|_{L^{\infty}(\Omega_1)}\leq 1,~~ u(0)=\cdots=|D^ku(0)|=0,\\
    &\|\beta\|_{C^{k-2,\alpha}(0)}\leq \delta,~~
    \|\gamma\|_{C^{k-2,\alpha}(0)}\leq \delta,\\
    &|f(x,t)|\leq \delta|(x,t)|^{k-2+\alpha},~~\forall ~(x,t)\in \Omega_1,\\
    &\|g\|_{C^{k,\alpha}(0)}\leq \delta,~~g(0)=\cdots=|D^kg(0)|=0,\\
    &\|(\partial \Omega)_1\|_{C^{1,\alpha}(0)} \leq \delta.
  \end{aligned}
\end{equation*}

Then there exists a $(k+1)$-form $H$ such that
\begin{equation*}
  \begin{aligned}
    & \|u-H\|_{L^{\infty}(\Omega_{\eta})}\leq \eta^{k+1+\alpha},\\
    &\|H\|\leq \bar{C}_k,\\
    &\left|H_t-F_0(D^2H,x,t)\right|\leq \bar{C}_k|(x,t)|^{k},~~\forall~~(x,t)\in \Omega_1,\\
  \end{aligned}
\end{equation*}
where $\bar{C}_k$ depends only on $k,n,\lambda,\Lambda,\alpha,K_0$ and $\omega$.
\end{lemma}
\proof As before, we prove the lemma by contradiction. Suppose that the conclusion is false. Then there exists $0<\alpha<1,K_0$, $\omega$ and sequences of $u_m,f_m,g_m,\Omega_m,F_m$ ($m\geq 1$) satisfying $u_m\in C^{k,\alpha}(0)$ and
\begin{equation*}
\left\{\begin{aligned}
&(u_m)_t-F_m(D^2u_m,x,t)=f_m&& ~~\mbox{in}~~\Omega_m\cap Q_1;\\
&u_m=g_m&& ~~\mbox{on}~~\partial \Omega_m\cap Q_1.
\end{aligned}\right.
\end{equation*}
In addition, $F_m \in C^{k-1,\alpha}(0)$ with $F_{0m}$, $\beta_m$, $\gamma_m$ and
\begin{equation}\label{e.4.F0m}
\|F_{0m}\|_{C^{k-1,1}(\bar{\textbf{B}}_r\times \overline{\Omega_m\cap Q_1})}\leq \omega(r),~\forall ~r>0.
\end{equation}
Furthermore,
\begin{equation*}
  \begin{aligned}
    &\|u_m\|_{L^{\infty}(\Omega_m\cap Q_1)}\leq 1,~~ u_m(0)=\cdots=|D^ku_m(0)|=0,\\
    &\|\beta_m\|_{C^{k-2,\alpha}(0)}\leq 1/m,~~
    \|\gamma_m\|_{C^{k-2,\alpha}(0)}\leq 1/m,\\
    &|f_m(x,t)|\leq |(x,t)|^{k-2+\alpha}/m,~~\forall ~(x,t)\in \Omega_m\cap Q_1,\\
    &\|g_m\|_{C^{k,\alpha}(0)}\leq 1/m,~~g_m(0)=\cdots=|D^kg_m(0)|=0,\\
    &\|\partial \Omega_m\cap Q_1\|_{C^{1,\alpha}(0)} \leq 1/m.
  \end{aligned}
\end{equation*}
But for any $(k+1)$-form $H$ satisfying $\|H\|\leq \bar{C}_k$ and
\begin{equation}\label{e.L.Cka.Q}
 \left|H_t-F_{0m}(D^2H,x,t)\right|\leq \bar{C}_k|(x,t)|^{k},
\end{equation}
we have
\begin{equation}\label{e.L.Cka.1}
  \|u_m-H\|_{L^{\infty}(\Omega_{m}\cap Q_{\eta})}> \eta^{k+1+\alpha},
\end{equation}
where $\bar{C}_k$ is to be specified later and $0<\eta<1$ is taken small such that
\begin{equation}\label{e.L.Cka.2}
\bar{C}_k\eta^{1-\alpha}<1/2.
\end{equation}

As before, there exists a subsequence (denoted by $u_m$ again) and $\bar u: Q_1^+\rightarrow \mathbb R$ such that $u_m\rightarrow \bar u$ uniformly in any compact subset of $Q^+_1$. In addition, since $F_{0m}\in C^{k-1,1}(\mathcal{S}^{n}\times \overline{\Omega_m\cap Q_1})$ with \cref{e.4.F0m}, there exist a subsequence of $F_{0m}$ (denoted by $F_{0m}$ again) and $\bar F\in C^{k-1,1}(\mathcal{S}^{n}\times \bar{Q}_1^+)$ such that $F_{0m}\rightarrow \bar F$ uniformly in $C^{k-1}$ in any compact subsets of $\mathcal {S}^{n}\times \bar{Q}^+_1$ and $\|\bar F\|_{C^{k-1,1}(\bar{\textbf{B}}_r\times \bar{Q}_1^+)}\leq \omega(r)$ for any $r>0$. Apparently, $\bar F$ is uniformly elliptic with constants $\lambda$ and $\Lambda$ and convex in $M$. Moreover, $\bar F$ satisfies \cref{e1.8} with $K_0$.

Furthermore, for any $Q\subset\subset Q_1^+$ and $\varphi\in C^2(\bar{Q})$, define $\psi_m(x,t)=-F_m(D^2\varphi,x,t)-f_m$ and $\psi(x,t)=-\bar F(D^2\varphi,x,t)$. Then
\begin{equation*}
  \begin{aligned}
&\|\psi_m-\psi\|_{L^{n+1}(Q)}\\
    &\leq \|F_m(D^2\varphi,x,t)-F_{0m}(D^2\varphi,x,t)\|_{L^{n+1}(Q)}\\
    &+\|F_{0m}(D^2\varphi,x,t)-\bar F(D^2\varphi,x,t)\|_{L^{n+1}(Q)}+\|f_m\|_{L^{n+1}(Q)}\\
    &\leq \|\beta_{m}\|_{L^{n+1}(Q)}\|D^2\varphi\|_{L^{\infty}(Q)}+\|\gamma_{m}\|_{L^{n+1}(Q)}\\
    &+\|F_{0m}(D^2\varphi,x,t)-\bar F(D^2\varphi,x,t)\|_{L^{n+1}(Q)}+\|f_m\|_{L^{n+1}(Q)}.
\end{aligned}
\end{equation*}
Thus, $\|\psi_m-\psi\|_{L^{n+1}(Q)}\rightarrow 0$ as $m\rightarrow \infty$. By \Cref{L21} and \Cref{L23}, $\bar u$ is a viscosity solution of
\begin{equation*}
\left\{\begin{aligned}
&\bar u_t-\bar F(D^2\bar u,x,t)=0&& ~~\mbox{in}~~Q_{1}^+;\\
&\bar u=0&& ~~\mbox{on}~~S_{1}.
\end{aligned}\right.
\end{equation*}

By induction, \Cref{t.Cka.1} holds for $k-1$, i.e., we have the boundary $C^{k,\alpha}$ regularity.
Note that $u_m(0)=\cdots=|D^{k}u_m(0)|=0$ and then
\begin{equation*}
\|u_m\|_{L^{\infty }(\Omega_m\cap Q_r)}\leq Cr^{k+\alpha} ,~~~~\forall~0<r<1.
\end{equation*}
Since $u_m$ converges to $\bar u$ uniformly, we have
\begin{equation*}
\|\bar u\|_{L^{\infty }(Q_r^+)}\leq Cr^{k+\alpha}, ~~~~\forall~0<r<1.
\end{equation*}
Hence, $\bar u(0)=\cdots=|D^{k}\bar u(0)|=0$. By the boundary estimate for $\bar u$ on a flat boundary (see \Cref{L.Cka.F0}), there exists a $(k+1)$-form $\bar{H}$ such that
\begin{equation}\label{e.cka-5}
\begin{aligned}
  &|\bar u(x,t)-\bar{H}(x,t)|\leq C_k |(x,t)|^{k+1+\alpha}, ~~\forall ~(x,t)\in Q_{1}^+,\\
  &\|\bar{H}\|\leq C_k,\\
  &|\bar{H}_t-\bar F(D^2\bar{H},x,t)|\leq C_k|(x,t)|^{k}, ~~\forall ~(x,t)\in Q_{1}^+,\\
\end{aligned}
\end{equation}
where $C_k$ depends only on $k,n,\lambda,\Lambda,\alpha,K_0$ and $\omega$.

Now, we try to construct a sequence of $k+1$-forms $\bar{H}_m$ such that \cref{e.L.Cka.Q} holds for $\bar{H}_m$ and $\bar{H}_m\rightarrow \bar{H}$ as $m\rightarrow \infty$. Let $H_m$ be $k+1$-forms ($m\geq 1$) with $\|H_m\|\leq 1$ to be specified later and $\bar{H}_m=H_m+\bar{H}$.  Let $G_m(x,t)=(\bar{H}_m)_t-F_{0m}(D^2\bar{H}_m,x,t)$. Since $G_m\in C^{k-1,1}(\overline{\Omega_m\cap{Q}_{1}})$, \cref{e.L.Cka.Q} holds for $\bar{H}_m$ if we show that $D^iG_m(0)=0$ for any $1\leq i\leq k-1$.

Note that $u_m(0)=\cdots=|D^{k}u_m(0)|=0$ and from the pointwise $C^{k,\alpha}$ regularity (note that \Cref{t.Cka.1} holds for $k-1$ by induction), we obtain (note $H=P_f=0$ in \cref{e.Cka.u2})
\begin{equation*}
|F_{0m}(0,x,t)|\leq C|(x,t)|^{k-1}.
\end{equation*}
Thus,
\begin{equation*}
  \begin{aligned}
&|G_m(x,t)|\\
&\leq |F_{0m}(D^2\bar{H}_m,x,t)-F_{0m}(0,x,t)|
    +|F_{0m}(0,x,t)|+|(\bar{H}_m)_t|\\
&=\left|\int_{0}^{1} F_{0m,ij}(\tau D^2\bar{H}_m,x,t)\bar{H}_{m,ij}d\tau \right|
    +|F_{0m}(0,x,t)|+|(\bar{H}_m)_t|\\
&\leq C|(x,t)|^{k-1}.
  \end{aligned}
\end{equation*}
That is, $D^iG_m(0)=0$ for any $1\leq i\leq k-2$. Hence, to verify \cref{e.L.Cka.Q} for $\bar{H}_m$, we only need to prove
\begin{equation*}
D^{k-1}G_m(0)=D^{k-1}\left((\bar{H}_m)_t- F_{0m}(D^2\bar{H}_m,x,t)\right)\bigg|_{(x,t)=0}=0.
\end{equation*}
Indeed, since $\bar{H}_{m}$ are $(k+1)$-forms,
\begin{equation*}
  \begin{aligned}
D^{k-1}G_m(0)&=D^{k-1}(\bar{H}_m)_t-F_{0m,ij}(0)D^{k-1}\bar{H}_{m,ij}-(D^{k-1}_{p}F_{0m})(0)\\
&=D^{k-1}\bar{H}_t-F_{0m,ij}(0)D^{k-1}\bar{H}_{ij}-(D^{k-1}_{p}F_{0m})(0)\\
&~~~~+D^{k-1}(H_m)_t-F_{0m,ij}(0)D^{k-1}H_{m,ij},
  \end{aligned}
\end{equation*}
where $D^{k-1}_{p}$ means taking $(k-1)$-th order derivatives only with respect
to $x$ or $t$. Since
\begin{equation*}
  \begin{aligned}
\lim_{m\rightarrow \infty} & D^{k-1}\bar{H}_t-F_{0m,ij}(0)D^{k-1}\bar{H}_{ij}-(D^{k-1}_{p}F_{0m})(0)\\ =&D^{k-1}\bar{H}_t-F_{0,ij}(0)D^{k-1}\bar{H}_{ij}-(D^{k-1}_{p}F_{0})(0)\\
=&D^{k-1}\left(\bar{H}_t-\bar F(D^2\bar{H},x,t)\right)\bigg |_{(x,t)=0}\\
=&0~(\mbox{by}~ \cref{e.cka-5})
  \end{aligned}
\end{equation*}
and $\lambda\leq F_{0m,ij}(0)\leq \Lambda$, we can choose proper $H_m$ such that $D^{k-1}G_m(0)=0$ for any $m\geq 1$ and $\|H_m\|\rightarrow 0$ as $m\rightarrow \infty$. Hence, there exists a large constant $\bar C_k$ depending only on $k,n,\lambda, \Lambda,\alpha,K_0$ and $\omega$ such that $\|\bar{H}_m\|\leq \bar C_k$ for any $m\geq 1$ and
\begin{equation*}
  \begin{aligned}
|G_m(x,t)|=|(\bar{H}_m)_t-F_{0m}(D^2\bar{H}_m,x,t)|\leq \bar C_k|(x,t)|^{k}.
  \end{aligned}
\end{equation*}

Let $m\rightarrow \infty$ in \cref{e.L.Cka.1} and we have
\begin{equation*}
    \|\bar u-\bar{H}\|_{L^{\infty}(Q_{\eta}^+)}\geq \eta^{k+1+\alpha},
\end{equation*}
However, by \cref{e.L.Cka.2} and \cref{e.cka-5},
\begin{equation*}
  \|\bar u-\bar{H}\|_{L^{\infty}(Q_{\eta}^+)}\leq \eta^{k+1+\alpha}/2,
\end{equation*}
which is a contradiction.  ~\qed~\\

As before, we now can prove the boundary pointwise $C^{k,\alpha}$ regularity by a scaling argument.
\begin{theorem}\label{t.Cka.S}
Let $0<\alpha <1$. Suppose that $F\in C^{k-1,\alpha}(0)$ is convex in $M$ and $u$ satisfies
\begin{equation*}
\left\{\begin{aligned}
&u_t-F(D^2u,x,t)=f&& ~~\mbox{in}~~\Omega_1;\\
&u=g&& ~~\mbox{on}~~(\partial \Omega)_1.
\end{aligned}\right.
\end{equation*}
Assume that $u\in C^{k,\alpha}(0)$ and
\begin{equation}\label{e.tCkas-be}
\begin{aligned}
&\|u\|_{L^{\infty}(\Omega_1)}\leq 1,~ u(0)=|Du(0)|=\cdots=|D^{k}u(0)|=0,\\
&|\beta(x,t)|\leq \frac{\delta}{2C_0}|(x,t)|^{k-1+\alpha}, ~~|\gamma(x,t)|\leq \frac{\delta}{2C_0}|(x,t)|^{k-1+\alpha},\\
&|f(x,t)|\leq \delta|(x,t)|^{k-1+\alpha}, ~~\forall ~(x,t)\in \Omega_1,\\
&|g(x,t)|\leq \frac{\delta}{2}|(x,t)|^{k+1+\alpha}, ~~\forall ~(x,t)\in (\partial \Omega)_1
~~\mbox{and}~~\|\partial \Omega\cap Q_1\|_{C^{1,\alpha}(0)} \leq \frac{\delta}{2C_0},\\
\end{aligned}
\end{equation}
where $\delta$ is as in \Cref{L.Cka.2} and $C_0$ depend only on $k,n,\lambda, \Lambda,\alpha,K_0$ and $\omega$.

Then $u\in C^{k+1,\alpha}(0)$ and there exists a $k+1$-form $H$ such that
\begin{equation}\label{e.tCkas-1}
  |u(x,t)-H(x,t)|\leq C |(x,t)|^{k+1+\alpha}, ~~\forall ~(x,t)\in \Omega_{1},
\end{equation}
\begin{equation}\label{e.tCkas-2}
\|H\|\leq C
\end{equation}
and
\begin{equation*}
  |H_t-F_0(D^2H,x,t)|\leq C |(x,t)|^{k}, ~~\forall ~(x,t)\in \Omega_{1},
\end{equation*}
where $C$ depends only on $k,n,\lambda, \Lambda,\alpha,K_0$ and $\omega$.
\end{theorem}

\proof As before, to prove that $u$ is $C^{k+1,\alpha}$ at $0$, we only need to prove the following. There exists a sequence of $k+1$-forms $H_m$ ($m\geq 0$) such that for all $m\geq 1$,

\begin{equation}\label{e.tCkas-6}
\|u-H_m\|_{L^{\infty }(\Omega _{\eta^{m}})}\leq \eta ^{m(k+1+\alpha )},
\end{equation}
\begin{equation}\label{e.tCkas-7}
\|H_m-H_{m-1}\|\leq \tilde{C}\eta ^{(m-1)\alpha}
\end{equation}
and
\begin{equation}\label{e.tCkas-9}
|(H_m)_t-F_0(D^2H_m,x,t)|\leq \tilde{C}|(x,t)|^{k}, ~~\forall ~(x,t)\in \Omega_{1},
\end{equation}
where $\tilde{C}$ and $\eta$ depend only on $k,n,\lambda, \Lambda,\alpha,K_0$ and $\omega$.

We prove the above by induction. For $m=1$, by \Cref{L.Cka.2}, there exists a $k+1$-form $H_1$ such that \crefrange{e.tCkas-6}{e.tCkas-9} hold for some $C_1$ and $\eta_1$ depending only on $k,n,\lambda, \Lambda,\alpha,K_0$ and $\omega$ where $H_0\equiv 0$. Take $\tilde{C}\geq C_1, \eta\leq \eta_1$ and then the conclusion holds for $m=1$. Suppose that the conclusion holds for $m=m_0$. We need to prove that the conclusion holds for $m=m_0+1$.

Let $r=\eta ^{m_{0}}$, $y=x/r$, $s=t/r^2$ and
\begin{equation}\label{e.tCkas-v}
  v(y,s)=\frac{u(x,t)-H_{m_0}(x,t)}{r^{k+1+\alpha}}.
\end{equation}
Then $v$ satisfies
\begin{equation}\label{e.Ckas-F}
\left\{\begin{aligned}
&v_s-\tilde{F}(D^2v,y,s)=\tilde{f}&& ~~\mbox{in}~~\tilde{\Omega}\cap Q_1;\\
&v=\tilde{g}&& ~~\mbox{on}~~\partial \tilde{\Omega}\cap Q_1,
\end{aligned}\right.
\end{equation}
where for $(M,y,s)\in \mathcal {S}^n\times \bar{\tilde\Omega}_1$,
\begin{equation*}
  \begin{aligned}
&\tilde{F}(M,y,s)=\frac{F(r^{k-1+\alpha}M+D^2H_{m_0},x,t)-(H_{m_0})_t}{r^{k-1+\alpha}},\\
&\tilde{f}(y,s)=\frac{f(x,t)}{r^{k-1+\alpha}},~
\tilde{g}(y,s)=\frac{g(x,t)-H_{m_0}(x,t)}{r^{k+1+\alpha}}~~\mbox{ and }~~\tilde{\Omega}=\frac{\Omega}{r}.\\
  \end{aligned}
\end{equation*}
In addition, define $\tilde{F}_0$ in a similar way to the definition of $\tilde{F}$ (only replacing $F$ by $F_0$).

In the following, we show that \cref{e.Ckas-F} satisfies the assumptions of \Cref{L.Cka.2}. First, by \cref{e.tCkas-7}, there exists a constant $C_{0}$ depending only on $k,n,\lambda, \Lambda,\alpha,K_0$ and $\omega$ such that $\|H_m\|\leq C_{0}$ ($\forall~0\leq m\leq m_0$). Then it is easy to verify that (note that $H_{m_0}$ is a $k+1$-form)
\begin{equation*}
\begin{aligned}
\|v\|_{L^{\infty}(\tilde{\Omega}_1)}&\leq 1, v(0)=\cdots=|D^{k}v(0)|=0~(\mathrm{by}~ \cref{e.tCkas-be},~ \cref{e.tCkas-6} ~\mbox{and}~\cref{e.tCkas-v}),\\
|\tilde{f}(y,s)|=& \frac{1}{r^{k-1+\alpha}}|f(x,t)|
\leq \frac{1}{r^{k-1+\alpha}} \delta |(x,t)|^{k-1+\alpha}
= \delta |(y,s)|^{k-1+\alpha},\\
|\tilde{g}(y,s)|=& \frac{1}{r^{k+1+\alpha}}(|g(x,t)|+|H_{m_0}(x,t)|)\\
\leq& \frac{1}{r^{k+1+\alpha}} \left(\frac{\delta}{2}|(x,t)|^{k+1+\alpha}
+\frac{\delta}{2C_0}\cdot C_{0}|(x,t)|^{k+1+\alpha}\right)\leq \delta |(y,s)|^{k+1+\alpha},\\
\|\partial \tilde{\Omega}\cap Q_1&\|_{C^{1,\alpha}(0)}\leq r^{\alpha}\|\partial\Omega\cap Q_1\|_{C^{1,\alpha}(0)}\leq \delta.
  \end{aligned}
\end{equation*}

Next, $\tilde{F}$ and $\tilde{F}_0$ are uniformly elliptic with constants $\lambda$ and $\Lambda$ obviously. In addition, for $(M,y,s)\in \mathcal{S}^{n}\times \bar{\tilde{\Omega}}_1$,
\begin{equation}\label{e.9.1}
\begin{aligned}
|&\tilde{F}(M,y,s)-\tilde{F}_0(M,y,s)|\\
&\leq \frac{1}{r^{k-1+\alpha}}\left(\beta(x,t)\left(r^{k-1+\alpha}\|M\|+C_0\right)
+\gamma(x,t)\right)\\
&= \beta(x,t)\|M\|+\frac{1}{r^{k-1+\alpha}}\left(C_0\beta(x,t)+\gamma(x,t)\right)\\
&\leq \beta(x,t)\|M\|+\frac{1}{r^{k-1+\alpha}}\left(C_0\cdot\frac{\delta}{2C_0}|(x,t)|^{k-1+\alpha}
+\frac{\delta}{2}|(x,t)|^{k-1+\alpha}\right)\\
&:= \beta(y,s)\|M\|+\tilde{\gamma}(y,s).
\end{aligned}
\end{equation}
Thus,
\begin{equation*}
\|\beta\|_{C^{k-1,\alpha}(0)}\leq \delta
~~~~\mbox{ and }~~\|\tilde{\gamma}\|_{C^{k-1,\alpha}(0)}\leq \delta.
\end{equation*}

Now, we prove that $\tilde{F}_0$ satisfies \cref{e1.8}. Indeed, for any $M_1,M_2\in \mathcal{S}^{n}$ and $(y_1,s_1),(y_2,s_2)\in \bar{\tilde{\Omega}}\cap \bar Q_{1}$,
\begin{equation*}
  \begin{aligned}
|\tilde{F}&_{0,ij}(M_1,y_1,s_1)-\tilde{F}_{0,ij}(M_2,y_2,s_2)|\\
=&|F_{0,ij}(r^{k-1+\alpha}M_1+D^2H_{m_0},x_1,t_1)
-F_{0,ij}(r^{k-1+\alpha}M_2+D^2H_{m_0},x_2,t_2)|\\
\leq &K_0\left(\|M_1-M_2\|+\|D^2H_{m_0}(x_1,t_1)-D^2H_{m_0}(x_2,t_2)\|
+|(x_1,t_1)-(x_2,t_2)|\right)\\
\leq &\tilde{K}_0\left(\|M_1-M_2\|+|(y_1,s_1)-(y_2,s_2)|\right),
  \end{aligned}
\end{equation*}
where $\tilde{K}_0$ depends only on $K_0$ and $C_0$.

Finally, with the aid of \cref{e.tCkas-9}, we show that
\begin{equation}\label{e.F0.omega}
\|\tilde{F}_0\|_{C^{k-1,1}(\bar{\textbf{B}}_\rho\times \bar{\tilde\Omega}_1)}\leq \tilde{\omega}(\rho),~\forall ~\rho>0,
\end{equation}
where $\tilde{\omega}$ depends only on $k,n,\lambda, \Lambda,\alpha,K_0$ and $\omega$. Note that
\begin{equation}
  \begin{aligned}
&\tilde{F}_0(M,y,s)\\
&=\tilde{F}_0(M,y,s)-\tilde{F}_0(0,y,s)+\tilde{F}_0(0,y,s)\\
    &=\int_{0}^{1} F_{0,ij}\left(\tau r^{k-1+\alpha} M+D^2H_{m_0},x,t\right)M_{ij}d\tau
    +\frac{F_0(D^2H_{m_0},x,t)-(H_{m_0})_t}{r^{k-1+\alpha}}\\
    &:=G_1(M,y,s)+G_2(y,s).
  \end{aligned}
\end{equation}
For $G_1$, since $F_0\in C^{k-1,1}$, $r,\tau\leq 1$ and $\|H_{m_0}\|\leq C_0$,
\begin{equation*}
\|D^{i}G_1\|_{L^{\infty}(\bar{\textbf{B}}_\rho\times \bar{\tilde\Omega}_1)}\leq \tilde{\omega}(\rho),~~\forall ~1\leq i\leq k-2.
\end{equation*}
From \cref{e.tCkas-9},
\begin{equation*}
|G_2(y,s)|\leq \tilde{C}|(y,s)|^k.
\end{equation*}
By combining with $G_2\in C^{k-1,1}$,
\begin{equation*}
\|D^{i}G_2\|_{L^{\infty}(\bar{\tilde\Omega}_1)}\leq C, ~~\forall ~ 1\leq i\leq k-1,
\end{equation*}
where $C$ depends only on $k,n,\lambda, \Lambda,\alpha,K_0$ and $\omega$. Hence, $\|\tilde{F}_0\|_{C^{k-2}(\bar{\textbf{B}}_\rho\times \bar{\tilde\Omega}_1)}$ has a uniform bound (independent of $m$).

From the definition of $\tilde{F}_0$, if any $(k-1)$-th derivative of $\tilde{F}_0$ involves one derivative with respect to $M$, it is bounded. Thus, we only need to take care of the $(k-1)$-th derivative with respect to $y$ or $s$. Note that $D^{k-1}F_0\in C^{0,1}$ and \cref{e.tCkas-9}. Then for $(M,y,s)\in \bar{\textbf{B}}_\rho\times \bar{\tilde\Omega}_1$,
\begin{equation*}
  \begin{aligned}
D^{k-1}_{\tilde{p}}&\tilde{F}_0(M,y,s)\\
=&D^{k-1}_{\tilde{p}}\left(\tilde{F}_0(M,y,s)-\tilde{F}_0(0,y,s)\right)
+D^{k-1}_{\tilde{p}}\tilde{F}_0(0,y,s)\\
=&r^{-\alpha}D^{k-1}_{p}\Big(F_0(r^{k-1+\alpha}M+D^2H_{m_0},x,t)
-F_0(D^2H_{m_0},x,t)\Big)\\
&+r^{-\alpha}D^{k-1}_{p}\left(F_0(D^2H_{m_0},x,t)-(H_{m_0})_t\right)\\
\leq & Cr^{-\alpha}\cdot r^{k-1+\alpha}\|M\|+C\\
\leq& \tilde{\omega}(\rho),
  \end{aligned}
\end{equation*}
where $\tilde{\omega}$ depends only on $k,n,\lambda, \Lambda,\alpha,K_0$ and $\omega$. That is, $\|\tilde{F}_0\|_{C^{k-1}(\bar{\textbf{B}}_\rho\times \bar{\tilde\Omega}_1)}$ has a uniform bound.

Similarly,
\begin{equation*}
  \begin{aligned}
&\left|D^{k-1}_{\tilde{p}}\tilde{F}_0(M,y_1,s_1)-D^{k-1}_{\tilde{p}}\tilde{F}_0(M,y_2,s_2)\right|\\
=&r^{-\alpha}\Big|D^{k-1}_{p}F_0(r^{k-1+\alpha}M+D^2H_{m_0},x_1,t_1)
-D^{k-1}_{p}F_0(r^{k-1+\alpha}M+D^2H_{m_0},x_2,t_2)\Big|\\
\leq & Cr^{-\alpha}\cdot r|(y_1,s_1)-(y_2,s_2)|\\
\leq& \tilde{\omega}(\rho)|(y_1,s_1)-(y_2,s_2)|.
  \end{aligned}
\end{equation*}
Thus, $\|\tilde{F}_0\|_{C^{k-1,1}(\bar{\textbf{B}}_\rho\times \bar Q_1)}$ is uniformly bounded. Therefore, \cref{e.F0.omega} holds.

Since \cref{e.Ckas-F} satisfies the assumptions of \Cref{L.Cka.2} (with $\omega$ replaced by $\tilde{\omega}$), there exists a $k+1$-form $\bar{H}$, constants $\tilde{C}\geq C_1$ and $\eta\leq \eta_1$ depending only on $k,n,\lambda, \Lambda,\alpha,K_0$ and $\omega$ such that
\begin{equation*}
\begin{aligned}
    \|v-\bar{H}\|_{L^{\infty }(\tilde{\Omega} _{\eta})}&\leq \eta^{k+1+\alpha},
\end{aligned}
\end{equation*}
\begin{equation*}
\|\bar{H}\|\leq \tilde{C}
\end{equation*}
and
\begin{equation*}
  |(\bar{H})_t-\tilde F_0(D^2\bar{H},y,s)|\leq \tilde{C}|(y,s)|^{k}
  , ~~\forall ~(y,s)\in \tilde{\Omega}_{1}.
\end{equation*}

Let $H_{m_0+1}(x,t)=H_{m_0}(x,t)+r^{\alpha}\bar{H}(x,t)$. Then \cref{e.tCkas-7} and \cref{e.tCkas-9} hold for $m_0+1$. Recalling \cref{e.tCkas-v}, we have
\begin{equation*}
  \begin{aligned}
\|u-H_{m_0+1}(x,t)\|_{L^{\infty}(\Omega_{\eta^{m_0+1}})}&= \|u-H_{m_0}(x,t)-r^{\alpha}\bar{H}(x,t)\|_{L^{\infty}(\Omega_{\eta r})}\\
&= \|r^{k+1+\alpha}v-r^{k+1+\alpha}\bar{H}(y,s)\|_{L^{\infty}(\tilde{\Omega}_{\eta})}\\
&\leq r^{k+1+\alpha}\eta^{k+1+\alpha}=\eta^{(m_0+1)(k+1+\alpha)}.
  \end{aligned}
\end{equation*}
Hence, \cref{e.tCkas-6} holds for $m=m_0+1$. By induction, the proof is completed.\qed~\\

Now, we give the~\\
\noindent\textbf{Proof of \Cref{t.Cka.1}.} As before, in the following proof,  we just make some necessary normalizations to satisfy the conditions of \Cref{t.Cka.S}. Throughout this proof, $C$ always denotes a constant depending only on $k,n,\lambda, \Lambda,\alpha,K_0,\omega$, $\|\beta\|_{C^{k-1,\alpha}(0)}$, $\|(\partial \Omega)_1\|_{C^{1,\alpha}(0)}$, $\|u\|_{L^{\infty}(\Omega_1)}$, $\|f\|_{C^{k-1,\alpha}(0)}$, $\|g\|_{C^{k+1,\alpha}(0)}$  and $\|\gamma\|_{C^{k-1,\alpha}(0)}$.

Let $F_1(M,x,t)=F(M,x,t)+P_f(x,t)$ for $(M,x,t)\in \mathcal {S}^{n}\times \bar{\Omega}_1$. Then $u$ satisfies
\begin{equation*}
\left\{\begin{aligned}
&u_t-F_1(D^2u,x,t)=f_1&& ~~\mbox{in}~~\Omega_1;\\
&u=g&& ~~\mbox{on}~~(\partial \Omega)_1,
\end{aligned}\right.
\end{equation*}
where $f_1(x,t)=f(x,t)-P_f(x,t)$. Thus,
\begin{equation*}
  |f_1(x,t)|\leq C|(x,t)|^{k-1+\alpha}, ~~\forall ~(x,t)\in \Omega_1.
\end{equation*}

Next, set $u_1(x,t)=u(x,t)-P_g(x,t)$ and $F_2(M,x,t)=F_1(M+D^2P_g,x,t)-(P_g)_t$. Then $u_1$ is a viscosity solution of
\begin{equation*}
\left\{\begin{aligned}
&(u_1)_t-F_2(D^2u_1,x,t)=f_1&& ~~\mbox{in}~~\Omega_1;\\
&u_1=g_1&& ~~\mbox{on}~~(\partial \Omega)_1,
\end{aligned}\right.
\end{equation*}
where $g_1(x,t)=g(x,t)-P_g(x,t)$. Hence,
\begin{equation}\label{e.cka.g}
  |g_1(x,t)|\leq C|(x,t)|^{k+1+\alpha}, ~~\forall ~(x,t)\in (\partial \Omega)_1.
\end{equation}

Finally, take $y=x/\rho$, $s=t/\rho^2$ and $u_2(y,s)=u_1(x,t)/\rho^2$, where $0<\rho<1$ is a constant to be specified later. Then $u_2$ satisfies
\begin{equation}\label{e.Cka.FF}
\left\{\begin{aligned}
&(u_2)_s-F_3(D^2u_2,y,s)=f_2&& ~~\mbox{in}~~\tilde{\Omega}_1;\\
&u_2=g_2&& ~~\mbox{on}~~(\partial \tilde{\Omega})_1,
\end{aligned}\right.
\end{equation}
where
\begin{equation*}
  \begin{aligned}
&F_3(M,y,s)=F_2(M,x,t),\\
&f_2(y,s)=f_1(x,t),~ g_2(y,s)=\frac{g_1(x,t)}{\rho^2}~~\mbox{ and }~~ \tilde{\Omega}=\frac{\Omega}{\rho}.
\end{aligned}
\end{equation*}
Define the fully nonlinear operator $F_{30}$ in the same way as $F_3$ (only replacing $F$ by $F_0$).

Now, we try to choose a proper $\rho$ such that \cref{e.Cka.FF} satisfies the conditions of \Cref{t.Cka.S}. First, by induction, \Cref{t.Cka.1} holds for $k-1$. Hence, by applying \Cref{t.Cka.1} for $k-1$ to $u_1$ (note $k\geq 2$ and $u_1(0)=\cdots=|D^{k}u_1(0)|=0$),
\begin{equation*}
  \begin{aligned}
    \|u_2\|_{L^{\infty}(\Omega_1)}&= \|u_1\|_{L^{\infty}(\Omega_\rho)}/\rho^2
    \leq C\rho^{k-2+\alpha}.
  \end{aligned}
\end{equation*}
Clearly, by the assumption,
\begin{equation*}
  u_2(0)=|Du_2(0)|=\cdots=|D^{k}u_2(0)|=|Dg_2(0)|=\cdots=|D^{k}g_2(0)|=0.
\end{equation*}
Next,
\begin{equation*}
\begin{aligned}
&|f_2(y,s)|=|f_1(x,t)|
\leq C\rho^{k-1+\alpha}
|(y,s)|^{k-1+\alpha},~~~~\forall~~(y,s)\in\tilde{\Omega}_1,\\
&|g_3(y,s)|=\frac{|g_1(x,t)|}{\rho^2}\leq C\rho^{k-1+\alpha}|(y,s)|^{k+1+\alpha},~~~~\forall~~(y,s)\in\partial \tilde{\Omega}_1,\\
&\|\partial \tilde{\Omega}\cap Q_1\|_{C^{1,\alpha}(0)}\leq \rho^{\alpha}\|\partial \Omega\cap Q_1\|_{C^{1,\alpha}(0)}.
\end{aligned}
\end{equation*}

The $F_3$ and $F_{30}$ are uniformly elliptic with $\lambda$ and $\Lambda$ and convex in $M$ obviously. Next, we show that $F_3\in C^{k-1,\alpha}(0)$. Indeed,
\begin{equation*}
  \begin{aligned}
&|F_3(M,y,s)-F_{30}(M,y,s)|\\
&=\left|F\left(M+D^2P_g,x,t\right)
-F_{0}\left(M+D^2P_g,x,t\right) \right|\\
&\leq  \beta(x,t)\left(\|M\|+C \right)+\gamma(x,t)\\
&=\beta(x,t)\|M\|+C\beta(x,t)+\gamma(x,t)\\
&\leq C\rho^{k-1+\alpha} |(y,s)|^{k-1+\alpha}\|M\|
+C\rho^{k-1+\alpha}|(y,s)|^{k-1+\alpha}\\
&:=\tilde{\beta}(y,s)\|M\|+\tilde{\gamma}(y,s).
\end{aligned}
\end{equation*}
In addition, for any $M_1,M_2\in \mathcal{S}^{n}$ and $(y_1,s_1),(y_2,s_2)\in \bar\Omega\cap \bar Q_{1}$,
\begin{equation*}
  \begin{aligned}
|F&_{30,ij}(M_1,y_1,s_1)-F_{30,ij}(M_2,y_2,s_2)|\\
=&|F_{0,ij}(M_1+D^2P_g(x_1,t_1),x_1,t_1)-F_{0,ij}(M_2+D^2P_g(x_2,t_2),x_2,t_2)|\\
\leq& K_0 (\|M_1-M_2\|+|D^2P_g(x_1,t_1)-D^2P_g(x_2,t_2)|+|(x_1,t_1)-(x_2,t_2)|)\\
\leq &\tilde{K}_0(\|M_1-M_2\|+|(x_1,t_1)-(x_2,t_2)|)\\
\leq &\tilde{K}_0(\|M_1-M_2\|+|(y_1,s_1)-(y_2,s_2)|),
  \end{aligned}
\end{equation*}
where $\tilde{K}_0$ depends only on $K_0$ and $\|g\|_{C^{k+1,\alpha}(0)}$.

Finally,
\begin{equation*}
  \begin{aligned}
    &\|F_{30}\|_{C^{k-1,1}(\bar{\textbf{B}}_r\times \bar{\tilde\Omega}_1)}\\
&=\left \| F_0\left(M +D^2P_g,x,t\right)
+P_f-(P_g)_t  \right \|_{C^{k-1,1}(\bar{\textbf{B}}_r\times \bar{\tilde\Omega}_1)}\\
&\leq \tilde{\omega}(r),~\forall ~r>0,
  \end{aligned}
\end{equation*}
where $\tilde{\omega}$ depends only on $k,n,\omega$, $\|f\|_{C^{k-1,\alpha}(0)}$ and $\|g\|_{C^{k+1,\alpha}(0)}$.

From above arguments, we can choose $\rho$ small enough (depending only on $k,n,\lambda, \Lambda,\alpha,K_0,\omega$, $\|\beta\|_{C^{k-1,\alpha}(0)}$, $\|(\partial \Omega)_1\|_{C^{1,\alpha}(0)}$, $\|u\|_{L^{\infty}(\Omega_1)}$, $\|f\|_{C^{k-1,\alpha}(0)}$, $\|g\|_{C^{k+1,\alpha}(0)}$  and $\|\gamma\|_{C^{k-1,\alpha}(0)}$) such that the assumptions of \Cref{t.Cka.S} are satisfied. Then $u_2$ and hence $u$ is $C^{k+1,\alpha}$ at $0$, and the estimates \crefrange{e.Cka.u1}{e.Cka.u2} hold. \qed~\\

~\\

Next, we prove \Cref{Ckla} based on \Cref{t.Cka.1}.\\
\noindent\textbf{Proof of \Cref{Ckla}.} Throughout this proof, $C$ always denotes a constant depending only on $ k, l, n, \lambda, \Lambda, \alpha,K_0,\omega$, $\|\beta\|_{C^{k+l-2,\alpha}(0)}$, $\|(\partial \Omega)_1\|_{C^{l,\alpha}(0)}$, $\|u\|_{L^{\infty}(\Omega_1)}$, $\|f\|_{C^{k+l-2,\alpha}(0)}$, $\|g\|_{C^{k+l,\alpha}(0)}$ and $\|\gamma\|_{C^{k+l-2,\alpha}(0)}$.

Since $g\in C^{k+l,\alpha}(0)$,
\begin{equation*}
  |g(x,t)-P_g(x,t)|\leq C|(x,t)|^{k+l+\alpha},~~~\forall ~(x,t)\in (\partial \Omega)_1.
\end{equation*}
Set $u_1=u-P_g$ and then $u_1$ is a viscosity solution of
\begin{equation*}
\left\{\begin{aligned}
&(u_1)_t-F_1(D^2u_1,x,t)=f&& ~~\mbox{in}~~\Omega_1;\\
&u_1=g_1&& ~~\mbox{on}~~(\partial \Omega)_1,
\end{aligned}\right.
\end{equation*}
where $F_1(M,x,t)=F(M+D^2P_g,x,t)-(P_g)_t$ and $g_1=g-P_g$. As before, $F_{10}$ is defined similarly to $F_1$. Hence,
\begin{equation}\label{e.Ckla.gnorm}
u_1(0)=|Du_1(0)|=\cdots|D^ku_1(0)|=|Dg_1(0)|=\cdots=|D^{k+l}g_1(0)|=0.
\end{equation}
In addition, it can be verified as before that $F_1$ and $F_{10}$ satisfies the conditions of \Cref{t.Cka.1}.
By \Cref{t.Cka.1}, $u_1\in  C^{k+1,\alpha}(0)$ and there exists a $(k+1)$-form $H_{k+1}$ such that
\begin{equation*}
|u_1(x,t)-H_{k+1}(x,t)|\leq C |(x,t)|^{k+1+\alpha},~~\forall ~(x,t)\in \Omega_{1}.
\end{equation*}

Let
\begin{equation*}
  u_2(x,t)=u_1(x,t)-H_{k+1}(x',x_n-P_{\Omega}(x',t),t),
\end{equation*}
where $P_{\Omega}\in \mathcal{P}_l$ corresponds to $\partial \Omega$ at $0$ (note that $\partial \Omega\in C^{l,\alpha}(0)$). Since $H_{k+1}$ is a $(k+1)$-form and $P_{\Omega}(0)=|DP_{\Omega}(0)|=0$,
\begin{equation*}
D^{k+1}(H_{k+1}(x',x_n-P_{\Omega}(x',t),t))(0)=D^{k+1}(H_{k+1}(x,t))(0).
\end{equation*}
Hence, $u_2(0)=\cdots= |D^{k+1}u_2(0)|=0$. In addition, $u_2$ satisfies
\begin{equation*}
\left\{\begin{aligned}
&(u_2)-F_2(D^2u_2,x,t)=f&& ~~\mbox{in}~~\Omega_1;\\
&u_2=g_2&& ~~\mbox{on}~~(\partial \Omega)_1,
\end{aligned}\right.
\end{equation*}
where $F_2(M,x,t)=F_1\left(M+D^2 H_{k+1},x,t\right)-(H_{k+1})_t$ and $g_2=g_1-H_{k+1}$. Thus,
\begin{equation*}
  |g_2(x,t)|\leq |g_1(x,t)|+|H_{k+1}(x',x_n-P_{\Omega}(x',t),t)|\leq C|(x,t)|^{k+l+\alpha},~~~\forall ~(x,t)\in (\partial \Omega)_1.
\end{equation*}
As before, $F_2 \in C^{k+l-2,\alpha}(0)$. By \Cref{t.Cka.1} again, $u_2\in C^{k+2,\alpha}(0)$ and there exists a $(k+2)$-form $H_{k+2}$ such that
\begin{equation*}
|u_2(x,t)-H_{k+2}(x,t)|\leq C |(x,t)|^{k+2+\alpha},~~\forall ~(x,t)\in \Omega_{1}.
\end{equation*}

Let
\begin{equation}\label{e.Ckla.u}
  \begin{aligned}
u_3(x,t)&=u_2(x,t)-H_{k+2}(x',x_n-P_{\Omega}(x',t),t)\\
  &=u_1(x,t)-H_{k+1}(x',x_n-P_{\Omega}(x',t),t)-H_{k+2}(x',x_n-P_{\Omega}(x',t),t).
  \end{aligned}
\end{equation}
Then $u_3(0)=\cdots=|D^{k+2}u_3(0)|=0$ and $u_3$ is a viscosity solution of
\begin{equation*}
\left\{\begin{aligned}
&(u_3)_t-F_3(D^2 u_3,x,t)=f&& ~~\mbox{in}~~\Omega_1;\\
&u_3=g_3&& ~~\mbox{on}~~(\partial \Omega)_1,
\end{aligned}\right.
\end{equation*}
where $F_3(M,x,t)=F_2\left(M+D^2 H_{k+2},(x,t)\right)-(H_{k+2})_t$ and $g_3=g_2-H_{k+2}$. Then,
\begin{equation*}
  |g_3(x,t)|\leq |g_2(x,t)|+|H_{k+2}(x',x_n-P_{\Omega}(x',t),t)|\leq C|(x,t)|^{k+l+\alpha},~~~\forall ~(x,t)\in (\partial \Omega)_1
\end{equation*}
and $F_3 \in C^{k+l-2,\alpha}(0)$.

By \Cref{t.Cka.1} again, $u_3\in C^{k+2,\alpha}(0)$ and hence $u\in C^{k+2,\alpha}(0)$. By similar arguments again and again, $u\in C^{k+l,\alpha}(0)$ eventually and \cref{e1.2} holds. Therefore, the proof of \Cref{Ckla} is completed.\qed~\\

\begin{remark}\label{re.Ckla.1}
By checking the proof, we know that the polynomial $P$ can be written as \cref{e1.5} (see \cref{e.Ckla.u}).
\end{remark}
~\\

Finally, we give the\\
\noindent\textbf{Proof of \Cref{Cka}.} As before, we only need to make some normalizations. In this proof, $C$ always denotes a constant depending only on $k, n, \lambda, \Lambda, \alpha,K_0,\omega$, $\|\beta\|_{C^{k-2,\alpha}(0)}$, $\|(\partial \Omega)_1\|_{C^{k,\alpha}(0)}$, $\|u\|_{L^{\infty}(\Omega_1)}$, $\|f\|_{C^{k-2,\alpha}(0)}$, $\|g\|_{C^{k,\alpha}(0)}$ and $\|\gamma\|_{C^{k-2,\alpha}(0)}$.

First, by \Cref{C2a}, $u\in C^{2,\tilde\alpha}(0)$ for some $0<\tilde{\alpha}<\bar{\alpha}$ and thus $u\in C^{1,\alpha}(0)$. Let
\begin{equation*}
v(x,t)=u(x,t)-P_g(x,t)+((P_{g})_n(0)-u_n(0))\left(x_n-P_{\Omega}(x',t)\right)
\end{equation*}
and
\begin{equation*}
\tilde{F}(M,x,t)=F(M+D^2P_g+((P_{g})_n(0)-u_n(0))D^2P_{\Omega},x,t)-(P_{g})_t
-((P_{g})_n(0)-u_n(0))(P_{\Omega})_t
\end{equation*}
for $(M,x,t) \in \mathcal{S}^n\times \Omega_1$.
Then $v$ satisfies
\begin{equation*}
\left\{\begin{aligned}
&v_t-\tilde{F}(D^2v,x,t)=f&& ~~\mbox{in}~~\Omega_1;\\
&v=\tilde{g}&& ~~\mbox{on}~~(\partial \Omega)_1,
\end{aligned}\right.
\end{equation*}
where
\begin{equation*}
  \tilde{g}(x,t)=g(x,t)-P_g(x,t)+((P_{g})_n(0)-u_n(0))\left(x_n-P_{\Omega}(x',t)\right).
\end{equation*}
By \cref{e.C2a.3-1} and noting that $P_{\Omega}\in \mathcal{P}_{k-1}$ with $P_{\Omega}(0)=|DP_{\Omega}(0)|=0$,
\begin{equation*}
v(0)=|Dv(0)|=|D\tilde g(0)|=0.
\end{equation*}
In addition, since $g\in C^{k,\alpha}(0)$ and $(\partial \Omega)_1\in C^{k,\alpha}(0)$,
$\tilde{g}\in C^{k,\alpha}(0)$. By \Cref{Ckla}, $v$ and hence $u\in C^{k,\alpha}(0)$. \qed~\\

\subsection{Application to the regularity of free boundaries}

Finally, we prove the higher regularity of free boundaries in obstacle-type problems with the aid of the boundary pointwise regularity.\\
\noindent\textbf{Proof of \Cref{FreeBd}.} For the interior regularity, it is easy to show that $u\in C^{\infty}(\Omega_1)$. In the following, we only need to take care of the regularity up to the boundary.

Assume that
\begin{equation*}
  \partial \Omega\cap Q_1=\left\{(x',x_n,t)\big| x_n=\varphi(x',t)\right\},
\end{equation*}
where $\varphi \in C^{1,\alpha}(S_1)$ and $\varphi(0)=|D\varphi(0)|=0$. By setting $(x_2,t_2)=0$ in \cref{e1.10}, $F$ satisfies \cref{e.C2a.beta} at $0$. In addition, note that $u=|Du|=0$ on $(\partial \Omega)_1$ and $(\partial \Omega)_1\in C^{1,\alpha}$. Then by \Cref{C2a-1}, $u\in C^{2,\alpha}(0)$. Similarly, $u\in C^{2,\alpha}(x_0,t_0)$ for any $(x_0,t_0)\in (\partial \Omega)_1$. By combining with the interior regularity, $u\in C^{2,\alpha}(\bar\Omega')$ for any $\Omega'\subset\subset \bar\Omega\cap Q_1$.

In the following proof, we will use the higher regularity \Cref{Ckla}. Since $u\in C^{2,\alpha}(\bar\Omega')$ for any $\Omega'\subset\subset \bar\Omega\cap Q_1$, $D^2u(\bar\Omega')$ is a compact set in $\mathcal{S}^n$. Note that $F$ is a smooth function. Hence, by the mean value theorem, there exists a constant $C$ (depending on $D^2u(\bar\Omega')$) such that for any $M_1,M_2\in D^2u(\bar\Omega')$ and $(x_1,t_1),(x_2,t_2)\in \bar\Omega\cap Q_1$,
\begin{equation}\label{e4.4}
|F_{ij}(M_1,x_1,t_1)-F_{ij}(M_2,x_2,t_2)|\leq C(\|M_1-M_2\|+|(x_1,t_1)-(x_2,t_2)|).
\end{equation}
Then we can redefine $F$ outside $D^2u(\bar\Omega')$ such that $F$ satisfies \cref{e4.4} for any $M\in \mathcal{S}^n$. Clearly, $u$ is still the solution of the modified $F$ in $\Omega'$. Therefore, we can and do assume that $F$ satisfies \cref{e4.4} (i.e. \cref{e1.8}) in the following proof.

From $|D\varphi(0)|=0$ and $u=|Du|=0$ on $(\partial \Omega)_1$ again, $u_{ij}(0)=0$ for $i+j<2n$ and $u_t(0)=0$ (see \cref{e1.5} with $P_g\equiv 0$ and $k=l=1$ and we have $u_t=0$ on $(\partial \Omega)_1$). Hence, $F(u_{nn}(0),0,0)=-1$ by the equation in \cref{e.FBd.1}. By virtue of the uniform ellipticity and $F(0,0,0)=0$, $u_{nn}(0)=-c$ for some $c>0$. Let
\begin{equation*}
v(x,t)=u(x,t)-P(x),~~\mbox{ where }~~P(x)=\frac{1}{2}u_{nn}(0)x_n^2.
\end{equation*}
Then $v$ satisfies
\begin{equation*}
\left\{\begin{aligned}
&   v_t-\tilde{F} (D^2v,x,t)=1~~~~\mbox{in}~~\Omega_1;\\
&   v=g~~~~\mbox{on}~~(\partial \Omega)_1;\\
&   v(0)=|Dv(0)|=|D^2v(0)|=0,
\end{aligned}\right.
\end{equation*}
where $\tilde{F}(M,x,t)=F(M+D^2P,x,t)$ and $g(x,t)=-u_{nn}(0)x_n^2/2$. Then
\begin{equation*}
|g(x,t)|=\left|\frac{u_{nn}(0)}{2}x_n^2 \right|\leq \frac{c}{2}\|(\partial \Omega)_1\|_{C^{1,\alpha}(0)}^2|(x,t)|^{2+2\alpha},~\forall ~(x,t)\in (\partial \Omega)_1.
\end{equation*}
That is, $g\in C^{2,2\alpha}(0)$ and $g(0)=|Dg(0)|=|D^2g(0)|=0$. By \Cref{Ckla}, $v\in C^{2,2\alpha}(0)$ and hence $u\in C^{2,2\alpha}(0)$. Similarly, for any $(x_0,t_0)\in (\partial \Omega)_1$, $u\in C^{2,2\alpha}(x_0,t_0)$. Hence, $u\in C^{2,2\alpha}(\bar\Omega')$ for any $\Omega'\subset\subset \bar\Omega\cap Q_1$.

Since $u_{nn}(0)=-c$, $u_{nn}(0)\leq -c/2$ in $\Omega_r$ for some $r>0$. Then $u_i/u_n\in C^{2\alpha} (\bar\Omega\cap Q_{r})$ ($1\leq i\leq n-1$). Note that $\varphi_i=-u_i/u_n$. Thus, $\varphi\in C^{1,2\alpha}(S_r)$, i.e., $(\partial \Omega)_r\in C^{1,2\alpha}$. By considering other $(x_0, t_0)\in (\partial \Omega)_1$ similarly, we have $(\partial\Omega)_1\in C^{1,2\alpha}$.

Consider $v$ again and $g\in C^{2,4\alpha}(0)$ now (since $(\partial\Omega)_1\in C^{1,2\alpha}$). From \Cref{Ckla}, $v\in C^{2,4\alpha}(0)$. By similar arguments as the above, $u\in C^{2,4\alpha}(\bar\Omega')$ for any $\Omega'\subset\subset \bar{\Omega}\cap Q_1$. Therefore, $(\partial \Omega)_1\in C^{1, 4\alpha}$. Consider $v$ again and again and we have
$u\in C^{3,\tilde{\alpha}}(\bar\Omega')$ for any $\Omega'\subset\subset \bar{\Omega}\cap Q_1$ for some $0<\tilde{\alpha}<1$ eventually. Then $\varphi_i\in C^{1,\tilde{\alpha}}$ for $1\leq i\leq n-1$ as before. Moreover, $\varphi_t=-u_t/u_n\in C^{\tilde{\alpha}}$. Hence, $(\partial \Omega)_1\in C^{2,\tilde{\alpha}}$.

Let $v(x,t)=u(x,t)-P(x,t)$ where $P(x,t)=u_{nn}(0)\left(x_n-P_{\Omega}(x',t) \right)^2/2$. Here, $P_{\Omega}\in \mathcal{HP}_2$ is the polynomial corresponding to $\partial \Omega$ at $0$ since $(\partial \Omega)_1\in C^{2,\tilde{\alpha}}$. Then $v$ satisfies
\begin{equation*}
\left\{\begin{aligned}
&   v_t-\tilde{F}(D^2 v,x,t)=f~~~~\mbox{in}~~\Omega_1;\\
&   v=g~~~~\mbox{on}~~(\partial \Omega)_1;\\
&   v(0)=|Dv(0)|=|D^2v(0)|=0,
\end{aligned}\right.
\end{equation*}
where $f=1-P_t \in \mathcal{P}_2$, $\tilde{F}(M,x,t)=F(M+D^2 P,x,t)$ and $g(x,t)=-u_{nn}(0)(x_n-P_{\Omega}(x',t))^2/2$. Then
\begin{equation*}
|g(x,t)|=|u_{nn}(0)(x_n-P_{\Omega}(x',t))^2/2|\leq C|(x,t)|^{4+2\tilde{\alpha}},~\forall ~(x,t)\in (\partial \Omega)_1.
\end{equation*}
As before, by \Cref{Ckla}, $v\in C^{4, \tilde{\alpha}}(0)$ and hence $u\in C^{4, \tilde{\alpha}}(0)$. Similar to previous arguments, $u\in C^{4,\tilde{\alpha}}(\bar\Omega')$ for any $\Omega'\subset\subset \bar{\Omega}\cap Q_1$ and then $(\partial \Omega)_1\in C^{3,\tilde{\alpha}}$.

Let $v(x,t)=u(x,t)-u_{nn}(0)(x_n-P_{\Omega}(x',t))^2/2$ where $P_{\Omega}\in \mathcal{P}_3$ now. Repeat above arguments and we have $u\in C^{\infty}(\bar\Omega\cap Q_1)$ and $(\partial \Omega)_1\in C^{\infty}$ eventually. \qed~\\

\noindent\textbf{Data availability statement} Data sharing not applicable to this article as no datasets were generated or analysed during the current study.

%
%

%
\bibliographystyle{model4-names}
\bibliography{PDE}

\begin{thebibliography}{27}
\expandafter\ifx\csname natexlab\endcsname\relax\def\natexlab#1{#1}\fi
\providecommand{\url}[1]{\texttt{#1}}
\providecommand{\href}[2]{#2}
\providecommand{\path}[1]{#1}
\providecommand{\DOIprefix}{doi:}
\providecommand{\ArXivprefix}{arXiv:}
\providecommand{\URLprefix}{URL: }
\providecommand{\Pubmedprefix}{pmid:}
\providecommand{\doi}[1]{\href{http://dx.doi.org/#1}{\path{#1}}}
\providecommand{\Pubmed}[1]{\href{pmid:#1}{\path{#1}}}
\providecommand{\bibinfo}[2]{#2}
\ifx\xfnm\undefined \def\xfnm[#1]{\unskip,\space#1}\fi
\bibitem[{Adimurthi et~al.(2020)Adimurthi, Banerjee and Verma}]{MR4073515}
\bibinfo{author}{Adimurthi\xfnm[ K.]}, \bibinfo{author}{Banerjee\xfnm[ A.]},
  \bibinfo{author}{Verma\xfnm[ R.B.]}.
\newblock \bibinfo{title}{Twice differentiability of solutions to fully
  nonlinear parabolic equations near the boundary}.
\newblock \bibinfo{journal}{Nonlinear Anal}
  \bibinfo{year}{2020};\bibinfo{volume}{197}:\bibinfo{pages}{111830, 16}.
\newblock \URLprefix \url{https://doi.org/10.1016/j.na.2020.111830}.
  \DOIprefix\doi{10.1016/j.na.2020.111830}.
\bibitem[{Caffarelli(1989)}]{MR1005611}
\bibinfo{author}{Caffarelli\xfnm[ L.A.]}.
\newblock \bibinfo{title}{Interior a priori estimates for solutions of fully
  nonlinear equations}.
\newblock \bibinfo{journal}{Ann of Math (2)}
  \bibinfo{year}{1989};\bibinfo{volume}{130}(\bibinfo{number}{1}):\bibinfo{pages}{189--213}.
\newblock \URLprefix \url{http://dx.doi.org/10.2307/1971480}.
  \DOIprefix\doi{10.2307/1971480}.
\bibitem[{Caffarelli and Cabr\'{e}(1995)}]{MR1351007}
\bibinfo{author}{Caffarelli\xfnm[ L.A.]}, \bibinfo{author}{Cabr\'{e}\xfnm[
  X.]}.
\newblock \bibinfo{title}{Fully nonlinear elliptic equations}.
\newblock volume~\bibinfo{volume}{43} of \textit{\bibinfo{series}{American
  Mathematical Society Colloquium Publications}}.
\newblock \bibinfo{publisher}{American Mathematical Society, Providence, RI},
  \bibinfo{year}{1995}.
\newblock \URLprefix \url{https://doi.org/10.1090/coll/043}.
  \DOIprefix\doi{10.1090/coll/043}.
\bibitem[{Chen(2003)}]{Chen}
\bibinfo{author}{Chen\xfnm[ Y.Z.]}.
\newblock \bibinfo{title}{Second order parabolic partial differential equations
  ({C}hinese)}.
\newblock \bibinfo{publisher}{Peking University Press}, \bibinfo{year}{2003}.
\bibitem[{Crandall et~al.(1992)Crandall, Ishii and Lions}]{MR1118699}
\bibinfo{author}{Crandall\xfnm[ M.G.]}, \bibinfo{author}{Ishii\xfnm[ H.]},
  \bibinfo{author}{Lions\xfnm[ P.L.]}.
\newblock \bibinfo{title}{User's guide to viscosity solutions of second order
  partial differential equations}.
\newblock \bibinfo{journal}{Bull Amer Math Soc (NS)}
  \bibinfo{year}{1992};\bibinfo{volume}{27}(\bibinfo{number}{1}):\bibinfo{pages}{1--67}.
\newblock \URLprefix \url{https://doi.org/10.1090/S0273-0979-1992-00266-5}.
  \DOIprefix\doi{10.1090/S0273-0979-1992-00266-5}.
\bibitem[{Crandall et~al.(2000)Crandall, Kocan and \'{S}wi\polhk
  ech}]{MR1789919}
\bibinfo{author}{Crandall\xfnm[ M.G.]}, \bibinfo{author}{Kocan\xfnm[ M.]},
  \bibinfo{author}{\'{S}wi\polhk ech\xfnm[ A.]}.
\newblock \bibinfo{title}{{$L^p$}-theory for fully nonlinear uniformly
  parabolic equations}.
\newblock \bibinfo{journal}{Comm Partial Differential Equations}
  \bibinfo{year}{2000};\bibinfo{volume}{25}(\bibinfo{number}{11-12}):\bibinfo{pages}{1997--2053}.
\newblock \URLprefix \url{https://doi.org/10.1080/03605300008821576}.
  \DOIprefix\doi{10.1080/03605300008821576}.
\bibitem[{De~Silva and Savin(2015)}]{MR3393271}
\bibinfo{author}{De~Silva\xfnm[ D.]}, \bibinfo{author}{Savin\xfnm[ O.]}.
\newblock \bibinfo{title}{A note on higher regularity boundary {H}arnack
  inequality}.
\newblock \bibinfo{journal}{Discrete Contin Dyn Syst}
  \bibinfo{year}{2015};\bibinfo{volume}{35}(\bibinfo{number}{12}):\bibinfo{pages}{6155--6163}.
\newblock \URLprefix \url{https://doi.org/10.3934/dcds.2015.35.6155}.
  \DOIprefix\doi{10.3934/dcds.2015.35.6155}.
\bibitem[{Figalli and Shahgholian(2014)}]{MR3198649}
\bibinfo{author}{Figalli\xfnm[ A.]}, \bibinfo{author}{Shahgholian\xfnm[ H.]}.
\newblock \bibinfo{title}{A general class of free boundary problems for fully
  nonlinear elliptic equations}.
\newblock \bibinfo{journal}{Arch Ration Mech Anal}
  \bibinfo{year}{2014};\bibinfo{volume}{213}(\bibinfo{number}{1}):\bibinfo{pages}{269--286}.
\newblock \URLprefix \url{https://doi.org/10.1007/s00205-014-0734-0}.
  \DOIprefix\doi{10.1007/s00205-014-0734-0}.
\bibitem[{Kinderlehrer and Nirenberg(1977)}]{MR440187}
\bibinfo{author}{Kinderlehrer\xfnm[ D.]}, \bibinfo{author}{Nirenberg\xfnm[
  L.]}.
\newblock \bibinfo{title}{Regularity in free boundary problems}.
\newblock \bibinfo{journal}{Ann Scuola Norm Sup Pisa Cl Sci (4)}
  \bibinfo{year}{1977};\bibinfo{volume}{4}(\bibinfo{number}{2}):\bibinfo{pages}{373--391}.
\newblock \URLprefix
  \url{http://www.numdam.org/item?id=ASNSP_1977_4_4_2_373_0}.
\bibitem[{Krylov(1983)}]{MR688919}
\bibinfo{author}{Krylov\xfnm[ N.V.]}.
\newblock \bibinfo{title}{Boundedly inhomogeneous elliptic and parabolic
  equations in a domain}.
\newblock \bibinfo{journal}{Izv Akad Nauk SSSR Ser Mat}
  \bibinfo{year}{1983};\bibinfo{volume}{47}(\bibinfo{number}{1}):\bibinfo{pages}{75--108}.
\bibitem[{Li et~al.(2022)Li, Li and Zhang}]{MR_LiLiZhang}
\bibinfo{author}{Li\xfnm[ D.]}, \bibinfo{author}{Li\xfnm[ X.]},
  \bibinfo{author}{Zhang\xfnm[ K.]}.
\newblock \bibinfo{title}{${W^{2,p}}$ estimates for elliptic equations on
  ${C^{1,\alpha}}$ domains}.
\newblock \bibinfo{journal}{Math Ann, to appear}
  \bibinfo{year}{2022};\href{http://arxiv.org/abs/2201.01975}{\tt
  arXiv:2201.01975}.
\bibitem[{Li and Zhang(2018)}]{MR3780142}
\bibinfo{author}{Li\xfnm[ D.]}, \bibinfo{author}{Zhang\xfnm[ K.]}.
\newblock \bibinfo{title}{Regularity for fully nonlinear elliptic equations
  with oblique boundary conditions}.
\newblock \bibinfo{journal}{Arch Ration Mech Anal}
  \bibinfo{year}{2018};\bibinfo{volume}{228}(\bibinfo{number}{3}):\bibinfo{pages}{923--967}.
\newblock \URLprefix \url{https://doi.org/10.1007/s00205-017-1209-x}.
  \DOIprefix\doi{10.1007/s00205-017-1209-x}.
\bibitem[{Lian(2022)}]{MR_cones}
\bibinfo{author}{Lian\xfnm[ Y.]}.
\newblock \bibinfo{title}{Boundary pointwise regularity and liouville theorems
  for fully nonlinear equations on cones: arxiv:2205.14291}.
\newblock \bibinfo{year}{2022}.
\newblock \href{http://arxiv.org/abs/2205.14291}{\tt arXiv:2205.14291}.
\bibitem[{Lian et~al.(2020)Lian, Wang and Zhang}]{lian2020pointwise}
\bibinfo{author}{Lian\xfnm[ Y.]}, \bibinfo{author}{Wang\xfnm[ L.]},
  \bibinfo{author}{Zhang\xfnm[ K.]}.
\newblock \bibinfo{title}{Pointwise regularity for fully nonlinear elliptic
  equations in general forms}.
\newblock \bibinfo{year}{2020}.
\newblock \href{http://arxiv.org/abs/2012.00324}{\tt arXiv:2012.00324}.
\bibitem[{Lian and Zhang(2018)}]{lian2018boundary}
\bibinfo{author}{Lian\xfnm[ Y.]}, \bibinfo{author}{Zhang\xfnm[ K.]}.
\newblock \bibinfo{title}{Boundary lipschitz regularity and the hopf lemma for
  fully nonlinear elliptic equations}.
\newblock \bibinfo{year}{2018}.
\newblock \href{http://arxiv.org/abs/1812.11357}{\tt arXiv:1812.11357}.
\bibitem[{Lian and Zhang(2020)}]{MR4088470}
\bibinfo{author}{Lian\xfnm[ Y.]}, \bibinfo{author}{Zhang\xfnm[ K.]}.
\newblock \bibinfo{title}{Boundary pointwise {$C ^{1,\alpha}$} and
  {$C^{2,\alpha}$} regularity for fully nonlinear elliptic equations}.
\newblock \bibinfo{journal}{J Differential Equations}
  \bibinfo{year}{2020};\bibinfo{volume}{269}(\bibinfo{number}{2}):\bibinfo{pages}{1172--1191}.
\newblock \URLprefix \url{https://doi.org/10.1016/j.jde.2020.01.006}.
  \DOIprefix\doi{10.1016/j.jde.2020.01.006}.
\bibitem[{Lian and Zhang(2022)}]{freeboundary}
\bibinfo{author}{Lian\xfnm[ Y.]}, \bibinfo{author}{Zhang\xfnm[ K.]}.
\newblock \bibinfo{title}{Boundary pointwise regularity and applications to the
  regularity of free boundaries, arxiv:2204.09304}.
\newblock \bibinfo{year}{2022}.
\newblock \href{http://arxiv.org/abs/2204.09304}{\tt arXiv:2204.09304}.
\bibitem[{Nornberg(2019)}]{MR3980853}
\bibinfo{author}{Nornberg\xfnm[ G.]}.
\newblock \bibinfo{title}{{$C^{1,\alpha}$} regularity for fully nonlinear
  elliptic equations with superlinear growth in the gradient}.
\newblock \bibinfo{journal}{J Math Pures Appl (9)}
  \bibinfo{year}{2019};\bibinfo{volume}{128}:\bibinfo{pages}{297--329}.
\newblock \URLprefix \url{https://doi.org/10.1016/j.matpur.2019.06.008}.
  \DOIprefix\doi{10.1016/j.matpur.2019.06.008}.
\bibitem[{Pogorelov(1978)}]{MR0478079}
\bibinfo{author}{Pogorelov\xfnm[ A.V.]}.
\newblock \bibinfo{title}{The {M}inkowski multidimensional problem}.
\newblock \bibinfo{publisher}{V. H. Winston \& Sons, Washington, D.C.; Halsted
  Press [John Wiley \&\ Sons], New York-Toronto-London}, \bibinfo{year}{1978}.
\newblock \bibinfo{note}{Translated from the Russian by Vladimir Oliker,
  Introduction by Louis Nirenberg, Scripta Series in Mathematics}.
\bibitem[{Savin(2007)}]{MR2334822}
\bibinfo{author}{Savin\xfnm[ O.]}.
\newblock \bibinfo{title}{Small perturbation solutions for elliptic equations}.
\newblock \bibinfo{journal}{Comm Partial Differential Equations}
  \bibinfo{year}{2007};\bibinfo{volume}{32}(\bibinfo{number}{4-6}):\bibinfo{pages}{557--578}.
\newblock \URLprefix \url{https://doi.org/10.1080/03605300500394405}.
  \DOIprefix\doi{10.1080/03605300500394405}.
\bibitem[{Savin(2013)}]{MR2983006}
\bibinfo{author}{Savin\xfnm[ O.]}.
\newblock \bibinfo{title}{Pointwise {$C^{2,\alpha}$} estimates at the boundary
  for the {M}onge-{A}mp\`ere equation}.
\newblock \bibinfo{journal}{J Amer Math Soc}
  \bibinfo{year}{2013};\bibinfo{volume}{26}(\bibinfo{number}{1}):\bibinfo{pages}{63--99}.
\newblock \URLprefix \url{https://doi.org/10.1090/S0894-0347-2012-00747-4}.
  \DOIprefix\doi{10.1090/S0894-0347-2012-00747-4}.
\bibitem[{Silvestre and Sirakov(2014)}]{MR3246039}
\bibinfo{author}{Silvestre\xfnm[ L.]}, \bibinfo{author}{Sirakov\xfnm[ B.]}.
\newblock \bibinfo{title}{Boundary regularity for viscosity solutions of fully
  nonlinear elliptic equations}.
\newblock \bibinfo{journal}{Comm Partial Differential Equations}
  \bibinfo{year}{2014};\bibinfo{volume}{39}(\bibinfo{number}{9}):\bibinfo{pages}{1694--1717}.
\newblock \URLprefix \url{https://doi.org/10.1080/03605302.2013.842249}.
  \DOIprefix\doi{10.1080/03605302.2013.842249}.
\bibitem[{Sirakov(2010)}]{MR2592289}
\bibinfo{author}{Sirakov\xfnm[ B.]}.
\newblock \bibinfo{title}{Solvability of uniformly elliptic fully nonlinear
  {PDE}}.
\newblock \bibinfo{journal}{Arch Ration Mech Anal}
  \bibinfo{year}{2010};\bibinfo{volume}{195}(\bibinfo{number}{2}):\bibinfo{pages}{579--607}.
\newblock \URLprefix \url{https://doi.org/10.1007/s00205-009-0218-9}.
  \DOIprefix\doi{10.1007/s00205-009-0218-9}.
\bibitem[{Wang(1992{\natexlab{a}})}]{MR1135923}
\bibinfo{author}{Wang\xfnm[ L.]}.
\newblock \bibinfo{title}{On the regularity theory of fully nonlinear parabolic
  equations. {I}}.
\newblock \bibinfo{journal}{Comm Pure Appl Math}
  \bibinfo{year}{1992}{\natexlab{a}};\bibinfo{volume}{45}(\bibinfo{number}{1}):\bibinfo{pages}{27--76}.
\newblock \URLprefix \url{http://dx.doi.org/10.1002/cpa.3160450103}.
  \DOIprefix\doi{10.1002/cpa.3160450103}.
\bibitem[{Wang(1992{\natexlab{b}})}]{MR1139064}
\bibinfo{author}{Wang\xfnm[ L.]}.
\newblock \bibinfo{title}{On the regularity theory of fully nonlinear parabolic
  equations. {II}}.
\newblock \bibinfo{journal}{Comm Pure Appl Math}
  \bibinfo{year}{1992}{\natexlab{b}};\bibinfo{volume}{45}(\bibinfo{number}{2}):\bibinfo{pages}{141--178}.
\newblock \URLprefix \url{http://dx.doi.org/10.1002/cpa.3160450202}.
  \DOIprefix\doi{10.1002/cpa.3160450202}.
\bibitem[{Wang(1992{\natexlab{c}})}]{MR1151267}
\bibinfo{author}{Wang\xfnm[ L.]}.
\newblock \bibinfo{title}{On the regularity theory of fully nonlinear parabolic
  equations. {III}}.
\newblock \bibinfo{journal}{Comm Pure Appl Math}
  \bibinfo{year}{1992}{\natexlab{c}};\bibinfo{volume}{45}(\bibinfo{number}{3}):\bibinfo{pages}{255--262}.
\newblock \URLprefix \url{http://dx.doi.org/10.1002/cpa.3160450302}.
  \DOIprefix\doi{10.1002/cpa.3160450302}.
\bibitem[{Wu et~al.(2021)Wu, Lian and Zhang}]{MR_Israel}
\bibinfo{author}{Wu\xfnm[ D.]}, \bibinfo{author}{Lian\xfnm[ Y.]},
  \bibinfo{author}{Zhang\xfnm[ K.]}.
\newblock \bibinfo{title}{Pointwise boundary differentiability for fully
  nonlinear elliptic equations: arxiv:2101.00228}.
\newblock \bibinfo{year}{2021}.
\newblock \href{http://arxiv.org/abs/2101.00228}{\tt arXiv:2101.00228}.

\end{thebibliography}

\end{document}